\documentclass{article}%
\usepackage{amsmath}
\usepackage{amsfonts}
\usepackage{amssymb}
\usepackage{graphicx}
\usepackage{color}%
\providecommand{\U}[1]{\protect\rule{.1in}{.1in}}
\newtheorem{theorem}{Theorem}

\newtheorem{corollary}[theorem]{Corollary}

\newtheorem{definition}[theorem]{Definition}
\newtheorem{example}[theorem]{Example}

\newtheorem{lemma}[theorem]{Lemma}
\newtheorem{notation}[theorem]{Notation}

\newtheorem{remark}[theorem]{Remark}

\newenvironment{proof}[1][Proof]{\noindent\textbf{#1.} }{\ \rule{0.5em}{0.5em}}
\begin{document}

\title{Crouzeix-Raviart elements on simplicial meshes in $d$ dimensions}
\author{N.-E. Bohne\thanks{(nis-erik.bohne@math.uzh.ch), Institut f\"{u}r Mathematik,
Universit\"{a}t Z\"{u}rich, Winterthurerstr 190, CH-8057 Z\"{u}rich,
Switzerland}
\and Patrick Ciarlet, Jr.\thanks{(patrick.ciarlet@ensta-paris.fr), POEMS, CNRS,
INRIA, ENSTA Paris, Institut Polytechnique de Paris, 828 Bd des Mar\'{e}chaux,
91762 Palaiseau Cedex, France.}
\and S. Sauter\thanks{(stas@math.uzh.ch), Institut f\"{u}r Mathematik,
Universit\"{a}t Z\"{u}rich, Winterthurerstr 190, CH-8057 Z\"{u}rich,
Switzerland}}
\maketitle

\begin{abstract}
In this paper we introduce Crouzeix-Raviart elements of general polynomial
order $k$ and spatial dimension $d\geq2$ for simplicial finite element meshes.
We give explicit representations of the non-conforming basis functions and
prove that the conforming companion space, i.e., the conforming finite element
space of polynomial order $k$ is contained in the Crouzeix-Raviart space. We
prove a direct sum decomposition of the Crouzeix-Raviart space into (a
subspace of) the conforming companion space and the span of the non-conforming
basis functions.

Degrees of freedom are introduced which are bidual to the basis functions and
give rise to the definition of a local approximation/interpolation operator.
In two dimensions or for $k=1$, these freedoms can be split into simplex and
$\left(  d-1\right)  $ dimensional facet integrals in such a way that, in a
basis representation of Crouzeix-Raviart functions, all coefficients which
belong to basis functions related to lower-dimensional faces in the mesh are
determined by these facet integrals. It will also be shown that such a set of
degrees of freedom does \textbf{not} exist in higher space dimension and $k>1$.

\end{abstract}

\noindent\textbf{AMS-Classification}: Primary 33C45, 33C50, 65N30; Secondary
65N12.\medskip


\section{Introduction}

The Crouzeix-Raviart (CR) finite elements spaces have been introduced in the
seminal paper \cite{CrouzeixRaviart} and allow for a very economic
non-conforming discretization of the Stokes equation. The original definition
in \cite{CrouzeixRaviart} is general but implicit by imposing certain moment
conditions across interelement facets. In the original paper, basis functions
are presented for the lowest order CR elements in dimension $d\in\left\{
2,3\right\}  $; this is the reason that most implementations and
methodological developments of CR elements are restricted to lowest order.
Explicit representations for $d=2$ of the non-conforming CR functions exist in
the literature, see \cite{Ainsworth_Rankin}, \cite{Baran_Stoyan},
\cite{CCSS_CR_1}, \cite{BaranCVD}, \cite[for $p=4,6.$]{ChaLeeLee},
\cite{ccss_2012} while for spatial dimension $d\geq3$ the lowest order CR
element has been introduced in \cite{CrouzeixRaviart}. For \textit{quadratic}
CR elements in 3D, an explicit basis has been introduced in \cite{Fortin_d3}.
In \cite{CDS}, a spanning set of functions is presented in 3D for a
\textit{maximal} CR space of any polynomial order $k\in\mathbb{N}$ which
allows for a \textit{local basis}. However, the question of linear
independence is subtle, in particular, the definition of a basis for a minimal
CR space. In \cite{SauterTorres_CR3D} basis functions of a CR space have been
introduced for general space dimension $d\geq2$ and polynomial order $k$ which
are the analoga of the 2D CR spaces. However, the corresponding space do not
always include its conforming \textit{companion space}, i.e., the conforming
finite element space of degree $k$ and suboptimal convergence rates must be expected.

In this paper we introduce a basis for CR elements for general polynomial
degree $k$ and space dimension $d$ which contains the conforming companion
space and its approximation properties are inherited. We introduce degrees of
freedom which are bidual to the CR basis functions and allow to define local
approximation operators.

The paper is structured as follows. In Section \ref{SecCanCR}, we introduce
the non-conforming CR space in $d$ dimension and general polynomial order $k$.
We define the non-conforming CR functions in such a way that the conforming
companion space is contained in the CR space. For the construction of a full
set of basis function we prove in Section \ref{DirSumDecomp} a direct sum
decomposition of the CR space into (a subset of) the conforming companions
space and the span of the linearly independent CR functions. This construction
is given in Section \ref{SecBasisFull} where a full basis of the CR space is
defined which consists of the product of orthogonal polynomials with bubble
functions on simplicial entities (vertices, edges, facets, etc.). This choice
of basis functions has impact on the construction of degrees of freedom and an
interpolation operator. In Section \ref{SecDofsGenD} we introduce degrees of
freedom which are bidual to the CR basis functions and (mostly) defined by
integrals over \textit{simplices} for certain weight functions. In\ Section
\ref{SecDofs} we construct for the two-dimensional case degrees of freedom for
the edge basis functions which are defined as integrals over the
\textit{edges} with appropriate weight functions. Interestingly, such a
construction is possible only in two dimensions or for $k=1$; the details and
a proof will be given in Section \ref{SecNonExist}.

\section{Crouzeix-Raviart Finite Elements\label{SecCanCR}}

In this section, we define Crouzeix-Raviart spaces in $d$ dimensions,
$d=2,3$,\ldots\ and polynomial degree $k=1,2,\ldots$.

Let $\Omega\subset\mathbb{R}^{d}$ be a bounded Lipschitz polytope. As usual
$L^{2}\left(  \Omega\right)  $ is the space of measurable square-integrable
functions with scalar product $\left(  u,v\right)  _{L^{2}\left(
\Omega\right)  }:=\int_{\Omega}uv$ and norm $\left\Vert u\right\Vert
_{L^{2}\left(  \Omega\right)  }:=\left(  u,u\right)  _{L^{2}\left(
\Omega\right)  }^{1/2}$. The Sobolev space $H^{1}\left(  \Omega\right)  $
contains all $L^{2}\left(  \Omega\right)  $ functions whose weak derivatives
(exist and) are square integrable. All function spaces are considered over the
field of real numbers.

Let $\mathbb{N}:=\left\{  1,2,3,\ldots\right\}  $ and $\mathbb{N}%
_{0}:=\mathbb{N}\cup\left\{  0\right\}  $. We employ the multi-index notation
in the usual way: for a vector $\mathbf{x}=\left(  x_{\ell}\right)  _{\ell
=1}^{m}\in\mathbb{R}_{\geq0}^{m}$ and $%
%TCIMACRO{\TeXButton{boldalphanew}{\boldsymbol{\alpha}}}%
%BeginExpansion
\boldsymbol{\alpha}%
%EndExpansion
\in\mathbb{N}_{0}^{m}$, we set%
\begin{align*}
\mathbf{x}^{%
%TCIMACRO{\TeXButton{boldalphanew}{\boldsymbol{\alpha}}}%
%BeginExpansion
\boldsymbol{\alpha}%
%EndExpansion
}  &  :=%
%TCIMACRO{\dprod _{\ell=1}^{m}}%
%BeginExpansion
{\displaystyle\prod_{\ell=1}^{m}}
%EndExpansion
x_{i}^{\alpha_{i}},\quad\left\vert \mathbf{x}\right\vert :=\sum_{\ell=1}%
^{m}x_{i},\quad%
%TCIMACRO{\TeXButton{boldalphanew}{\boldsymbol{\alpha}}}%
%BeginExpansion
\boldsymbol{\alpha}%
%EndExpansion
!:=%
%TCIMACRO{\dprod _{\ell=1}^{m}}%
%BeginExpansion
{\displaystyle\prod_{\ell=1}^{m}}
%EndExpansion
\alpha_{\ell}!,\\
\mathbb{N}_{\leq k}^{m}  &  :=\left\{
%TCIMACRO{\TeXButton{boldalphanew}{\boldsymbol{\alpha}}}%
%BeginExpansion
\boldsymbol{\alpha}%
%EndExpansion
\in\mathbb{N}_{0}^{m}\mid\left\vert
%TCIMACRO{\TeXButton{boldalphanew}{\boldsymbol{\alpha}}}%
%BeginExpansion
\boldsymbol{\alpha}%
%EndExpansion
\right\vert \leq k\right\}  .
\end{align*}

Let $\mathcal{T}$ be a conforming finite element mesh for $\Omega$ consisting
of closed simplices $K\in\mathcal{T}$. By $\mathcal{F}\left(  K\right)  $ we
denote the set of $\left(  d-1\right)  $-dimensional closed facets in
$\partial K$. The set of vertices is denoted by $\mathcal{V}\left(  K\right)
$ and we employ this notation also for lower dimensional simplicial entities
in $\partial K$, e.g., for $F\in\mathcal{F}\left(  K\right)  $, the set
$\mathcal{V}\left(  F\right)  $ is the set of vertices of $F$.

We denote by $\hat{K}$ the (closed) reference simplex with vertices
\begin{equation}%
\begin{array}
[c]{ll}%
\mathbf{\hat{z}}_{0}:=\mathbf{0}, & \text{for }1\leq i\leq d\text{,
}\mathbf{\hat{z}}_{i}:=\mathbf{e}_{i},
\end{array}
\label{Defzhat}%
\end{equation}
where $\mathbf{e}_{i}\in\mathbb{R}^{d}$ is the $i$-th canonical unit vector in
$\mathbb{R}^{d}$. Moreover, let $\mathcal{F}$ ($\mathcal{V}$, resp.) be the
set of all $\left(  d-1\right)  $-dimensional facets (vertices, resp.) in the
mesh and let%
\begin{equation}%
\begin{array}
[c]{ll}%
\mathcal{F}_{\partial\Omega}:=\left\{  F\in\mathcal{F}\mid F\subset
\partial\Omega\right\}  , & \mathcal{F}_{\Omega}:=\mathcal{F}\backslash
\mathcal{F}_{\partial\Omega},\\
\mathcal{V}_{\partial\Omega}:=\left\{  \mathbf{z}\in\mathcal{V}\mid
\mathbf{z}\in\partial\Omega\right\}  , & \mathcal{V}_{\Omega}:=\mathcal{V}%
\backslash\mathcal{V}_{\partial\Omega}.
\end{array}
\label{Defboundinfacets}%
\end{equation}
For $F\in\mathcal{F}$, $\mathbf{z}\in\mathcal{V}$, we define facet and nodal
patches by%
\begin{equation}%
\begin{array}
[c]{ll}%
\mathcal{T}_{\mathbf{z}}:=\left\{  K\in\mathcal{T}:\mathbf{z}\in K\right\}
, & \omega_{\mathbf{z}}:=\bigcup_{K\in\mathcal{T}_{\mathbf{z}}}K,\\
\mathcal{F}_{\mathbf{z}}:=\left\{  F\in\mathcal{F}:\mathbf{z}\in F\right\}
, & \\
\mathcal{T}_{F}:=\left\{  K\in\mathcal{T}:F\subset K\right\}  , & \omega
_{F}:=\bigcup_{K\in\mathcal{T}_{F}}K.
\end{array}
\label{DefPatches}%
\end{equation}

For a measurable subset $D\subset\mathbb{R}^{d}$, we denote by $\left\vert
D\right\vert _{d}$ the $d$-dimensional volume of $D$ and skip the index $d$ if
the dimension is clear from the context, e.g., $\left\vert K\right\vert $
denotes the $d$-dimensional volume of a simplex $K\in\mathcal{T}$, $\left\vert
F\right\vert $ the $\left(  d-1\right)  $-dimensional volume of a facet
$F\in\mathcal{F}$, etc. For a finite discrete set $\mathcal{I}$ its
cardinality is denoted by $\left\vert \mathcal{I}\right\vert $.

For a conforming simplicial mesh $\mathcal{T}$ of the domain $\Omega$, let
\begin{equation}
H^{1}\left(  \mathcal{T}\right)  :=\left\{  u\in L^{2}\left(  \Omega\right)
\mid\forall K\in\mathcal{T}:~\left.  u\right\vert _{\overset{\circ}{K}}\in
H^{1}\left(  \overset{\circ}{K}\right)  \right\}  . \label{DefH1Tau}%
\end{equation}

For $k\in\mathbb{N}_{0}$ and a simplex $K$, we denote by $\mathbb{P}_{k}(K)$
the space of polynomials of maximal degree $k$ on $K$ and set $\mathbb{P}%
_{-1}\left(  K\right)  :=\left\{  0\right\}  $. For a $\left(  d-1\right)  $
dimensional facet $F$, we denote by $\mathbb{P}_{k}\left(  F\right)  $ the
space of $\left(  d-1\right)  $-variate polynomials in the local variables of
$F$. We also we need the subspace%
\begin{equation}
\overset{\bullet}{\mathbb{P}_{k}}\left(  K\right)  :=\left\{  v\in
\mathbb{P}_{k}\left(  K\right)  \mid\forall\mathbf{y}\in\mathcal{V}\left(
K\right)  :v\left(  \mathbf{y}\right)  =0\right\}  \label{defPkstrich}%
\end{equation}
and employ this notation also for lower dimensional simplicial entities in
$\partial K$, e.g., for $F\in\mathcal{F}\left(  K\right)  $:%
\begin{equation}
\overset{\bullet}{\mathbb{P}_{k}}\left(  F\right)  :=\left\{  v\in
\mathbb{P}_{k}\left(  F\right)  \mid\forall\mathbf{y}\in\mathcal{V}\left(
F\right)  :v\left(  \mathbf{y}\right)  =0\right\}  . \label{PkFprime}%
\end{equation}

For the definition of the non-conforming Crouzeix-Raviart space, the functions
$B_{k}^{\operatorname*{CR},K}$, $K\in\mathcal{T}$, and $B_{k}%
^{\operatorname*{CR},F}$, $F\in\mathcal{F}_{\Omega}$, will play an important
role; their definition require some preliminaries.

\begin{notation}
\label{NotBary}For a simplex $K$ and $\mathbf{z}\in\mathcal{V}\left(
K\right)  $, the affine function $\lambda_{K,\mathbf{z}}\in\mathbb{P}%
_{1}\left(  K\right)  $ is the \emph{barycentric coordinate} for the vertex
$\mathbf{z}$ and characterized by%
\[
\lambda_{K,\mathbf{z}}\left(  \mathbf{y}\right)  =\delta_{\mathbf{y}%
,\mathbf{z}}\quad\forall\mathbf{y}\in\mathcal{V}\left(  K\right)  ,
\]
where $\delta_{\mathbf{y},\mathbf{z}}$ denotes the Kronecker delta.
\end{notation}

Let $\alpha,\beta>-1$ and $n\in\mathbb{N}_{0}$. The \emph{Jacobi polynomial}
$P_{n}^{\left(  \alpha,\beta\right)  }$ is a polynomial of degree $n$ such
that
\begin{equation}
\int_{-1}^{1}P_{n}^{\left(  \alpha,\beta\right)  }\left(  x\right)  \,q\left(
x\right)  \left(  1-x\right)  ^{\alpha}\left(  1+x\right)  ^{\beta}\,dx=0
\label{charJac}%
\end{equation}
for all polynomials $q$ of degree less than $n$, one has (cf. \cite[Table
18.6.1]{NIST:DLMF}):%
\begin{equation}
P_{n}^{\left(  \alpha,\beta\right)  }\left(  1\right)  =\frac{\left(
\alpha+1\right)  _{n}}{n!},\qquad P_{n}^{\left(  \alpha,\beta\right)  }\left(
-1\right)  =\left(  -1\right)  ^{n}\frac{\left(  \beta+1\right)  _{n}}{n!}.
\label{Pnormalization}%
\end{equation}
Here, the \emph{shifted factorial} is defined by $\left(  a\right)
_{n}:=a\left(  a+1\right)  \ldots\left(  a+n-1\right)  $ for $n>0$ and
$\left(  a\right)  _{0}:=1$. Note that $P_{n}^{\left(  0,0\right)  }$ are the
Legendre polynomials (see \cite[18.7.9]{NIST:DLMF}). In our application, the
Jacobi polynomials with $\alpha=0$ and $\beta=d-2$ are relevant. In this case,
the endpoint values are given by
%(cf. \cite[Table 18.6.1]{NIST:DLMF})%
\begin{equation}
P_{n}^{\left(  0,d-2\right)  }\left(  1\right)  =1\quad\text{and\quad}%
P_{n}^{\left(  0,d-2\right)  }\left(  -1\right)  =\left(  -1\right)  ^{n}%
\rho_{n}\text{ with }\rho_{n}:=\binom{n+d-2}{n}. \label{defrho}%
\end{equation}

Next, we define the non-conforming shape functions for the Crouzeix-Raviart
space. They differ from the definition in \cite[(3.8)]{SauterTorres_CR3D} and
lead to a Crouzeix-Raviart space which, in contrast to the space in
\cite{SauterTorres_CR3D}, always contains the conforming space as a subspace
(see Corollary \ref{CorCont}).

\begin{definition}
\label{DefNCShapeFct}Let $\mathcal{T}$ be a conforming simplicial finite
element mesh.

\begin{enumerate}
\item For any $K\in\mathcal{T}$, the non-conforming\emph{ simplex function}
$B_{k}^{\operatorname*{CR},K}\in\mathbb{P}_{k}\left(  \mathcal{T}\right)  $ is
given by%
\begin{equation}
B_{k}^{\operatorname*{CR},K}:=\left\{
\begin{array}
[c]{ll}%
\frac{1}{d}\left(  -1+%
%TCIMACRO{\dsum \limits_{\mathbf{z}\in\mathcal{V}\left(  K\right)  }}%
%BeginExpansion
{\displaystyle\sum\limits_{\mathbf{z}\in\mathcal{V}\left(  K\right)  }}
%EndExpansion
P_{k}^{\left(  0,d-2\right)  }\left(  1-2\lambda_{K,\mathbf{z}}\right)
\right)  & \text{on }K\text{,}\\
0 & \text{otherwise.}%
\end{array}
\right.  \label{def:CRtetrahedron}%
\end{equation}
The space spanned by these functions is%
\[
V_{k}^{\operatorname*{nc}}\left(  \mathcal{T}\right)  :=\operatorname*{span}%
\left\{  B_{k}^{\operatorname*{CR},K}:K\in\mathcal{T}\right\}  .
\]

\item For any $F\in\mathcal{F}$, the non-conforming \emph{facet function}
$B_{k}^{\operatorname*{CR},F}\in\mathbb{P}_{k}\left(  \mathcal{T}\right)  $ is
given by%
\begin{equation}
B_{k}^{\operatorname*{CR},F}:=\left\{
\begin{array}
[c]{ll}%
P_{k}^{\left(  0,d-2\right)  }\left(  1-2\lambda_{K,F}\right)  -B_{k}%
^{\operatorname*{CR},K} & \text{for }K\in\mathcal{T}_{F},\\
0 & \text{otherwise,}%
\end{array}
\right.  \label{def:CRfacet}%
\end{equation}
where $\lambda_{K,F}$ denotes the barycentric coordinate for a simplex $K$
adjacent to $F$ corresponding to the vertex $\mathbf{z}\in\mathcal{V}\left(
K\right)  $ opposite to $F$. The spaces spanned by these functions are%
\begin{align*}
V_{k}^{\operatorname*{nc}}\left(  \mathcal{F}\right)   &
:=\operatorname*{span}\left\{  B_{k}^{\operatorname*{CR},F}:F\in
\mathcal{F}\right\}  ,\\
V_{k}^{\operatorname*{nc}}\left(  \mathcal{F}_{\Omega}\right)   &
:=\operatorname*{span}\left\{  B_{k}^{\operatorname*{CR},F}:F\in
\mathcal{F}_{\Omega}\right\}  .
\end{align*}

\end{enumerate}
\end{definition}

\begin{remark}
\label{Remk=1}Note that $\operatorname*{supp}B_{k}^{\operatorname*{CR}%
,F}=\omega_{F}$ and $\operatorname*{supp}B_{k}^{\operatorname*{CR},K}=K$.
Since $\lambda_{K,F}=0$ on $F$, the value of $\left.  P_{k}^{\left(
0,d-2\right)  }\left(  1-2\lambda_{K,F}\right)  \right\vert _{F}$ in
(\ref{def:CRfacet}) is equal to $P_{k}^{\left(  0,d-2\right)  }\left(
1\right)  =1$.

For even $k$, the span $V_{k}^{\operatorname*{nc}}\left(  \mathcal{T}\right)
$ of the non-conforming shape functions $B_{k}^{\operatorname*{CR},K}$ is the
same as the one in \cite[before Def. 3.6]{SauterTorres_CR3D}. For $k=1$, the
space $V_{1}^{\operatorname*{nc}}\left(  \mathcal{F}_{\Omega}\right)  $ equals
the one in \cite[before Def. 3.6]{SauterTorres_CR3D}. The reason is that
$B_{1}^{\operatorname*{CR},K}=0$ which will be proved as Case 2 of Lemma
\ref{LemBasisPartial}.

For $d=2$ and odd $k$ the subtraction of $B_{k}^{\operatorname*{CR},K}$ in
(\ref{def:CRfacet}) can be removed since $B_{k}^{\operatorname*{CR},K}$ is
continuous in this case. On a facet $F\subset\partial K$ it holds%
\begin{equation}
\left.  \left(  P_{k}^{\left(  0,0\right)  }\left(  1-2\lambda_{K,F}\right)
-B_{k}^{\operatorname*{CR},K}\right)  \right\vert _{F}=\left.  \left(
P_{k}^{\left(  0,0\right)  }\left(  1-2\lambda_{K,F}\right)  \right)
\right\vert _{F}. \label{bkbound}%
\end{equation}

\end{remark}

Next, we introduce the relevant finite element spaces. For $k\in\mathbb{N}$,
let $S_{k}\left(  \mathcal{T}\right)  $ denote the space of globally
continuous, piecewise polynomials of maximal degree $k$:%
\[
S_{k}\left(  \mathcal{T}\right)  :=\left\{  v\in C^{0}(\overline{\Omega}%
)\mid\left.  v\right\vert _{K}\in\mathbb{P}_{k}(K)\quad\forall K\in
\mathcal{T}\right\}
\]
and the subspace with zero traces on the boundary:%
\[
S_{k,0}\left(  \mathcal{T}\right)  :=S_{k}\left(  \mathcal{T}\right)  \cap
H_{0}^{1}\left(  \Omega\right)  .
\]

The space of discontinuous polynomials of maximal degree $k$ is%
\[
\mathbb{P}_{k}\left(  \mathcal{T}\right)  :=\left\{  p\in L^{2}\left(
\Omega\right)  \mid\left.  p\right\vert _{\overset{\circ}{K}}\in\mathbb{P}%
_{k}(\overset{\circ}{K})\quad\forall K\in\mathcal{T}\right\}  ,
\]
where we employed the following notation.

\begin{notation}
A function $v\in H^{1}\left(  \mathcal{T}\right)  $ may be discontinuous
across simplex facets. In this way the values on a facet are not uniquely
defined. For $K\in\mathcal{T}$ the function $\left.  v\right\vert
_{\overset{\circ}{K}}$ is well defined and the function $\left.  v\right\vert
_{K}$ equals $\left.  v\right\vert _{\overset{\circ}{K}}$ in the interior of
$K$ and is defined on $\partial K$ as the continuous continuation of $\left.
v\right\vert _{\overset{\circ}{K}}$ to $\partial K$.

In a similar fashion the sum of functions $v_{1},v_{2}\in H^{1}\left(
\mathcal{T}\right)  $ on some $K\in\mathcal{T}$ is defined by%
\[
\left.  \left(  v_{1}+v_{2}\right)  \right\vert _{K}:=\left.  v_{1}\right\vert
_{K}+\left.  v_{2}\right\vert _{K}.
\]
If the functions $v,v_{1},v_{2}$ belong to $C^{0}\left(  \overline{\Omega
}\right)  $ these definitions coincide with the standard restriction to $K$.
\end{notation}

Next we define the non-conforming Crouzeix-Raviart space in the original
implicit way (see \cite{CrouzeixRaviart}). For a function $v\in H^{1}\left(
\mathcal{T}\right)  $, we denote by $\left[  v\right]  _{F}$ the jump across
the facet $F\in\mathcal{F}$. The Crouzeix-Raviart space of order $k$ is
defined implicitly by%
\[
\operatorname*{CR}\nolimits_{k}^{\max}\left(  \mathcal{T}\right)  :=\left\{
v\in\mathbb{P}_{k}\left(  \mathcal{T}\right)  \mid\forall F\in\mathcal{F}%
_{\Omega}\quad\left[  v\right]  _{F}\perp\mathbb{P}_{k-1}\left(  F\right)
\right\}
\]
and the subspace with (approximate) zero boundary conditions by%
\begin{equation}
\operatorname*{CR}\nolimits_{k,0}^{\max}\left(  \mathcal{T}\right)  :=\left\{
v\in\operatorname*{CR}\nolimits_{k}^{\max}\left(  \mathcal{T}\right)
\mid\forall F\in\mathcal{F}_{\partial\Omega}\quad v\perp\mathbb{P}%
_{k-1}\left(  F\right)  \right\}  . \label{DefCRmax}%
\end{equation}
Here the symbol \textquotedblleft$\perp$\textquotedblright\ always refers to
$L^{2}$ orthogonality.

\begin{lemma}
\label{LemCRprop}The functions $B_{k}^{\operatorname*{CR},F}$ and
$B_{k}^{\operatorname*{CR},K}$ in Definition \ref{DefNCShapeFct} belong to
$\operatorname*{CR}\nolimits_{k}^{\max}\left(  \mathcal{T}\right)  $.
\end{lemma}

%

%TCIMACRO{\TeXButton{Proof}{\proof}}%
%BeginExpansion
\proof
%EndExpansion
This follows from \cite[Theorem 3.5]{SauterTorres_CR3D} by showing that
$P_{k}^{\left(  0,d-2\right)  }$ satisfies \cite[Condition (A.1)]%
{SauterTorres_CR3D}:%
%TCIMACRO{\TeXButton{orthofacetprop}{\begin{subequations}
%\label{orthofacetprop}
%\end{subequations}}}%
%BeginExpansion
\begin{subequations}
\label{orthofacetprop}
\end{subequations}%
%EndExpansion%
\begin{align}
\left.  P_{k}^{\left(  0,d-2\right)  }\left(  1-2\lambda_{K,\mathbf{z}%
}\right)  \right\vert _{F_{\mathbf{z}}}  &  =1\tag{%
%TCIMACRO{\TeXButton{orthofacetprop}{\ref{orthofacetprop}}}%
%BeginExpansion
\ref{orthofacetprop}%
%EndExpansion
a}\label{orthofacetpropa}\\
\forall F\in\mathcal{F}\left(  K\right)  \backslash\left\{  F_{\mathbf{z}%
}\right\}  \quad\int_{F}P_{k}^{\left(  0,d-2\right)  }\left(  1-2\lambda
_{K,\mathbf{z}}\right)  q  &  =0\quad\forall q\in\mathbb{P}_{k-1}\left(
F\right)  \tag{%
%TCIMACRO{\TeXButton{orthofacetprop}{\ref{orthofacetprop}}}%
%BeginExpansion
\ref{orthofacetprop}%
%EndExpansion
b}\label{orthofacetpropb}%
\end{align}
for all $K\in\mathcal{T}$, $\mathbf{z}\in\mathcal{V}\left(  K\right)  $, where
$F_{\mathbf{z}}$ denotes the facet of $K$ opposite to $\mathbf{z}$.

The first condition follows from $\left.  \lambda_{K,\mathbf{z}}\right\vert
_{F_{\mathbf{z}}}=0$ and $P_{k}^{\left(  0,d-2\right)  }\left(  1\right)  =1$
(cf. (\ref{Pnormalization})).

The second condition follows from \cite[(2.5.21)]{Yuan_inproc} by choosing $%
%TCIMACRO{\TeXButton{boldalphanew}{\boldsymbol{\alpha}}}%
%BeginExpansion
\boldsymbol{\alpha}%
%EndExpansion
=\left(  k,0,\ldots,0\right)  $, $\kappa_{1}=\ldots=\kappa_{d+1}=1/2$ therein
and by using
\[
P_{k}^{\left(  d-2,0\right)  }\left(  2\lambda_{K,\mathbf{z}}-1\right)
=\left(  -1\right)  ^{k}P_{k}^{\left(  0,d-2\right)  }\left(  1-2\lambda
_{K,\mathbf{z}}\right)
\]
(see \cite[Table 18.6.1]{NIST:DLMF}). Here, we present an elementary proof to
give an intuition of this specific choice of the Jacobi weights.

By using the affine equivalence of $K$ to the reference simplex it is
sufficient to show the second condition for $\widehat{K}$, $\widehat
{F}:=\left\{  \mathbf{x}=\left(  x_{j}\right)  _{j=1}^{d}\in\widehat{K}\mid
x_{d}=0\right\}  $, and $\mathbf{z}=\left(  1,0,\ldots,0\right)  ^{T}$: for
all $%
%TCIMACRO{\TeXButton{boldalphanew}{\boldsymbol{\alpha}}}%
%BeginExpansion
\boldsymbol{\alpha}%
%EndExpansion
=\left(  \alpha_{i}\right)  _{i=1}^{d-1}\in\mathbb{N}_{\leq k-1}^{d-1}$ and
corresponding monomials%
\[
g\left(  \mathbf{x}\right)  :=%
%TCIMACRO{\dprod _{j=1}^{d-1}}%
%BeginExpansion
{\displaystyle\prod_{j=1}^{d-1}}
%EndExpansion
x_{j}^{\alpha_{j}}\quad\forall\mathbf{x}=\left(  x_{j}\right)  _{j=1}^{d-1}%
\]
it holds%
\[
%0=
\int_{\widehat{F}}P_{k}^{\left(  0,d-2\right)  }\left(  1-2x_{1}\right)
g\left(  \mathbf{x}\right)  d\mathbf{x}=\int_{0}^{1}P_{k}^{\left(
0,d-2\right)  }\left(  1-2x_{1}\right)  G\left(  x_{1}\right)  dx_{1},
\]
where%
\[
G\left(  x_{1}\right)  =\int_{0}^{1-x_{1}}\int_{0}^{1-x_{1}-x_{2}}\ldots
\int_{0}^{1-x_{1}-\ldots-x_{d-2}}g\left(  \mathbf{x}\right)  dx_{d-1}%
dx_{d-2}\ldots dx_{2}.
\]
Let $%
%TCIMACRO{\TeXButton{boldalphanew}{\boldsymbol{\alpha}}}%
%BeginExpansion
\boldsymbol{\alpha}%
%EndExpansion
^{\prime}=\left(  \alpha_{j}\right)  _{j=2}^{d-1}$. From the proof of
\cite[Lemma A.1]{SauterTorres_CR3D} it follows%
\[
G\left(  x_{1}\right)  =c_{d,%
%TCIMACRO{\TeXButton{boldalphanew}{\boldsymbol{\alpha}}}%
%BeginExpansion
\boldsymbol{\alpha}%
%EndExpansion
}x_{1}^{\alpha_{1}}\left(  1-x_{1}\right)  ^{d-2+\left\vert
%TCIMACRO{\TeXButton{boldalphanew}{\boldsymbol{\alpha}}}%
%BeginExpansion
\boldsymbol{\alpha}%
%EndExpansion
^{\prime}\right\vert }\quad\text{for }c_{d,%
%TCIMACRO{\TeXButton{boldalphanew}{\boldsymbol{\alpha}}}%
%BeginExpansion
\boldsymbol{\alpha}%
%EndExpansion
}:=\frac{%
%TCIMACRO{\TeXButton{boldalphanew}{\boldsymbol{\alpha}}}%
%BeginExpansion
\boldsymbol{\alpha}%
%EndExpansion
^{\prime}!}{\left(  d-2+\left\vert
%TCIMACRO{\TeXButton{boldalphanew}{\boldsymbol{\alpha}}}%
%BeginExpansion
\boldsymbol{\alpha}%
%EndExpansion
^{\prime}\right\vert \right)  !}.
\]
Hence,%
\begin{align*}
\int_{\widehat{F}}P_{k}^{\left(  0,d-2\right)  }\left(  1-2x_{1}\right)
g\left(  \mathbf{x}\right)  d\mathbf{x}  &  =c_{d,%
%TCIMACRO{\TeXButton{boldalphanew}{\boldsymbol{\alpha}}}%
%BeginExpansion
\boldsymbol{\alpha}%
%EndExpansion
}\int_{0}^{1}P_{k}^{\left(  0,d-2\right)  }\left(  1-2x_{1}\right)
x_{1}^{\alpha_{1}}\left(  1-x_{1}\right)  ^{d-2+\left\vert
%TCIMACRO{\TeXButton{boldalphanew}{\boldsymbol{\alpha}}}%
%BeginExpansion
\boldsymbol{\alpha}%
%EndExpansion
^{\prime}\right\vert }dx_{1}\\
&  =\frac{c_{d,%
%TCIMACRO{\TeXButton{boldalphanew}{\boldsymbol{\alpha}}}%
%BeginExpansion
\boldsymbol{\alpha}%
%EndExpansion
}}{2^{d-1}}\int_{-1}^{1}\omega\left(  t\right)  P_{k}^{\left(  0,d-2\right)
}\left(  t\right)  \left(  \frac{1-t}{2}\right)  ^{\alpha_{1}}\left(
\frac{t+1}{2}\right)  ^{\left\vert
%TCIMACRO{\TeXButton{boldalphanew}{\boldsymbol{\alpha}}}%
%BeginExpansion
\boldsymbol{\alpha}%
%EndExpansion
^{\prime}\right\vert }dt
\end{align*}
for the weight function $\omega\left(  t\right)  =\left(  t+1\right)  ^{d-2}$.
Since $\left(  \frac{1-t}{2}\right)  ^{\alpha_{1}}\left(  \frac{t+1}%
{2}\right)  ^{\left\vert
%TCIMACRO{\TeXButton{boldalphanew}{\boldsymbol{\alpha}}}%
%BeginExpansion
\boldsymbol{\alpha}%
%EndExpansion
^{\prime}\right\vert }\in\mathbb{P}_{\left\vert
%TCIMACRO{\TeXButton{boldalphanew}{\boldsymbol{\alpha}}}%
%BeginExpansion
\boldsymbol{\alpha}%
%EndExpansion
\right\vert }\subset\mathbb{P}_{k-1}$, the orthogonality properties of the
Jacobi polynomials imply%
\[
\int_{\widehat{F}}P_{k}^{\left(  0,d-2\right)  }\left(  1-2x_{1}\right)
g\left(  \mathbf{x}\right)  d\mathbf{x}=0.
\]%
%TCIMACRO{\TeXButton{End Proof}{\endproof}}%
%BeginExpansion
\endproof
%EndExpansion

Next, we introduce canonical Crouzeix-Raviart functions on simplicial meshes
in $d$ dimensions by modifying the basic space in \cite{SauterTorres_CR3D}.

\begin{definition}
\label{DefCanCR} The scalar \emph{canonical Crouzeix-Raviart space} of order
$k$ for conforming simplicial finite element meshes $\mathcal{T}$ is given by%
\[
\operatorname*{CR}\nolimits_{k}\left(  \mathcal{T}\right)  :=\left\{
\begin{array}
[c]{cc}%
S_{k}\left(  \mathcal{T}\right)  +V_{k}^{\operatorname*{nc}}\left(
\mathcal{T}\right)  & \text{if }k\text{ is even,}\\
S_{k}\left(  \mathcal{T}\right)  +V_{k}^{\operatorname*{nc}}\left(
\mathcal{F}\right)  & \text{if }k\text{ is odd.}%
\end{array}
\right.
\]
and the subspace for approximating problems with zero boundary conditions by:%
\begin{equation}
\operatorname*{CR}\nolimits_{k,0}\left(  \mathcal{T}\right)  :=\left\{
\begin{array}
[c]{cc}%
S_{k,0}\left(  \mathcal{T}\right)  +V_{k}^{\operatorname*{nc}}\left(
\mathcal{T}\right)  & \text{if }k\text{ is even,}\\
S_{k,0}\left(  \mathcal{T}\right)  +V_{k}^{\operatorname*{nc}}\left(
\mathcal{F}_{\Omega}\right)  & \text{if }k\text{ is odd.}%
\end{array}
\right.  \label{defCRcan}%
\end{equation}

\end{definition}

We emphasize that the sums in (\ref{defCRcan}), in general, are not direct.

\section{A direct sum decomposition of the canonical Crouzeix-Raviart
space\label{DirSumDecomp}}

The Lagrange interpolation points for the reference element are given by%
\begin{equation}
\mathcal{N}_{k}\left(  \hat{K}\right)  :=\left\{
\begin{array}
[c]{ll}%
\frac{1}{d+1}\mathbf{1} & k=0,\\
\left\{  \dfrac{%
%TCIMACRO{\TeXButton{boldmue}{\boldsymbol{\mu}}}%
%BeginExpansion
\boldsymbol{\mu}%
%EndExpansion
}{k}:%
%TCIMACRO{\TeXButton{boldmue}{\boldsymbol{\mu}}}%
%BeginExpansion
\boldsymbol{\mu}%
%EndExpansion
\in\mathbb{N}_{k}^{d}\right\}  & k\geq1,
\end{array}
\right.  \label{defNkhut}%
\end{equation}
where $\mathbf{1}\in\mathbb{R}^{d}$ is the vector with constant coefficients
$1$. As usual they are lifted to $K\in\mathcal{T}$ via an affine map
$\varphi_{K}:\hat{K}\rightarrow K$ to obtain the interpolation points on the
mesh element $K$ via $\mathcal{N}_{k}\left(  K\right)  :=\left\{  \varphi
_{K}\left(  \mathbf{\hat{y}}\right)  :\mathbf{\hat{y}}\in\mathcal{N}%
_{k}\left(  \hat{K}\right)  \right\}  $. We also will need the subset
\begin{equation}
\overset{\bullet}{\mathcal{N}_{k}}\left(  K\right)  :=\mathcal{N}_{k}\left(
K\right)  \backslash\mathcal{V}\left(  K\right)  \label{NlkprimeK}%
\end{equation}
and the corresponding subspaces of $S_{k}\left(  \mathcal{T}\right)  $,
$S_{k,0}\left(  \mathcal{T}\right)  $ consisting of functions which vanish at
the vertices of the mesh are given by%

\begin{align}
\overset{\bullet}{S_{k}}\left(  \mathcal{T}\right)   &  :=\left\{  v\in
S_{k}\left(  \mathcal{T}\right)  \mid\forall\mathbf{z}\in\mathcal{V}:\quad
v\left(  \mathbf{z}\right)  =0\right\}  ,\label{DefSkprime0}\\
\overset{\bullet}{S_{k,0}}\left(  \mathcal{T}\right)   &  :=\left\{  v\in
S_{k,0}\left(  \mathcal{T}\right)  \mid\forall\mathbf{z}\in\mathcal{V}:\quad
v\left(  \mathbf{z}\right)  =0\right\}  . \label{DefSkprime}%
\end{align}

The Lagrange interpolation points for $\mathcal{T}$ are given by
$\mathcal{N}_{k}\left(  \mathcal{T}\right)  :=%
%TCIMACRO{\dbigcup \limits_{K\in\mathcal{T}}}%
%BeginExpansion
{\displaystyle\bigcup\limits_{K\in\mathcal{T}}}
%EndExpansion
\mathcal{N}_{k}\left(  K\right)  $. Let $\mathcal{N}_{k,0}\left(
\mathcal{T}\right)  \subset\mathcal{N}_{k}\left(  \mathcal{T}\right)  $ be the
subset, where the Lagrange points on the boundary $\partial\Omega$ are
removed. Further, let $\overset{\bullet}{\mathcal{N}_{k}}\left(
\mathcal{T}\right)  :=\mathcal{N}_{k}\left(  \mathcal{T}\right)
\backslash\mathcal{V}$ and $\overset{\bullet}{\mathcal{N}_{k,0}}\left(
\mathcal{T}\right)  :=\mathcal{N}_{k,0}\left(  \mathcal{T}\right)
\backslash\mathcal{V}_{\Omega}$. As usual, for $\mathbf{z}\in\mathcal{N}%
_{k}\left(  \mathcal{T}\right)  $, the corresponding Lagrange basis
$B_{k}^{\mathbf{z}}\in S_{k}\left(  \mathcal{T}\right)  $ is defined
implicitly by the conditions%
\[
B_{k}^{\mathbf{z}}\left(  \mathbf{y}\right)  =\delta_{\mathbf{z},\mathbf{y}%
}\quad\forall\mathbf{y},\mathbf{z}\in\mathcal{N}_{k}\left(  \mathcal{T}%
\right)  .
\]
The basis functions on the reference element $\hat{K}$ are denoted by $\hat
{B}_{k}^{\mathbf{\hat{z}}}\in\mathbb{P}_{k}\left(  \hat{K}\right)  $ for
$\mathbf{\hat{z}}\in\mathcal{N}_{k}\left(  \hat{K}\right)  $, and
characterized by: $\hat{B}_{k}^{\mathbf{\hat{z}}}\left(  \mathbf{\hat{y}%
}\right)  =\delta_{\mathbf{\hat{y}},\mathbf{\hat{z}}}$ for all $\mathbf{\hat
{y}},\mathbf{\hat{z}}\in\mathcal{N}_{k}\left(  \hat{K}\right)  $.

Next, some useful properties of the non-conforming functions $B_{k}%
^{\operatorname*{CR},K}$ and $B_{k}^{\operatorname*{CR},F}$ and corresponding
spaces $V_{k}^{\operatorname*{nc}}\left(  \mathcal{T}\right)  $ $V_{k}%
^{\operatorname*{nc}}\left(  \mathcal{F}\right)  $, and $V_{k}%
^{\operatorname*{nc}}\left(  \mathcal{F}_{\Omega}\right)  $ are derived.

The endpoint properties of the Jacobi polynomials (cf. (\ref{defrho})) allow
us to compute the values of the non-conforming functions at the simplex
vertices.
%TCIMACRO{\TeXButton{BkCRvertex}{\begin{subequations}
%\label{BkCRvertex}
%\end{subequations}}}%
%BeginExpansion
\begin{subequations}
\label{BkCRvertex}
\end{subequations}%
%EndExpansion
For $K\in\mathcal{T}$ it holds:%
\begin{equation}
\left.  B_{k}^{\operatorname*{CR},K}\right\vert _{K}\left(  \mathbf{y}\right)
=\frac{d-1+\left(  -1\right)  ^{k}\rho_{k}}{d}\quad\forall\mathbf{y}%
\in\mathcal{V}\left(  K\right)  . \tag{%
%TCIMACRO{\TeXButton{BkCRvertex}{\ref{BkCRvertex}}}%
%BeginExpansion
\ref{BkCRvertex}%
%EndExpansion
a}\label{BkCRvertexa}%
\end{equation}
For $F\in\mathcal{F}$, $K\in\mathcal{T}_{F}$ it holds:%
\begin{equation}
\left.  B_{k}^{\operatorname*{CR},F}\right\vert _{K}\left(  \mathbf{y}\right)
=\left\{
\begin{array}
[c]{ll}%
\frac{1+\left(  -1\right)  ^{k+1}\rho_{k}}{d} & \text{if }\mathbf{y}%
\in\mathcal{V}\left(  F\right)  ,\\
\frac{1+\left(  -1\right)  ^{k+1}\rho_{k}}{d}\left(  1-d\right)  & \text{if
}\mathbf{y\in}\mathcal{V}\left(  K\right)  \text{ is opposite to }F.
\end{array}
\right.  \tag{%
%TCIMACRO{\TeXButton{BkCRvertex}{\ref{BkCRvertex}}}%
%BeginExpansion
\ref{BkCRvertex}%
%EndExpansion
b}\label{BkCRvertexb}%
\end{equation}
We need the following technical lemma with regard to the vertex values in
(\ref{BkCRvertex}).

\begin{lemma}
\label{LemmaTechnical}Let $d\geq2$ and $k\geq1$. It holds%
%TCIMACRO{\TeXButton{TechThings}{\begin{subequations}
%\label{TechThings}
%\end{subequations}}}%
%BeginExpansion
\begin{subequations}
\label{TechThings}
\end{subequations}%
%EndExpansion%
\begin{align}
d-1+\left(  -1\right)  ^{k}\rho_{k}  &  =0\iff\left(  k=1\right)  \vee\left(
d=2\wedge k\text{ is odd}\right)  ,\tag{%
%TCIMACRO{\TeXButton{TechThings}{\ref{TechThings}}}%
%BeginExpansion
\ref{TechThings}%
%EndExpansion
a}\label{TechThingsa}\\
d+\left(  -1\right)  ^{k}\rho_{k}  &  \neq0,\tag{%
%TCIMACRO{\TeXButton{TechThings}{\ref{TechThings}}}%
%BeginExpansion
\ref{TechThings}%
%EndExpansion
b}\label{TechThingsb}\\
\rho_{k}  &  =1\iff d=2. \tag{%
%TCIMACRO{\TeXButton{TechThings}{\ref{TechThings}}}%
%BeginExpansion
\ref{TechThings}%
%EndExpansion
c}\label{TechThingsc}%
\end{align}

\end{lemma}%

%TCIMACRO{\TeXButton{Proof}{\proof}}%
%BeginExpansion
\proof
%EndExpansion
We first prove (\ref{TechThingsa}).

The direction \textquotedblleft$\impliedby$\textquotedblright\ in
(\ref{TechThingsa})\ follows by direct computation since $\rho_{k}=d-1$ for
$k=1$ and $\rho_{k}=1$ for $d=2$. In the right-hand side of (\ref{TechThingsa}%
) $k$ is odd so that $\theta\left(  k\right)  :=d-1+\left(  -1\right)
^{k}\rho_{k}=0$.

Next we prove \textquotedblleft$\implies$\textquotedblright\ in
(\ref{TechThingsa}). For even $k$, we have $\theta\left(  k\right)
=d-1+\rho_{k}\geq d>0$. Hence $\theta\left(  k\right)  =0$ implies $k$ is odd
and $\theta\left(  k\right)  =d-1-\rho_{k}$.

For $k=1$, we have $\rho_{k}=d-1$ so that $\theta\left(  1\right)  =0$ and it
remains to consider odd $k\geq3$. It is easy to verify that $\binom{k+d-2}{k}$
is monotonically growing with respect to $k$. For odd $k\geq3$ it holds%
\[
\theta\left(  k\right)  =d-1-\binom{k+d-2}{k}\leq d-1-\binom{d+1}{3}=-\frac
{1}{6}\left(  d+3\right)  \left(  d-1\right)  \left(  d-2\right)  .
\]
Since $d\geq2$ the condition $\theta\left(  k\right)  =0$ implies $d=2$.

Next we prove (\ref{TechThingsb}). For even $k$ it holds $d+\left(  -1\right)
^{k}\rho_{k}=d+\rho_{k}>0$. For $d=2$ it holds $d+\left(  -1\right)  ^{k}%
\rho_{k}=2+\left(  -1\right)  ^{k}>0$ while for $k=1$ we get $d+\left(
-1\right)  ^{k}\rho_{k}=1$.

It remains to consider the case $d\geq3$ and odd $k\geq3$. A similar
monotonicity argument leads to%
\[
d-\rho_{k}=d-\binom{k+d-2}{k}\overset{k\geq3}{\leq}d-\binom{d+1}{3}=-\frac
{1}{6}d\left(  d^{2}-7\right)  <0.
\]
Hence, (\ref{TechThingsa}) is proved. The equivalence (\ref{TechThingsc}) is
trivial.%
%TCIMACRO{\TeXButton{End Proof}{\endproof}}%
%BeginExpansion
\endproof
%EndExpansion

\begin{lemma}
\label{LemDefPsikz}Let $k\in\mathbb{N}$. For $\mathbf{z}\in\mathcal{V}$, the
function%
\begin{equation}
\psi_{k}^{\mathbf{z}}:=\sum_{F\in\mathcal{F}_{\mathbf{z}}}B_{k}%
^{\operatorname*{CR},F} \label{psikz}%
\end{equation}
satisfies%
\begin{equation}%
\begin{array}
[c]{l}%
\psi_{k}^{\mathbf{z}}\in\left\{
\begin{array}
[c]{ll}%
S_{k}\left(  \mathcal{T}\right)  & \text{for any }\mathbf{z}\in\mathcal{V},\\
S_{k,0}\left(  \mathcal{T}\right)  & \text{if }\mathbf{z}\in\mathcal{V}%
_{\Omega},
\end{array}
\right. \\
\forall\mathbf{y}\in\mathcal{V}:\psi_{k}^{\mathbf{z}}\left(  \mathbf{y}%
\right)  =\left.  \psi_{k}^{\mathbf{z}}\right\vert _{K}\left(  \mathbf{y}%
\right)  =\left(  1+\left(  -1\right)  ^{k+1}\rho_{k}\right)  \delta
_{\mathbf{y},\mathbf{z}}%
\end{array}
\label{LagrangeProp}%
\end{equation}

\end{lemma}%

%TCIMACRO{\TeXButton{Proof}{\proof}}%
%BeginExpansion
\proof
%EndExpansion
Let $\mathbf{z}\in\mathcal{V}$. The function $\psi_{k}^{\mathbf{z}}$ has
support in $\omega_{\mathbf{z}}$ and for any $F\in\mathcal{F}_{\mathbf{z}}$
and adjacent simplex $K\in\mathcal{T}$ it holds%
\[
\left.  \psi_{k}^{\mathbf{z}}\right\vert _{K}=\sum_{F\in\mathcal{F}%
_{\mathbf{z}}}\left.  B_{k}^{\operatorname*{CR},F}\right\vert _{K}%
=\sum_{\mathbf{y}\in\mathcal{V}\left(  K\right)  \backslash\left\{
\mathbf{z}\right\}  }\left.  \left(  P_{k}^{\left(  0,d-2\right)  }\left(
1-2\lambda_{K,\mathbf{y}}\right)  -B_{k}^{\operatorname*{CR},K}\right)
\right\vert _{K}.
\]

For $K$, $K^{\prime}\in\mathcal{T}_{F}$ and $\mathbf{y}\in\mathcal{V}\left(
F\right)  $ it holds $\left.  \lambda_{K,\mathbf{y}}\right\vert _{F}=\left.
\lambda_{K^{\prime},\mathbf{y}}\right\vert _{F}$ so that the function
$\psi_{k}^{\mathbf{z}}$ is continuous across $F$ so that $\left.  \psi
_{k}^{\mathbf{z}}\right\vert _{\omega_{\mathbf{z}}}$ is continuous.

On the facet $F_{\mathbf{z}}$ in $K$ opposite to $\mathbf{z}$ it holds%
\begin{align*}
\left.  \left.  \psi_{k}^{\mathbf{z}}\right\vert _{K}\right\vert
_{F_{\mathbf{z}}}  &  =\sum_{\mathbf{y}\in\mathcal{V}\left(  F_{\mathbf{z}%
}\right)  }P_{k}^{\left(  0,d-2\right)  }\left(  1-2\left.  \lambda
_{K,\mathbf{y}}\right\vert _{F_{\mathbf{z}}}\right)  -\left(
-1_{F_{\mathbf{z}}}+\sum_{\mathbf{y}\in\mathcal{V}\left(  K\right)  }%
P_{k}^{\left(  0,d-2\right)  }\left(  1-2\left.  \lambda_{K,\mathbf{y}%
}\right\vert _{F_{\mathbf{z}}}\right)  \right) \\
&  =\sum_{\mathbf{y}\in\mathcal{V}\left(  F_{\mathbf{z}}\right)  }%
P_{k}^{\left(  0,d-2\right)  }\left(  1-2\left.  \lambda_{K,\mathbf{y}%
}\right\vert _{F_{\mathbf{z}}}\right)  -\left(  \sum_{\mathbf{y}\in
\mathcal{V}\left(  F_{\mathbf{z}}\right)  }P_{k}^{\left(  0,d-2\right)
}\left(  1-2\left.  \lambda_{K,\mathbf{y}}\right\vert _{F_{\mathbf{z}}%
}\right)  \right)  =0.
\end{align*}
Hence, in general $\psi_{k}^{\mathbf{z}}\in S_{k}\left(  \mathcal{T}\right)  $
with $\operatorname*{supp}\psi_{k}^{\mathbf{z}}\subset\omega_{\mathbf{z}}$ and
for inner vertices $\mathbf{z}\in\mathcal{V}_{\Omega}$ it follows $\psi
_{k}^{\mathbf{z}}\in S_{k,0}\left(  \mathcal{T}\right)  $.

We use the vertex values of the non-conforming function (\ref{BkCRvertex}) to
evaluate $\psi_{k}^{\mathbf{z}}$ at the simplex vertices.%
%TCIMACRO{\TeXButton{End Proof}{\endproof}}%
%BeginExpansion
\endproof
%EndExpansion

\begin{lemma}
\label{LemRedRep}For any $k\in\mathbb{N}$, it holds%
\begin{equation}
\operatorname*{CR}\nolimits_{k}\left(  \mathcal{T}\right)  =\left\{
\begin{array}
[c]{ll}%
\overset{\bullet}{S_{k}}\left(  \mathcal{T}\right)  +V_{k}^{\operatorname*{nc}%
}\left(  \mathcal{F}\right)  & k\text{ odd,}\\
S_{k}\left(  \mathcal{T}\right)  +V_{k}^{\operatorname*{nc}}\left(
\mathcal{T}\right)  & k\text{ even,}%
\end{array}
\right.  \label{CRequal0}%
\end{equation}
and%
\begin{equation}
\operatorname*{CR}\nolimits_{k,0}\left(  \mathcal{T}\right)  =\left\{
\begin{array}
[c]{ll}%
\overset{\bullet}{S_{k,0}}\left(  \mathcal{T}\right)  +V_{k}%
^{\operatorname*{nc}}\left(  \mathcal{F}_{\Omega}\right)  & k\text{ odd,}\\
S_{k,0}\left(  \mathcal{T}\right)  +V_{k}^{\operatorname*{nc}}\left(
\mathcal{T}\right)  & k\text{ even.}%
\end{array}
\right.  \label{CRequal}%
\end{equation}

\end{lemma}%

%TCIMACRO{\TeXButton{Proof}{\proof}}%
%BeginExpansion
\proof
%EndExpansion
The second statements in (\ref{CRequal0}) and (\ref{CRequal}) are simply
repetitions of the definition.

We only prove the first statement in (\ref{CRequal}); the first one in
(\ref{CRequal0}) follows by the same arguments. We employ the function
$\psi_{k}^{\mathbf{z}}$ as in Lemma \ref{LemDefPsikz}. The combination of
(\ref{TechThingsb}) and (\ref{LagrangeProp}) implies $\psi_{k}^{\mathbf{z}%
}\left(  \mathbf{z}\right)  \neq0$ so that the scaled version $\tilde{\psi
}_{k}^{\mathbf{z}}=\psi_{k}^{\mathbf{z}}/\psi_{k}^{\mathbf{z}}\left(
\mathbf{z}\right)  $ is well defined, belongs to $S_{k,0}\left(
\mathcal{T}\right)  $, and satisfies for $\mathbf{z}\in\mathcal{V}_{\Omega}$:%
\[
\tilde{\psi}_{k}^{\mathbf{z}}\left(  \mathbf{y}\right)  =\delta_{\mathbf{y}%
,\mathbf{z}}\quad\forall\mathbf{y}\in\mathcal{V}\text{.}%
\]
Any $v\in\operatorname*{CR}_{k,0}\left(  \mathcal{T}\right)  $ can be written
in the form%
\[
v=w^{\operatorname{c}}+w_{\mathcal{F}}^{\operatorname*{nc}}\text{\quad
for\quad}w^{\operatorname{c}}\in S_{k,0}\left(  \mathcal{T}\right)  \text{,
}w_{\mathcal{F}}^{\operatorname*{nc}}\in V_{k}^{\operatorname*{nc}}\left(
\mathcal{F}_{\Omega}\right)  \text{.}%
\]
Then,
\[
w^{\operatorname{c}}=w_{1}^{\operatorname{c}}+w_{0}^{\operatorname{c}}%
\quad\text{for }w_{1}^{\operatorname{c}}:=\sum_{\mathbf{z}\in\mathcal{V}%
_{\Omega}}w^{\operatorname{c}}\left(  \mathbf{z}\right)  \tilde{\psi}%
_{k}^{\mathbf{z}}\quad\text{and\quad}w_{0}^{\operatorname{c}}\in
\overset{\bullet}{S_{k,0}}\left(  \mathcal{T}\right)  \text{.}%
\]
Since $\tilde{\psi}_{k}^{\mathbf{z}}\in V_{k}^{\operatorname*{nc}}\left(
\mathcal{F}_{\Omega}\right)  $ it follows $w_{1}^{\operatorname{c}}\in
V_{k}^{\operatorname*{nc}}\left(  \mathcal{F}_{\Omega}\right)  $ and we have
derived the representation $v=w_{0}^{\operatorname{c}}+w_{1}^{\operatorname{c}%
}+w_{\mathcal{F}}^{\operatorname*{nc}}$ with $w_{0}^{\operatorname{c}}%
\in\overset{\bullet}{S_{k,0}}\left(  \mathcal{T}\right)  $ and $w_{1}%
^{\operatorname{c}}+w_{\mathcal{F}}^{\operatorname*{nc}}\in V_{k}%
^{\operatorname*{nc}}\left(  \mathcal{F}_{\Omega}\right)  $.%
%TCIMACRO{\TeXButton{End Proof}{\endproof}}%
%BeginExpansion
\endproof
%EndExpansion

\begin{corollary}
\label{CorCont}The conforming finite element space are contained in the
Crouzeix-Raviart spaces for $k\geq1$:%
\[
S_{k}\subset\operatorname*{CR}\nolimits_{k}\left(  \mathcal{T}\right)
\text{\quad and }S_{k,0}\subset\operatorname*{CR}\nolimits_{k,0}\left(
\mathcal{T}\right)  .
\]

\end{corollary}

\begin{lemma}
\label{LemBasisPartial}For $k\geq2$, a basis for $V_{k}^{\operatorname*{nc}%
}\left(  \mathcal{T}\right)  $ is given by $\mathcal{B}_{k}%
^{\operatorname*{nc}}\left(  \mathcal{T}\right)  :=\left\{  B_{k}%
^{\operatorname*{CR},K},K\in\mathcal{T}\right\}  $ while for $k=1$ it holds
$V_{1}^{\operatorname*{nc}}\left(  \mathcal{T}\right)  =\left\{  0\right\}  $.

For odd $k\geq1$, a basis for $V_{k}^{\operatorname*{nc}}\left(
\mathcal{F}_{\Omega}\right)  $ is given by $\mathcal{B}_{k}%
^{\operatorname*{nc}}\left(  \mathcal{F}_{\Omega}\right)  :=\left\{
B_{k}^{\operatorname*{CR},F},F\in\mathcal{F}_{\Omega}\right\}  $ and a basis
for $V_{k}^{\operatorname*{nc}}\left(  \mathcal{F}\right)  $ by $\mathcal{B}%
_{k}^{\operatorname*{nc}}\left(  \mathcal{F}\right)  :=\left\{  B_{k}%
^{\operatorname*{CR},F},F\in\mathcal{F}\right\}  $.
\end{lemma}

%

%TCIMACRO{\TeXButton{Proof}{\proof}}%
%BeginExpansion
\proof
%EndExpansion
From the definition of $V_{k}^{\operatorname*{nc}}\left(  \mathcal{T}\right)
$ and $V_{k}^{\operatorname*{nc}}\left(  \mathcal{F}_{\Omega}\right)  $ it
follows that it is sufficient to prove that each of the sets $\mathcal{B}%
_{k}^{\operatorname*{nc}}\left(  \mathcal{T}\right)  $, $\mathcal{B}%
_{k}^{\operatorname*{nc}}\left(  \mathcal{F}_{\Omega}\right)  $ consists of
linearly independent functions.

\textbf{Case 1}, \textbf{the space }$V_{k}^{\operatorname*{nc}}\left(
\mathcal{F}_{\Omega}\right)  $.

It is sufficient to prove that
\begin{equation}
\left(  0=\sum_{F\in\mathcal{F}_{\Omega}}\alpha_{F}B_{k}^{\operatorname*{CR}%
,F}\right)  \implies\left(  \forall F\in\mathcal{F}_{\Omega}\quad\alpha
_{F}=0\right)  \text{.} \label{0linunabh}%
\end{equation}
To simplify the notation in the sequel we define $\alpha_{F}:=0$ for
$F\in\mathcal{F}_{\partial\Omega}$ so that%
\[
\sum_{F\in\mathcal{F}_{\Omega}}\alpha_{F}B_{k}^{\operatorname*{CR},F}%
=\sum_{F\in\mathcal{F}}\alpha_{F}B_{k}^{\operatorname*{CR},F}.
\]

Let $K\in\mathcal{T}$ and note that%
\begin{equation}
\left.  \left(  \sum_{F\in\mathcal{F}_{\Omega}}\alpha_{F}B_{k}%
^{\operatorname*{CR},F}\right)  \right\vert _{K}=\sum_{F\in\mathcal{F}\left(
K\right)  }\alpha_{F}\left.  B_{k}^{\operatorname*{CR},F}\right\vert _{K}.
\label{sumfloc}%
\end{equation}
The vertex values in (\ref{BkCRvertex}) combined with the assumption in
(\ref{0linunabh}) and (\ref{sumfloc}) imply the following vertex conditions:%
\begin{equation}
-d\alpha_{F_{\mathbf{y}}}+\sum_{F\in\mathcal{F}\left(  K\right)  }\alpha
_{F}=0\quad\forall\mathbf{y}\in\mathcal{V}\left(  K\right)  , \label{eqfy}%
\end{equation}
where $F_{\mathbf{y}}$ is the facet of $K$ opposite to $\mathbf{y}$. We
collect the coefficients in $%
%TCIMACRO{\TeXButton{boldalphanew}{\boldsymbol{\alpha}}}%
%BeginExpansion
\boldsymbol{\alpha}%
%EndExpansion
=\left(  \alpha_{F}\right)  _{F\in\mathcal{F}\left(  K\right)  }$ and define
the matrix%
\begin{equation}
\mathbf{Q}_{d}\left(  s\right)  =\left[
\begin{array}
[c]{cccc}%
-s & 1 & \ldots & 1\\
1 & -s & \ddots & \vdots\\
\vdots & \ddots & \ddots & 1\\
1 & \ldots & 1 & -s
\end{array}
\right]  \in\mathbb{R}^{d\times d}, \label{defQd}%
\end{equation}
so that (\ref{eqfy}) becomes equivalent to%
\begin{equation}
\mathbf{Q}_{d+1}\left(  d-1\right)
%TCIMACRO{\TeXButton{boldalphanew}{\boldsymbol{\alpha}}}%
%BeginExpansion
\boldsymbol{\alpha}%
%EndExpansion
=\mathbf{0}. \label{solQd}%
\end{equation}
In Appendix \ref{AppDetForm} we will prove that the determinant is given by%
\begin{equation}
\det\mathbf{Q}_{d}\left(  s\right)  =\left(  -1\right)  ^{d+1}\left(
1+s\right)  ^{d-1}\left(  d-1-s\right)  \label{defformula}%
\end{equation}
so that%
\begin{equation}
\det\mathbf{Q}_{d+1}\left(  d-1\right)  =\left(  -d\right)  ^{d}\neq0.
\label{Qdregular}%
\end{equation}

Hence, the system (\ref{solQd}) has the unique solution $%
%TCIMACRO{\TeXButton{boldalphanew}{\boldsymbol{\alpha}}}%
%BeginExpansion
\boldsymbol{\alpha}%
%EndExpansion
=\mathbf{0}$. Since $K\in\mathcal{T}$ was arbitrary we conclude that
(\ref{0linunabh}) holds.

The proof that the functions in $\mathcal{B}_{k}^{\operatorname*{nc}}\left(
\mathcal{F}\right)  $ are linearly independent is verbatim.

\textbf{Case 2}, \textbf{the space }$V_{k}^{\operatorname*{nc}}\left(
\mathcal{T}\right)  $.

The supports of the functions $B_{k}^{\operatorname*{CR},K}$ have pairwise
disjoint interior and so the functions in $\mathcal{B}_{k}^{\operatorname*{nc}%
}\left(  \mathcal{T}\right)  $ are $L^{2}\left(  \Omega\right)  $-orthogonal.
Hence, it suffices to show that $B_{k}^{\operatorname*{CR},K}$ is not the zero
function for $k\geq2$ and $B_{1}^{\operatorname*{CR},K}=0$. For the reference
element, it holds%
\[
d\left.  B_{k}^{\operatorname*{CR},\hat{K}}\right\vert _{\hat{K}}\left(
\mathbf{x}\right)  =\left(  \sum_{\ell=1}^{d}P_{k}^{\left(  0,d-2\right)
}\left(  1-2x_{\ell}\right)  \right)  +P_{k}^{\left(  0,d-2\right)  }\left(
-\left(  1-2\left\vert \mathbf{x}\right\vert \right)  \right)  -1\quad
\forall\mathbf{x}=\left(  x_{\ell}\right)  _{\ell=1}^{d}\in\hat{K}.
\]
We employ the representation (see \cite[18.5.7]{NIST:DLMF}):%

\begin{align*}
P_{k}^{\left(  0,d-2\right)  }\left(  1-2s\right)   &  =\sum_{\ell=0}%
^{k}\left(  -1\right)  ^{\ell}\binom{k+\ell+d-2}{\ell}\binom{k}{\ell}s^{\ell
}\\
P_{k}^{\left(  0,d-2\right)  }\left(  -\left(  1-2s\right)  \right)   &
\overset{\text{\cite[Table 18.6.1]{NIST:DLMF}}}{=}\left(  -1\right)  ^{k}%
P_{k}^{\left(  d-2,0\right)  }\left(  1-2s\right) \\
&  =\left(  -1\right)  ^{k}\sum_{\ell=0}^{k}\left(  -1\right)  ^{\ell}%
\binom{k+\ell+d-2}{\ell}\binom{k+d-2}{\ell+d-2}s^{\ell}%
\end{align*}
to obtain%
\begin{align}
d\left.  B_{k}^{\operatorname*{CR},\hat{K}}\right\vert _{\hat{K}}\left(
\mathbf{x}\right)   &  =-1+\sum_{\ell=0}^{k}\left(  -1\right)  ^{\ell}%
\binom{k}{\ell}\binom{k+\ell+d-2}{\ell}\sum_{j=1}^{d}x_{j}^{\ell
}\label{RepBkpoly}\\
&  +\left(  -1\right)  ^{k}\sum_{\ell=0}^{k}\left(  -1\right)  ^{\ell}%
\binom{k+\ell+d-2}{\ell}\binom{k+d-2}{\ell+d-2}\left\vert \mathbf{x}%
\right\vert ^{\ell}\\
&  =\binom{2k+d-2}{k}\left(  \left(  \left(  -1\right)  ^{k}\sum_{j=1}%
^{d}x_{j}^{k}\right)  +\left\vert \mathbf{x}\right\vert ^{k}\right)
+p_{k-1}\left(  \mathbf{x}\right)  ,\nonumber
\end{align}
where $p_{k-1}\in\mathbb{P}_{k-1}$. It is easy to see that for $k\geq2$ the
highest order term $\left(  \left(  \left(  -1\right)  ^{k}\sum_{j=1}^{d}%
x_{j}^{k}\right)  +\left\vert \mathbf{x}\right\vert ^{k}\right)  $ is not the
zero function and $\left.  B_{k}^{\operatorname*{CR},\hat{K}}\right\vert
_{\hat{K}}\neq0$.

For $k=1$, we obtain from (\ref{RepBkpoly})%
\[
d\left.  B_{1}^{\operatorname*{CR},\hat{K}}\right\vert _{\hat{K}}\left(
\mathbf{x}\right)  =-1+d-d\left\vert \mathbf{x}\right\vert -\left(
d-1\right)  +d\left\vert \mathbf{x}\right\vert =0.
\]
%

%TCIMACRO{\TeXButton{End Proof}{\endproof}}%
%BeginExpansion
\endproof
%EndExpansion

\begin{lemma}
\label{LemBasis}With the Notation as in Lemma \ref{LemBasisPartial}, a basis
for $\operatorname*{CR}_{k}\left(  \mathcal{T}\right)  $ is given by%
\begin{equation}
\mathcal{B}_{k}^{\operatorname*{CR}}\left(  \mathcal{T}\right)  :=\left\{
\begin{array}
[c]{lc}%
\left\{  B_{k}^{\mathbf{z}}\mid\forall\mathbf{z}\in\overset{\bullet
}{\mathcal{N}_{k}}\left(  \mathcal{T}\right)  \right\}  \cup\mathcal{B}%
_{k}^{\operatorname*{nc}}\left(  \mathcal{F}\right)  & k\text{ odd,}\\
\left\{  B_{k}^{\mathbf{z}}\mid\forall\mathbf{z}\in\mathcal{N}_{k}\left(
\mathcal{T}\right)  \right\}  \cup\mathcal{B}_{k}^{\operatorname*{nc}}\left(
\mathcal{T}^{\prime}\right)  & k\text{ even,}%
\end{array}
\right.  \label{generalbasisdef}%
\end{equation}
where $\mathcal{T}^{\prime}\subset\mathcal{T}$ is any submesh with $\left\vert
\mathcal{T}^{\prime}\right\vert =\left\vert \mathcal{T}\right\vert -1$.

A basis for the subspace $\operatorname*{CR}_{k,0}\left(  \mathcal{T}\right)
$ is given by%
\begin{equation}
\mathcal{B}_{k,0}^{\operatorname*{CR}}\left(  \mathcal{T}\right)  :=\left\{
\begin{array}
[c]{lc}%
\left\{  B_{k}^{\mathbf{z}}\mid\forall\mathbf{z}\in\overset{\bullet
}{\mathcal{N}_{k,0}}\left(  \mathcal{T}\right)  \right\}  \cup\mathcal{B}%
_{k}^{\operatorname*{nc}}\left(  \mathcal{F}_{\Omega}\right)  & k\text{
odd,}\\
\left\{  B_{k}^{\mathbf{z}}\mid\forall\mathbf{z}\in\mathcal{N}_{k,0}\left(
\mathcal{T}\right)  \right\}  \cup\mathcal{B}_{k}^{\operatorname*{nc}}\left(
\mathcal{T}\right)  & k\text{ even.}%
\end{array}
\right.  \label{generalbasisdef0}%
\end{equation}

\end{lemma}%

%TCIMACRO{\TeXButton{Proof}{\proof}}%
%BeginExpansion
\proof
%EndExpansion
\textbf{Case 1: }$k\geq1$ is odd.

For this case, we only prove that $\mathcal{B}_{k}^{\operatorname*{CR}}\left(
\mathcal{T}\right)  $ is a basis of $\operatorname*{CR}_{k}\left(
\mathcal{T}\right)  $ while the proof of (\ref{generalbasisdef0}) is verbatim.
From (\ref{CRequal}) and the definition of $V_{k}^{\operatorname*{nc}}\left(
\mathcal{F}\right)  $ and $V_{k}^{\operatorname*{nc}}\left(  \mathcal{T}%
\right)  $ we conclude that it is sufficient to prove%

\begin{equation}%
\begin{array}
[c]{cl}
& \underset{=:w_{\operatorname{c}}}{\underbrace{\sum_{\mathbf{z}\in
\overset{\bullet}{\mathcal{N}_{k}}\left(  \mathcal{T}\right)  }\gamma
_{\mathbf{z}}B_{k}^{\mathbf{z}}}}+\underset{=:w_{\mathcal{F}}}{\underbrace
{\sum_{F\in\mathcal{F}}\alpha_{F}B_{k}^{\operatorname*{CR},F}}}=0\\
\Longrightarrow & \text{all coefficients }\gamma_{\mathbf{z}}\text{, }%
\alpha_{F}\text{ are zero.}%
\end{array}
\label{equalzerocond}%
\end{equation}
Let $K\in\mathcal{T}$ and%
\begin{equation}
w_{K}:=\sum_{F\in\mathcal{F}\left(  K\right)  }\alpha_{F}\left.
B_{k}^{\operatorname*{CR},F}\right\vert _{K}. \label{equalzerofull}%
\end{equation}
Since $w_{\operatorname{c}}\left(  \mathbf{y}\right)  =0$ for all
$\mathbf{y}\in\mathcal{V}$, condition (\ref{equalzerocond}) implies the vertex
condition $w_{K}\left(  \mathbf{y}\right)  =0$ for all $\mathbf{y}%
\in\mathcal{V}\left(  K\right)  $. Reasoning as in Case 1 of the proof of
Lemma \ref{LemBasisPartial} we conclude that all coefficients $\alpha_{F}$ are
zero. It is well known that the conforming basis functions $B_{k}^{\mathbf{z}%
}$ are linearly independent so that also all coefficients $\gamma_{\mathbf{z}%
}$ in (\ref{equalzerocond}) are zero.\medskip

\textbf{Case 2: }$k\geq2$ is even.

\textbf{Part 1: }We first prove that $\mathcal{B}_{k}^{\operatorname*{CR}%
}\left(  \mathcal{T}\right)  $ in (\ref{generalbasisdef}) defines a basis of
$\operatorname*{CR}_{k}\left(  \mathcal{T}\right)  $. The statement is trivial
if $\mathcal{T}$ contains only a single simplex since then $\operatorname*{CR}%
_{k}\left(  \mathcal{T}\right)  =S_{k}\left(  \mathcal{T}\right)  $ and
$\mathcal{B}_{k}^{\operatorname*{CR}}\left(  \mathcal{T}\right)  =\left\{
B_{k}^{\mathbf{z}}\mid\forall\mathbf{z}\in\mathcal{N}_{k}\left(
\mathcal{T}\right)  \right\}  $ is the standard Lagrange basis. It remains to
consider the case $\left\vert \mathcal{T}\right\vert \geq2$.

\textbf{1st step:} We show%
\begin{equation}
V_{k}^{\operatorname*{nc}}\left(  \mathcal{T}\right)  \cap S_{k}\left(
\mathcal{T}\right)  =\operatorname*{span}\left\{  \Psi_{k}\right\}
\quad\text{with\quad}\Psi_{k}:=\sum_{K\in\mathcal{T}}B_{k}^{\operatorname*{CR}%
,K}. \label{Vknccup}%
\end{equation}
Let $w\in V_{k}^{\operatorname*{nc}}\left(  \mathcal{T}\right)  $. For any
$K\in\mathcal{T}$ it holds $\left.  w\right\vert _{K}=\beta_{K}\left.
B_{k}^{\operatorname*{CR},K}\right\vert _{K}$. Let $K^{\prime}$ be an adjacent
simplex: $F:=K\cap K^{\prime}\in\mathcal{F}_{\Omega}$. From (\ref{BkCRvertexa}%
) and (\ref{TechThingsa}) it follows $\left.  \left.  B_{k}%
^{\operatorname*{CR},K}\right\vert _{K}\right\vert _{F}=\left.  \left.
B_{k}^{\operatorname*{CR},K^{\prime}}\right\vert _{K^{\prime}}\right\vert
_{F}\neq0$ so that the condition $w\in S_{k}\left(  \mathcal{T}\right)  $
implies that $\beta_{K^{\prime}}=\beta_{K}$. By iterating this argument over
adjacent simplices we obtain $\beta_{K^{\prime\prime}}=\beta_{K}$ for all
$K^{\prime\prime}\in\mathcal{T}$ and (\ref{Vknccup}) follows.

\textbf{2nd step: }Since $\Psi_{k}\notin V_{k}^{\operatorname*{nc}}\left(
\mathcal{T}^{\prime}\right)  $ we conclude that $V_{k}^{\operatorname*{nc}%
}\left(  \mathcal{T}^{\prime}\right)  \cap S_{k}\left(  \mathcal{T}\right)
=\left\{  0\right\}  $ and the sum in $\widetilde{\operatorname*{CR}_{k}%
}:=V_{k}^{\operatorname*{nc}}\left(  \mathcal{T}^{\prime}\right)  \oplus
S_{k}\left(  \mathcal{T}\right)  $ is direct. From Lemma \ref{LemBasisPartial}
and the well-known fact that the functions in $\left\{  B_{k}^{\mathbf{z}}%
\mid\forall\mathbf{z}\in\mathcal{N}_{k}\left(  \mathcal{T}\right)  \right\}  $
are linearly independent we conclude that $\mathcal{B}_{k}^{\operatorname*{CR}%
}\left(  \mathcal{T}\right)  $ is a basis of $\widetilde{\operatorname*{CR}%
_{k}}$. It remains to show $\widetilde{\operatorname*{CR}_{k}}%
=\operatorname*{CR}_{k}\left(  \mathcal{T}\right)  $. Since $\widetilde
{\operatorname*{CR}_{k}}\subset\operatorname*{CR}_{k}\left(  \mathcal{T}%
\right)  $ it is sufficient to prove $B_{k}^{\operatorname*{CR},K}%
\in\widetilde{\operatorname*{CR}_{k}}$ for $\left\{  K\right\}  =\mathcal{T}%
\backslash\mathcal{T}^{\prime}$. Let $\tilde{\Psi}_{k}:=\sum_{K^{\prime}%
\in\mathcal{T}^{\prime}}B_{k}^{\operatorname*{CR},K^{\prime}}\in
\widetilde{\operatorname*{CR}_{k}}$. Then%
\[
\Psi_{k}=B_{k}^{\operatorname*{CR},K}+\tilde{\Psi}_{k}\in S_{k}\left(
\mathcal{T}\right)  ,
\]
i.e., $B_{k}^{\operatorname*{CR},K}=\Psi_{k}-\tilde{\Psi}_{k}$ with $\Psi
_{k}\in S_{k}\left(  \mathcal{T}\right)  \subset\widetilde{\operatorname*{CR}%
_{k}}$ and $\tilde{\Psi}_{k}\in\widetilde{\operatorname*{CR}_{k}}$. \medskip

\textbf{Part 2: }It remains to prove that $\mathcal{B}_{k,0}%
^{\operatorname*{CR}}\left(  \mathcal{T}\right)  $ in (\ref{generalbasisdef0})
defines a basis of $\operatorname*{CR}_{k,0}\left(  \mathcal{T}\right)  $. We
generalize the proof \cite[Thm. 22]{ccss_2012} in 2D to general dimension $d$.
It is sufficient to prove the conclusion%
\begin{equation}%
\begin{array}
[c]{cl}
& \underset{=:w_{\operatorname{c}}}{\underbrace{\sum_{\mathbf{z}\in
\mathcal{N}_{k,0}\left(  \mathcal{T}\right)  }\gamma_{\mathbf{z}}%
B_{k}^{\mathbf{z}}}}+\underset{=:w_{\mathcal{T}}}{\underbrace{\sum
_{K\in\mathcal{T}}\beta_{K}B_{k}^{\operatorname*{CR},K}}}=0\\
\Longrightarrow & \text{all coefficients }\gamma_{\mathbf{z}}\text{, }%
\beta_{F}\text{ are zero.}%
\end{array}
\label{linindep}%
\end{equation}
Consider some $K\in\mathcal{T}$ which has at least one facet, say $F$, on
$\partial\Omega$. Then%
\[
\left.  \left(  w_{\operatorname{c}}+w_{\mathcal{T}}\right)  \right\vert
_{F}=\left.  w_{\mathcal{T}}\right\vert _{F}=\beta_{K}\left.  B_{k}%
^{\operatorname*{CR},K}\right\vert _{F}.
\]
From (\ref{def:CRtetrahedron}), (\ref{BkCRvertex}), and (\ref{TechThingsa}) we
conclude that for $\mathbf{y}\in\mathcal{V}\left(  F\right)  :$%
\[
\left.  B_{k}^{\operatorname*{CR},K}\right\vert _{F}\left(  \mathbf{y}\right)
=\frac{d-1+\left(  -1\right)  ^{k}\rho_{k}}{d}\neq0,
\]
so that $\left.  \left(  w_{\operatorname{c}}+w_{\mathcal{T}}\right)
\right\vert _{F}=0$ implies $\beta_{K}=0$ and, in turn%
\begin{equation}
\left(  \left.  \left(  w_{\operatorname{c}}+w_{\mathcal{T}}\right)
\right\vert _{K}=\left.  w_{\operatorname{c}}\right\vert _{K}=0\right)
\implies\left(  \gamma_{\mathbf{z}}=0\quad\forall\mathbf{z}\in\mathcal{N}%
_{k}\cap K\right)  . \label{2ndpartLinIndep}%
\end{equation}

This argument allows us to set up an induction step. We set $\Omega^{\prime
}:=\Omega\backslash K$ and $\mathcal{T}^{\prime}:=\mathcal{T}\backslash
\left\{  K\right\}  $. Then, the function $\left(  w_{\operatorname{c}%
}+w_{\mathcal{T}}\right)  _{\Omega^{\prime}}$ vanishes on $\partial
\Omega^{\prime}$; we again choose some $K^{\prime}\in\mathcal{T}^{\prime}$
with one facet on $\partial\Omega^{\prime}$ and conclude as before that the
corresponding coefficients $\beta_{K^{\prime}}$ and $\gamma_{\mathbf{z}}$,
$\mathbf{z}\in K^{\prime}$ are zero.
%TCIMACRO{\TeXButton{End Proof}{\endproof}}%
%BeginExpansion
\endproof
%EndExpansion

These findings are collected as the following theorem.

\begin{theorem}
\label{ThmCRk0drsplit}The sums in (\ref{CRequal}) are direct: For any
$k\in\mathbb{N}$, it holds%
\[
\operatorname*{CR}\nolimits_{k}\left(  \mathcal{T}\right)  =\left\{
\begin{array}
[c]{ll}%
\overset{\bullet}{S_{k}}\left(  \mathcal{T}\right)  \oplus V_{k}%
^{\operatorname*{nc}}\left(  \mathcal{F}\right)  & k\text{ odd,}\\
S_{k}\left(  \mathcal{T}\right)  \oplus V_{k}^{\operatorname*{nc}}\left(
\mathcal{T}^{\prime}\right)  & k\text{ even,}%
\end{array}
\right.
\]
where $\mathcal{T}^{\prime}\subset\mathcal{T}$ is any submesh with $\left\vert
\mathcal{T}^{\prime}\right\vert =\left\vert \mathcal{T}\right\vert -1$ and%
\[
\operatorname*{CR}\nolimits_{k,0}\left(  \mathcal{T}\right)  =\left\{
\begin{array}
[c]{ll}%
\overset{\bullet}{S_{k,0}}\left(  \mathcal{T}\right)  \oplus V_{k}%
^{\operatorname*{nc}}\left(  \mathcal{F}_{\Omega}\right)  & k\text{ odd,}\\
S_{k,0}\left(  \mathcal{T}\right)  \oplus V_{k}^{\operatorname*{nc}}\left(
\mathcal{T}\right)  & k\text{ even.}%
\end{array}
\right.
\]

For odd $k\geq1$, a basis for $V_{k}^{\operatorname*{nc}}\left(
\mathcal{F}\right)  $ is given by $\mathcal{B}_{k}^{\operatorname*{nc}}\left(
\mathcal{F}\right)  =\left\{  B_{k}^{\operatorname*{CR},F},F\in\mathcal{F}%
\right\}  $ and for $V_{k}^{\operatorname*{nc}}\left(  \mathcal{F}_{\Omega
}\right)  $ by $\mathcal{B}_{k}^{\operatorname*{nc}}\left(  \mathcal{F}%
_{\Omega}\right)  =\left\{  B_{k}^{\operatorname*{CR},F},F\in\mathcal{F}%
_{\Omega}\right\}  $ while, for $k\geq2$, a basis for $V_{k}%
^{\operatorname*{nc}}\left(  \mathcal{T}\right)  $ is given by $\mathcal{B}%
_{k}^{\operatorname*{nc}}\left(  \mathcal{T}\right)  =\left\{  B_{k}%
^{\operatorname*{CR},K},K\in\mathcal{T}\right\}  $. For $k=1$ it holds
$V_{1}^{\operatorname*{nc}}\left(  \mathcal{T}\right)  =\left\{  0\right\}  $.

It holds%
\begin{equation}
S_{k}\left(  \mathcal{T}\right)  \subset\operatorname*{CR}\nolimits_{k}\left(
\mathcal{T}\right)  \subset\operatorname*{CR}\nolimits_{k}^{\max}\left(
\mathcal{T}\right)  \label{Eq:Inclusion Sk in CRk}%
\end{equation}
and%
\begin{equation}
S_{k,0}\left(  \mathcal{T}\right)  \subset\operatorname*{CR}\nolimits_{k,0}%
\left(  \mathcal{T}\right)  \subset\operatorname*{CR}\nolimits_{k,0}^{\max
}\left(  \mathcal{T}\right)  . \label{inclusions}%
\end{equation}
The first inclusion in (\ref{inclusions}) is strict if $k$ is even or if $k$
is odd and $\mathcal{T}$ consists of more than a single simplex. The second
inclusion in (\ref{inclusions}) is strict, e.g., for $d=3$ and $k\geq2$.
\end{theorem}

%

%TCIMACRO{\TeXButton{Proof}{\proof}}%
%BeginExpansion
\proof
%EndExpansion
The first claims are direct consequences of the previous lemmata. The first
inclusion in (\ref{inclusions}) is strict if $\operatorname*{CR}_{k,0}\left(
\mathcal{T}\right)  $ contains at least one non-conforming shape function.
This is always the case if $k\geq2$ is even and for odd $k$ if $\mathcal{T}$
contains at least one inner facet, i.e., if $\mathcal{T}$ consists of more
than a single simplex. The second inclusion is strict, e.g., for $d=3$ and any
$k\geq2$ (see \cite[Sec. 5.3.1, Table 5.3.2]{CDS}). For $d=2$ however, the
second inclusion becomes an equality (see \cite{Baran_Stoyan},
\cite{ccss_2012}).%
%TCIMACRO{\TeXButton{End Proof}{\endproof}}%
%BeginExpansion
\endproof
%EndExpansion

\section{A basis for Crouzeix-Raviart spaces $\operatorname*{CR}_{k}\left(
\mathcal{T}\right)  $ and $\operatorname*{CR}_{k,0}\left(  \mathcal{T}\right)
$\label{SecBasisFull}}

In this section we will present a basis for the Crouzeix-Raviart space which
is composed of orthogonal polynomials multiplied by face bubble functions.
This allows for a simple construction of the degrees of freedom for the face
basis functions as facet and simplex integrals with weight functions such they
are well defined on $\operatorname*{CR}_{k,0}\left(  \mathcal{T}\right)  $.

Theorem \ref{ThmCRk0drsplit} directly give guidelines for a construction: any
(standard) basis for $\overset{\bullet}{S_{k}}\left(  \mathcal{T}\right)  $,
$\overset{\bullet}{S_{k,0}}\left(  \mathcal{T}\right)  $ appended by non
conforming facet functions $B_{k}^{\operatorname*{CR},F}$ for odd $k$ and by
non-conforming simplex functions $B_{k}^{\operatorname*{CR},K}$ form a basis
of the Crouzeix-Raviart spaces. However, the construction of the local bidual
degrees of freedom is a non-trivial task and depends on the choice of basis
for the conforming parts. Since the functions in $\operatorname*{CR}%
_{k,0}\left(  \mathcal{T}\right)  $ are in general discontinuous across facets
while their moments up to a polynomial degree $k-1$ are continuous, the
degrees of freedom associated to the inner facets must be moments of the form
$J_{\mu}^{F}\left(  u\right)  =\int_{F}g_{\mu}^{F}u$ for polynomial weights
$g_{\mu}^{F}$ up to a degree $k-1$. Orthogonality relations on facets will
become important for the choice of $g_{\mu}^{F}$ and this is the reason to
define a basis based on the products of face bubbles with orthogonal
polynomials on the simplicial entities (edges, facets, etc.).

\subsection{The simplicial complex of a simplicial finite element mesh}

Let $\mathcal{T}$ be a conforming finite element mesh for $\Omega$ consisting
of closed simplices $K\in\mathcal{T}$. Let $\mathcal{S}$ denote the associate
simplicial complex, i.e., the set of $\ell$-(dimensional) simplices (faces),
$\ell\in\left\{  0,1,\ldots,d\right\}  $, that satisfies the following conditions:

\begin{enumerate}
\item Every face\footnote{A \textit{face} of a simplex is the convex hull of a
non-empty subset of the simplex vertices. The \textit{facets} of a simplex are
the ($d-1$ dimensional) faces with exactly $d$ vertices.} of a simplex
$\tau\in\mathcal{S}$ is also in $\mathcal{S}$.

\item The non-empty intersection of any two simplices $\tau_{1},\tau_{2}%
\in\mathcal{S}$ is a face of both, $\tau_{1}$ and $\tau_{2}$.
\end{enumerate}

We distinguish between the boundary faces $\mathcal{S}_{\partial\Omega
}:=\left\{  \tau\in\mathcal{S}\mid\tau\in\partial\Omega\right\}  $ and the
inner faces $\mathcal{S}_{\Omega}:=\mathcal{S}\backslash\mathcal{S}%
_{\partial\Omega}$. The subset of $\mathcal{S}$ which contains all $\ell
$-simplices in $\mathcal{S}$ is $\mathcal{S}_{\ell}$. Specifically, the set of
vertices $\mathcal{V}$ is $\mathcal{S}_{0}$, the set of edges $\mathcal{E}$ is
$\mathcal{S}_{1}$ and the set of facets $\mathcal{F}$ equals $\mathcal{S}%
_{d-1}$. The subset of $\mathcal{S}_{\ell}$ which contains all those $\ell
$-simplices which are \textit{not} a subset of the boundary $\partial\Omega$
is denoted by $\mathcal{S}_{\ell,\Omega}$. The corresponding
\textit{skeletons} $\Sigma_{\ell}$, $\Sigma_{\ell,\Omega}$ are given by%
\begin{equation}
\Sigma_{\ell}:=%
%TCIMACRO{\dbigcup \limits_{\tau\in\mathcal{S}_{\ell}}}%
%BeginExpansion
{\displaystyle\bigcup\limits_{\tau\in\mathcal{S}_{\ell}}}
%EndExpansion
\tau\quad\text{and\quad}\Sigma_{\ell,\Omega}:=%
%TCIMACRO{\dbigcup \limits_{\tau\in\mathcal{S}_{\ell,\Omega}}}%
%BeginExpansion
{\displaystyle\bigcup\limits_{\tau\in\mathcal{S}_{\ell,\Omega}}}
%EndExpansion
\tau. \label{skeletondef}%
\end{equation}

In analogy to (\ref{DefPatches}) we define for $\tau\in\mathcal{S}_{\ell}$ the
adjacent simplex patch by%
\[
\mathcal{T}_{\tau}:=\left\{  K\in\mathcal{T}\mid\tau\subset K\right\}
\quad\text{and\quad}\omega_{\tau}:=%
%TCIMACRO{\tbigcup \nolimits_{K\in\mathcal{T}_{\tau}}}%
%BeginExpansion
{\textstyle\bigcup\nolimits_{K\in\mathcal{T}_{\tau}}}
%EndExpansion
K.
\]

\subsection{The face bubble functions}

For $\mathbf{z}\in\mathcal{V}$, let $\varphi_{\mathbf{z}}^{1}\in S_{1}\left(
\mathcal{T}\right)  $ be the \textquotedblleft hat function\textquotedblright%
\ for the vertex $\mathbf{z}$ characterized by the condition $\varphi
_{\mathbf{z}}^{1}\left(  \mathbf{y}\right)  =\delta_{\mathbf{z},\mathbf{y}}$
for all $\mathbf{z},\mathbf{y}\in\mathcal{V}$. For a face $\tau$, the set of
vertices is denoted by $\mathcal{V}\left(  \tau\right)  $, the face bubble
$W_{\tau}\in S_{k,0}\left(  \mathcal{T}\right)  $ is given by%
\[
W_{\tau}=%
%TCIMACRO{\dprod \limits_{\mathbf{z}\in\mathcal{V}\left(  \tau\right)  }}%
%BeginExpansion
{\displaystyle\prod\limits_{\mathbf{z}\in\mathcal{V}\left(  \tau\right)  }}
%EndExpansion
\varphi_{\mathbf{z}}^{1}\qquad\text{and satisfies }\operatorname*{supp}%
W_{\tau}=\mathcal{T}_{\tau}\text{.}%
\]

\subsection{Orthogonal polynomials on the reference simplex}

In \cite[\S 2.5.2]{Yuan_inproc} orthogonal polynomials on the $\ell
$-dimensional reference element%
\[
\hat{K}_{\ell}:=\left\{  \mathbf{x}=\left(  x_{j}\right)  _{j=1}^{\ell}%
\in\mathbb{R}_{\geq0}^{\ell}\mid%
%TCIMACRO{\tsum \nolimits_{j=1}^{\ell}}%
%BeginExpansion
{\textstyle\sum\nolimits_{j=1}^{\ell}}
%EndExpansion
x_{j}\leq1\right\}
\]
are defined which are generalizations of the univariate Jacobi polynomials as
we will explain in the following. Let $\mathbf{\hat{z}}_{0}:=\mathbf{0}$,
$\mathbf{\hat{z}}_{j}:=\mathbf{e}_{j}$, $1\leq j\leq\ell$ denote the vertices
of $\hat{K}_{\ell}$ (see (\ref{Defzhat})), where $\mathbf{e}_{j}$ denotes the
$j$-th canonical unit vector in $\mathbb{R}^{\ell}$, and let $\hat{\lambda
}_{j}$ denote the barycentric coordinate in $\hat{K}_{\ell}$ for the vertex
$\mathbf{\hat{z}}_{j}$. In \cite[\S 2.5.2]{Yuan_inproc}, the polynomials $P_{%
%TCIMACRO{\TeXButton{boldalphanew}{\boldsymbol{\alpha}}}%
%BeginExpansion
\boldsymbol{\alpha}%
%EndExpansion
}$ are defined for $%
%TCIMACRO{\TeXButton{boldalphanew}{\boldsymbol{\alpha}}}%
%BeginExpansion
\boldsymbol{\alpha}%
%EndExpansion
\in\mathbb{N}_{0}^{\ell}$, by first introducing the formal expression for
$\mathbf{y}=\left(  y_{j}\right)  _{j=0}^{\ell}\in\mathbb{R}^{\ell+1}$%
\[
P_{%
%TCIMACRO{\TeXButton{boldalphanew}{\boldsymbol{\alpha}}}%
%BeginExpansion
\boldsymbol{\alpha}%
%EndExpansion
}^{\operatorname*{bary}}\left(  \mathbf{y}\right)  =%
%TCIMACRO{\dprod \limits_{j=1}^{\ell}}%
%BeginExpansion
{\displaystyle\prod\limits_{j=1}^{\ell}}
%EndExpansion
\left(  y_{0}+\sum_{m=j}^{\ell}y_{m}\right)  ^{\alpha_{j}}P_{\alpha_{j}%
}^{\left(  s_{j},1\right)  }\left(  \frac{2y_{j}}{y_{0}+\sum_{m=j}^{\ell}%
y_{m}}-1\right)
\]
with%
\[
s_{j}=2\left(  \sum_{m=j+1}^{\ell}\alpha_{m}\right)  +2\left(  \ell-j\right)
+1
\]
and then inserting the barycentric coordinates $\widehat{%
%TCIMACRO{\TeXButton{boldlambda}{\boldsymbol{\lambda}}}%
%BeginExpansion
\boldsymbol{\lambda}%
%EndExpansion
}=\left(  \hat{\lambda}_{j}\right)  _{j=0}^{\ell}:$%
\[
P_{%
%TCIMACRO{\TeXButton{boldalphanew}{\boldsymbol{\alpha}}}%
%BeginExpansion
\boldsymbol{\alpha}%
%EndExpansion
}\left(  \mathbf{\hat{x}}\right)  :=P_{%
%TCIMACRO{\TeXButton{boldalphanew}{\boldsymbol{\alpha}}}%
%BeginExpansion
\boldsymbol{\alpha}%
%EndExpansion
}^{\operatorname*{bary}}\left(  \widehat{%
%TCIMACRO{\TeXButton{boldlambda}{\boldsymbol{\lambda}}}%
%BeginExpansion
\boldsymbol{\lambda}%
%EndExpansion
}\left(  \mathbf{\hat{x}}\right)  \right)  \quad\forall\mathbf{\hat{x}}\in
\hat{K}_{\ell}.
\]
They are orthogonal with respect to the scalar product
\[
\left(  u,v\right)  _{\hat{K}_{\ell}}:=\int_{\hat{K}_{\ell}}W_{\hat{K}_{\ell}%
}uv\quad\text{with weight function }W_{\hat{K}_{\ell}}:=%
%TCIMACRO{\dprod \limits_{\ell=0}^{d}}%
%BeginExpansion
{\displaystyle\prod\limits_{\ell=0}^{d}}
%EndExpansion
\hat{\lambda}_{\ell},
\]
i.e.,%
\[
\left(  P_{%
%TCIMACRO{\TeXButton{boldalphanew}{\boldsymbol{\alpha}}}%
%BeginExpansion
\boldsymbol{\alpha}%
%EndExpansion
},P_{%
%TCIMACRO{\TeXButton{boldbeta}{\boldsymbol{\beta}}}%
%BeginExpansion
\boldsymbol{\beta}%
%EndExpansion
}\right)  _{\hat{K}_{\ell}}=\delta_{%
%TCIMACRO{\TeXButton{boldalphanew}{\boldsymbol{\alpha}}}%
%BeginExpansion
\boldsymbol{\alpha}%
%EndExpansion
,%
%TCIMACRO{\TeXButton{boldbeta}{\boldsymbol{\beta}}}%
%BeginExpansion
\boldsymbol{\beta}%
%EndExpansion
}c_{%
%TCIMACRO{\TeXButton{boldalphanew}{\boldsymbol{\alpha}}}%
%BeginExpansion
\boldsymbol{\alpha}%
%EndExpansion
}%
\]
for some constants $c_{%
%TCIMACRO{\TeXButton{boldalphanew}{\boldsymbol{\alpha}}}%
%BeginExpansion
\boldsymbol{\alpha}%
%EndExpansion
}\neq0$ (the value of $c_{%
%TCIMACRO{\TeXButton{boldalphanew}{\boldsymbol{\alpha}}}%
%BeginExpansion
\boldsymbol{\alpha}%
%EndExpansion
}$ is given in \cite[\S 2.5.2]{Yuan_inproc}). In this way, the set $\left\{
P_{%
%TCIMACRO{\TeXButton{boldalphanew}{\boldsymbol{\alpha}}}%
%BeginExpansion
\boldsymbol{\alpha}%
%EndExpansion
}\mid%
%TCIMACRO{\TeXButton{boldalphanew}{\boldsymbol{\alpha}}}%
%BeginExpansion
\boldsymbol{\alpha}%
%EndExpansion
\in\mathbb{N}_{\leq m}^{d}\right\}  $ forms an orthogonal basis of the space
$\mathbb{P}_{m}\left(  \hat{K}_{\ell}\right)  $ of $\ell-$variate polynomials
on $\hat{K}_{\ell}$ of maximal total degree $m$.

\begin{remark}
\label{Rembubbledeg}Let $\ell\in\left\{  1,2,\ldots,d\right\}  $. Let
$m\geq\ell+1$. The set
\[
\left\{  W_{\hat{K}_{\ell}}P_{%
%TCIMACRO{\TeXButton{boldalphanew}{\boldsymbol{\alpha}}}%
%BeginExpansion
\boldsymbol{\alpha}%
%EndExpansion
}\mid%
%TCIMACRO{\TeXButton{boldalphanew}{\boldsymbol{\alpha}}}%
%BeginExpansion
\boldsymbol{\alpha}%
%EndExpansion
\in\mathbb{N}_{\leq m-\left(  \ell+1\right)  }^{\ell}\right\}
\]
forms a basis of%
\begin{equation}
\mathbb{P}_{m,0}\left(  \hat{K}_{\ell}\right)  :=\left\{  v\in\mathbb{P}%
_{m}\left(  \hat{K}_{\ell}\right)  \mid\left.  v\right\vert _{\partial\hat
{K}_{\ell}}=0\right\}  . \label{DefPl0}%
\end{equation}

\end{remark}

\begin{example}
The orthogonal polynomials $P_{%
%TCIMACRO{\TeXButton{boldalphanew}{\boldsymbol{\alpha}}}%
%BeginExpansion
\boldsymbol{\alpha}%
%EndExpansion
}\left(  \mathbf{\hat{x}}\right)  $ have the explicit form

\begin{enumerate}
\item for $\ell=1$:%
\[
P_{\alpha}\left(  \hat{x}\right)  =P_{\alpha}^{\left(  1,1\right)  }\left(
2\hat{\lambda}_{1}\left(  \hat{x}\right)  -1\right)
\]

\item and for $\ell=2$:%
\[
P_{\left(  \alpha_{1},\alpha_{2}\right)  }\left(  \mathbf{\hat{x}}\right)
=P_{\alpha_{1}}^{\left(  2\alpha_{2}+3,1\right)  }\left(  2\hat{\lambda}%
_{1}\left(  \mathbf{\hat{x}}\right)  -1\right)  \left(  1-\hat{\lambda}%
_{1}\left(  \mathbf{\hat{x}}\right)  \right)  ^{\alpha_{2}}P_{\alpha_{2}%
}^{\left(  1,1\right)  }\left(  \frac{2\hat{\lambda}_{2}\left(  \mathbf{\hat
{x}}\right)  }{1-\hat{\lambda}_{1}\left(  \mathbf{\hat{x}}\right)  }-1\right)
.
\]

\end{enumerate}
\end{example}

\subsection{Basis functions associated with simplex faces}

To define basis functions associated to a face $\tau\in\mathcal{S}_{\ell}$ we
fix a numbering of its vertices $\mathbf{A}_{j}\left(  \tau\right)  $, $0\leq
j\leq\ell$. For any adjacent simplex $K\in\mathcal{T}_{\tau}$, we choose an
affine pullback $\chi_{K}:\hat{K}_{d}\rightarrow K$ such that the restriction
to $\hat{K}_{\ell}$ satisfies%
\[
\chi_{K}\left(  \mathbf{\hat{x}}\right)  :=\mathbf{A}_{0}\left(  \tau\right)
+\sum_{j=1}^{\ell}\hat{x}_{j}\left(  \mathbf{A}_{j}\left(  \tau\right)
-\mathbf{A}_{0}\left(  \tau\right)  \right)  \quad\forall\mathbf{\hat{x}%
}=\left(  \hat{x}_{j}\right)  _{j=1}^{d}\in\hat{K}_{d}\text{ with }\left(
x_{j}\right)  _{j=\ell+1}^{d}=\mathbf{0}.
\]
It is important to note that the restrictions $\left.  \chi_{K}\right\vert
_{\hat{K}_{\ell}}$ coincide for all $K\in\mathcal{T}_{\tau}$ due to the chosen
fixed numbering of vertices in $\tau$, i.e., $\tau$ is parametrized in the
same way for all adjacent simplices. According to Remark \ref{Rembubbledeg}
basis functions associated to a face $\tau\in\mathcal{S}_{\ell}$ exist only if
$k\geq\ell+1$ which we assume for the following. We employ the index sets%
\[%
\begin{array}
[c]{ll}%
\mathcal{I}_{k,\ell}:=\left\{  \left(  \tau,%
%TCIMACRO{\TeXButton{boldalphanew}{\boldsymbol{\alpha}}}%
%BeginExpansion
\boldsymbol{\alpha}%
%EndExpansion
\right)  :\tau\in\mathcal{S}_{\ell}\quad%
%TCIMACRO{\TeXButton{boldalphanew}{\boldsymbol{\alpha}}}%
%BeginExpansion
\boldsymbol{\alpha}%
%EndExpansion
\in\mathbb{N}_{k-\left(  \ell+1\right)  }^{\ell}\right\}  , & \mathcal{I}%
_{k,\ell,\Omega}:=\left\{  \left(  \tau,%
%TCIMACRO{\TeXButton{boldalphanew}{\boldsymbol{\alpha}}}%
%BeginExpansion
\boldsymbol{\alpha}%
%EndExpansion
\right)  \in\mathcal{I}_{k,\ell}:\tau\in\mathcal{S}_{\ell,\Omega}\right\}  ,\\
\mathcal{I}_{k}:=\left\{
\begin{array}
[c]{ll}%
%TCIMACRO{\dbigcup \limits_{\ell=1}^{\min\left\{  k-1,d\right\}  }}%
%BeginExpansion
{\displaystyle\bigcup\limits_{\ell=1}^{\min\left\{  k-1,d\right\}  }}
%EndExpansion
\mathcal{I}_{k,\ell} & \text{if }k\text{ is odd,}\\%
%TCIMACRO{\dbigcup \limits_{\ell=0}^{\min\left\{  k-1,d\right\}  }}%
%BeginExpansion
{\displaystyle\bigcup\limits_{\ell=0}^{\min\left\{  k-1,d\right\}  }}
%EndExpansion
\mathcal{I}_{k,\ell} & \text{if }k\text{ is even,}%
\end{array}
\right.  & \mathcal{I}_{k,\Omega}:=\left\{
\begin{array}
[c]{ll}%
%TCIMACRO{\dbigcup \limits_{\ell=1}^{\min\left\{  k-1,d\right\}  }}%
%BeginExpansion
{\displaystyle\bigcup\limits_{\ell=1}^{\min\left\{  k-1,d\right\}  }}
%EndExpansion
\mathcal{I}_{k,\ell,\Omega} & \text{if }k\text{ is odd,}\\%
%TCIMACRO{\dbigcup \limits_{\ell=0}^{\min\left\{  k-1,d\right\}  }}%
%BeginExpansion
{\displaystyle\bigcup\limits_{\ell=0}^{\min\left\{  k-1,d\right\}  }}
%EndExpansion
\mathcal{I}_{k,\ell,\Omega} & \text{if }k\text{ is even,}%
\end{array}
\right. \\
\mathcal{I}_{k,\partial\Omega}:=\mathcal{I}_{k}\backslash\mathcal{I}%
_{k,\Omega}. &
\end{array}
\]

\begin{definition}
\label{DefConfBasis}Let $k\in\mathbb{N}_{\geq1}$. For $\ell\in\left\{
0,1,\ldots,k-1\right\}  $, the conforming basis functions associated to a face
$\tau\in\mathcal{S}_{\ell}$ are given for

\begin{enumerate}
\item $\ell=0$ by%
\[
B_{\mathbf{z},%
%TCIMACRO{\TeXButton{boldalphanew}{\boldsymbol{\alpha}}}%
%BeginExpansion
\boldsymbol{\alpha}%
%EndExpansion
}=\varphi_{\mathbf{z}}^{1}\quad\forall\mathbf{z}\in\mathcal{V}\quad%
%TCIMACRO{\TeXButton{boldalphanew}{\boldsymbol{\alpha}}}%
%BeginExpansion
\boldsymbol{\alpha}%
%EndExpansion
\in\mathbb{N}_{\leq k-1}^{0}:=\left\{  0\right\}  ,
\]

\item $\ell\geq1$ by%
\begin{equation}
B_{\tau,%
%TCIMACRO{\TeXButton{boldalphanew}{\boldsymbol{\alpha}}}%
%BeginExpansion
\boldsymbol{\alpha}%
%EndExpansion
}:=\left\{
\begin{array}
[c]{ll}%
W_{\tau}P_{%
%TCIMACRO{\TeXButton{boldalphanew}{\boldsymbol{\alpha}}}%
%BeginExpansion
\boldsymbol{\alpha}%
%EndExpansion
}\circ\chi_{K}^{-1} & K\in\mathcal{T}_{\tau},\\
0 & \text{otherwise}%
\end{array}
\right\vert \quad\forall%
%TCIMACRO{\TeXButton{boldalphanew}{\boldsymbol{\alpha}}}%
%BeginExpansion
\boldsymbol{\alpha}%
%EndExpansion
\in\mathbb{N}_{\leq k-\left(  \ell+1\right)  }^{\ell} \label{defBtau}%
\end{equation}

\end{enumerate}

and collected in%
\begin{align*}
\mathcal{B}_{\ell}\left(  \mathcal{T}\right)   &  :=\left\{  B_{\tau,%
%TCIMACRO{\TeXButton{boldalphanew}{\boldsymbol{\alpha}}}%
%BeginExpansion
\boldsymbol{\alpha}%
%EndExpansion
}\mid\forall\left(  \tau,%
%TCIMACRO{\TeXButton{boldalphanew}{\boldsymbol{\alpha}}}%
%BeginExpansion
\boldsymbol{\alpha}%
%EndExpansion
\right)  \in\mathcal{I}_{k,\ell}\right\}  ,\\
\mathcal{B}_{\ell,0}\left(  \mathcal{T}\right)   &  :=\left\{  B_{\tau,%
%TCIMACRO{\TeXButton{boldalphanew}{\boldsymbol{\alpha}}}%
%BeginExpansion
\boldsymbol{\alpha}%
%EndExpansion
}\mid\forall\left(  \tau,%
%TCIMACRO{\TeXButton{boldalphanew}{\boldsymbol{\alpha}}}%
%BeginExpansion
\boldsymbol{\alpha}%
%EndExpansion
\right)  \in\mathcal{I}_{k,\ell,\Omega}\right\}  .
\end{align*}

\end{definition}

\begin{remark}
The function $B_{\tau,%
%TCIMACRO{\TeXButton{boldalphanew}{\boldsymbol{\alpha}}}%
%BeginExpansion
\boldsymbol{\alpha}%
%EndExpansion
}$ in (\ref{defBtau}) can be expressed without the pullback $\chi_{K}$ by
setting $%
%TCIMACRO{\TeXButton{boldvarphi}{\mbox{\boldmath$ \varphi$}}}%
%BeginExpansion
\mbox{\boldmath$ \varphi$}%
%EndExpansion
_{\tau}^{1}:=\left(  \varphi_{\mathbf{A}_{i}\left(  \tau\right)  }^{1}\right)
_{i=0}^{\ell}$ and%
\[
B_{\tau,%
%TCIMACRO{\TeXButton{boldalphanew}{\boldsymbol{\alpha}}}%
%BeginExpansion
\boldsymbol{\alpha}%
%EndExpansion
}=W_{\tau}P_{%
%TCIMACRO{\TeXButton{boldalphanew}{\boldsymbol{\alpha}}}%
%BeginExpansion
\boldsymbol{\alpha}%
%EndExpansion
}^{\operatorname*{bary}}\left(
%TCIMACRO{\TeXButton{boldvarphi}{\mbox{\boldmath$ \varphi$}}}%
%BeginExpansion
\mbox{\boldmath$ \varphi$}%
%EndExpansion
_{\tau}^{1}\right)  .
\]

\end{remark}

\subsection{A basis for $\operatorname*{CR}_{k}\left(  \mathcal{T}\right)  $
and $\operatorname*{CR}_{k,0}\left(  \mathcal{T}\right)  $}

In this section, we will prove that the union of the basis functions in
$\mathcal{B}_{\ell}\left(  \mathcal{T}\right)  $ with the non-conforming
facet/simplex functions yields a basis of the Crouzeix-Raviart spaces.

\begin{definition}
\label{DefCRBasis}Let $k\in\mathbb{N}_{\geq1}$. A basis for
$\operatorname*{CR}_{k}\left(  \mathcal{T}\right)  $ is given by%
\begin{equation}
\mathcal{B}_{k}^{\operatorname*{CR}}\left(  \mathcal{T}\right)  :=\left\{
\begin{array}
[c]{ll}%
\left(
%TCIMACRO{\dbigcup \limits_{\ell=1}^{\min\left\{  k-1,d\right\}  }}%
%BeginExpansion
{\displaystyle\bigcup\limits_{\ell=1}^{\min\left\{  k-1,d\right\}  }}
%EndExpansion
\mathcal{B}_{\ell}\left(  \mathcal{T}\right)  \right)  \cup\left\{
B_{k}^{\operatorname*{CR},F}:F\in\mathcal{F}\right\}  & \text{if }k\text{ is
odd,}\\
\left(
%TCIMACRO{\dbigcup \limits_{\ell=0}^{\min\left\{  k-1,d\right\}  }}%
%BeginExpansion
{\displaystyle\bigcup\limits_{\ell=0}^{\min\left\{  k-1,d\right\}  }}
%EndExpansion
\mathcal{B}_{\ell}\left(  \mathcal{T}\right)  \right)  \cup\left\{
B_{k}^{\operatorname*{CR},K}:K\in\mathcal{T}\right\}  & \text{if }k\text{ is
even}%
\end{array}
\right.  \label{basisfull}%
\end{equation}
and for $\operatorname*{CR}_{k,0}\left(  \mathcal{T}\right)  $ by%
\begin{equation}
\mathcal{B}_{k,0}^{\operatorname*{CR}}\left(  \mathcal{T}\right)  :=\left\{
\begin{array}
[c]{ll}%
\left(
%TCIMACRO{\dbigcup \limits_{\ell=1}^{\min\left\{  k-1,d\right\}  }}%
%BeginExpansion
{\displaystyle\bigcup\limits_{\ell=1}^{\min\left\{  k-1,d\right\}  }}
%EndExpansion
\mathcal{B}_{\ell,0}\left(  \mathcal{T}\right)  \right)  \cup\left\{
B_{k}^{\operatorname*{CR},F}:F\in\mathcal{F}_{\Omega}\right\}  & \text{if
}k\text{ is odd,}\\
\left(
%TCIMACRO{\dbigcup \limits_{\ell=0}^{\min\left\{  k-1,d\right\}  }}%
%BeginExpansion
{\displaystyle\bigcup\limits_{\ell=0}^{\min\left\{  k-1,d\right\}  }}
%EndExpansion
\mathcal{B}_{\ell,0}\left(  \mathcal{T}\right)  \right)  \cup\left\{
B_{k}^{\operatorname*{CR},K}:K\in\mathcal{T}\right\}  & \text{if }k\text{ is
even.}%
\end{array}
\right.  \label{basiszero}%
\end{equation}

\end{definition}

The following lemma states that the notion \textquotedblleft
basis\textquotedblright\ in Definition \ref{DefCRBasis} indeed is justified.

\begin{lemma}
\label{LemLinInd}The set $\mathcal{B}_{k}^{\operatorname*{CR}}\left(
\mathcal{T}\right)  $ forms a basis of $\operatorname*{CR}_{k}\left(
\mathcal{T}\right)  $ and $\mathcal{B}_{k,0}^{\operatorname*{CR}}\left(
\mathcal{T}\right)  $ of $\operatorname*{CR}_{k,0}\left(  \mathcal{T}\right)
$.
\end{lemma}

%

%TCIMACRO{\TeXButton{Proof}{\proof}}%
%BeginExpansion
\proof
%EndExpansion
We only prove the first statement. Due to Theorem \ref{ThmCRk0drsplit} it is
sufficient to prove that
\[
\mathcal{\tilde{B}}_{k}^{\operatorname*{CR}}\left(  \mathcal{T}\right)
:=\left\{
\begin{array}
[c]{ll}%
%TCIMACRO{\dbigcup \limits_{\ell=1}^{\min\left\{  k-1,d\right\}  }}%
%BeginExpansion
{\displaystyle\bigcup\limits_{\ell=1}^{\min\left\{  k-1,d\right\}  }}
%EndExpansion
\mathcal{B}_{\ell}\left(  \mathcal{T}\right)  & \text{if }k\text{ is odd,}\\%
%TCIMACRO{\dbigcup \limits_{\ell=0}^{\min\left\{  k-1,d\right\}  }}%
%BeginExpansion
{\displaystyle\bigcup\limits_{\ell=0}^{\min\left\{  k-1,d\right\}  }}
%EndExpansion
\mathcal{B}_{\ell}\left(  \mathcal{T}\right)  & \text{if }k\text{ is even}%
\end{array}
\right.
\]
is a basis for $\overset{\bullet}{S_{k}}\left(  \mathcal{T}\right)  $ for odd
$k$ and for $S_{k}\left(  \mathcal{T}\right)  $ for even $k$. We start with
the linear independence and assume $0\leq\ell\leq\min\left\{  k-1,d\right\}
$. The construction of the basis functions $B_{\tau,%
%TCIMACRO{\TeXButton{boldalphanew}{\boldsymbol{\alpha}}}%
%BeginExpansion
\boldsymbol{\alpha}%
%EndExpansion
}$ via the face bubbles $W_{\tau}$ imply that%
\[
\left.  B_{\tau,%
%TCIMACRO{\TeXButton{boldalphanew}{\boldsymbol{\alpha}}}%
%BeginExpansion
\boldsymbol{\alpha}%
%EndExpansion
}\right\vert _{\tau^{\prime}}=0\quad\forall\tau\in\mathcal{S}_{\ell}%
\quad\forall\tau^{\prime}\in\mathcal{S}_{\ell^{\prime}}\quad\forall
\ell^{\prime}\leq\ell.
\]
Since $\left.  B_{\tau,%
%TCIMACRO{\TeXButton{boldalphanew}{\boldsymbol{\alpha}}}%
%BeginExpansion
\boldsymbol{\alpha}%
%EndExpansion
}\right\vert _{\tau}$ is not the zero function the direct sum representation
follows%
\[
\sum_{\ell=0}^{\min\left\{  k-1,d\right\}  }\sum_{\tau\in\mathcal{S}_{\ell}%
}\operatorname*{span}\left\{  B_{\tau,%
%TCIMACRO{\TeXButton{boldalphanew}{\boldsymbol{\alpha}}}%
%BeginExpansion
\boldsymbol{\alpha}%
%EndExpansion
}:%
%TCIMACRO{\TeXButton{boldalphanew}{\boldsymbol{\alpha}}}%
%BeginExpansion
\boldsymbol{\alpha}%
%EndExpansion
\in\mathbb{N}_{k-\left(  \ell+1\right)  }^{\ell}\right\}  =%
%TCIMACRO{\dbigoplus \limits_{\ell=0}^{\min\left\{  k-1,d\right\}  }}%
%BeginExpansion
{\displaystyle\bigoplus\limits_{\ell=0}^{\min\left\{  k-1,d\right\}  }}
%EndExpansion%
%TCIMACRO{\dbigoplus \limits_{\tau\in\mathcal{S}_{\ell}}}%
%BeginExpansion
{\displaystyle\bigoplus\limits_{\tau\in\mathcal{S}_{\ell}}}
%EndExpansion
\operatorname*{span}\left\{  B_{\tau,%
%TCIMACRO{\TeXButton{boldalphanew}{\boldsymbol{\alpha}}}%
%BeginExpansion
\boldsymbol{\alpha}%
%EndExpansion
}:%
%TCIMACRO{\TeXButton{boldalphanew}{\boldsymbol{\alpha}}}%
%BeginExpansion
\boldsymbol{\alpha}%
%EndExpansion
\in\mathbb{N}_{k-\left(  \ell+1\right)  }^{\ell}\right\}  .
\]
Remark \ref{Rembubbledeg} implies via a simple counting argument:%
\[
\dim P_{k-\left(  \ell+1\right)  ,0}\left(  \tau\right)  =\binom{k-1}{\ell
}=\left\vert \mathbb{N}_{k-\left(  \ell+1\right)  }^{\ell}\right\vert
\]
so that the functions $B_{\tau,%
%TCIMACRO{\TeXButton{boldalphanew}{\boldsymbol{\alpha}}}%
%BeginExpansion
\boldsymbol{\alpha}%
%EndExpansion
}$, $%
%TCIMACRO{\TeXButton{boldalphanew}{\boldsymbol{\alpha}}}%
%BeginExpansion
\boldsymbol{\alpha}%
%EndExpansion
\in\mathbb{N}_{k-\left(  \ell+1\right)  }^{\ell}$ are also linearly
independent and, in turn, the set $%
%TCIMACRO{\dbigcup \nolimits_{\ell=0}^{\min\left\{  k-1,d\right\}  }}%
%BeginExpansion
{\displaystyle\bigcup\nolimits_{\ell=0}^{\min\left\{  k-1,d\right\}  }}
%EndExpansion
\mathcal{B}_{\ell}\left(  \mathcal{T}\right)  $ consists of linearly
independent functions. It remains to show that any function in $\overset
{\bullet}{S_{k}}\left(  \mathcal{T}\right)  $, $S_{k}\left(  \mathcal{T}%
\right)  $ can be represented as a linear combination of $\mathcal{\tilde{B}%
}_{k}^{\operatorname*{CR}}\left(  \mathcal{T}\right)  $.

\textbf{Case 1: }$k$ is even.

We employ the Lagrange basis as in Lemma \ref{LemBasis}. Recall the definition
of the skeleton $\Sigma_{\ell}$ as in (\ref{skeletondef}) and note that
$\mathcal{N}_{k}\left(  \mathcal{T}\right)  \cap\Sigma_{0}=\mathcal{V}$. Any
function $u\in S_{k}\left(  \mathcal{T}\right)  $ can be written in the form%
\begin{equation}
u=u_{0}+v_{1}\quad\text{for }u_{0}=\sum_{\mathbf{z}\in\mathcal{N}_{k}\left(
\mathcal{T}\right)  \cap\Sigma_{0}}u\left(  \mathbf{z}\right)  B_{k}%
^{\mathbf{z}} \label{usplit}%
\end{equation}
and $v_{1}$ vanishes in all simplex vertices. Remark \ref{Rembubbledeg}
implies that there is a function $w_{0}\in\operatorname*{span}\mathcal{\tilde
{B}}_{k}^{\operatorname*{CR}}\left(  \mathcal{T}\right)  $ such that
\[
u_{0}\left(  \mathbf{z}\right)  =w_{0}\left(  \mathbf{z}\right)  \quad
\forall\mathbf{z\in}\mathcal{V}.
\]
Hence, the difference $u-w_{0}$ vanish in all simplex vertices and $u-w_{0}$
has the representation\footnote{Note that $\mathcal{N}_{k}\left(
\mathcal{T}\right)  \cap\left(  \Sigma_{1}\backslash\Sigma_{0}\right)  $ is
the set of nodal points which lie in the interiors of the edges $E\in
\mathcal{E}$.}%
\[
u-w_{0}=u_{1}+v_{2}\quad\text{for }v_{2}=\sum_{\mathbf{z}\in\mathcal{N}%
_{k}\left(  \mathcal{T}\right)  \cap\left(  \Sigma_{1}\backslash\Sigma
_{0}\right)  }u\left(  \mathbf{z}\right)  B_{k}^{\mathbf{z}}.
\]
Again Remark \ref{Rembubbledeg} implies that the restriction of $u_{1}$ to an
edge can be represented by some $w_{1}\in\operatorname*{span}\mathcal{\tilde
{B}}_{k}^{\operatorname*{CR}}\left(  \mathcal{T}\right)  :$%
\[
\left.  u_{1}\right\vert _{E}=\left.  w_{1}\right\vert _{E}\quad\forall
E\mathbf{\in}\mathcal{S}_{1}\left(  \mathcal{T}\right)
\]
and $u-w_{0}-w_{1}$ vanishes on $\Sigma_{1}$. In this way, an induction
argument with respect to the dimension $\ell$ implies%
\[
u\in\operatorname*{span}\mathcal{\tilde{B}}_{k}^{\operatorname*{CR}}\left(
\mathcal{T}\right)  \text{.}%
\]

\textbf{Case 2:} $k$ is odd.

This case can be concluded by the same arguments as before by observing the
$u_{0}=0$ in (\ref{usplit}) for $u\in\overset{\bullet}{S_{k}}\left(
\mathcal{T}\right)  $.
%TCIMACRO{\TeXButton{End Proof}{\endproof}}%
%BeginExpansion
\endproof
%EndExpansion

\section{Degrees of freedom for the case $d\geq2$ and $k$ is
odd\label{SecDofsGenD}}

In this section, we will define degrees of freedom for the Crouzeix-Raviart
space for general dimension $d\in\left\{  2,3,\ldots\right\}  $ and odd $k$
which form a dual basis. We recall the formal definition.

\begin{definition}
Let $\left\{  b_{i}:i\in\left\{  1,2,\ldots,N\right\}  \right\}  $ be a basis
of a finite dimensional vector space $V_{N}$. A \emph{bidual basis} is a set
of functionals $J_{i}$, $i\in\left\{  1,2,\ldots,N\right\}  $, in the dual
space $V_{N}^{\prime}$ such that%
\[
J_{i}\left(  b_{j}\right)  =\delta_{i,j}\quad\forall i,j\in\left\{
1,2,\ldots,N\right\}  .
\]

\end{definition}

As a very mild assumption on the mesh we require that
\begin{equation}
\mathcal{T}\text{ contains more than only one tetrahedron.}
\label{condTnonsingle}%
\end{equation}

The construction starts with some preliminaries for a single simplex
$K\in\mathcal{T}$. We set%
\[
\left.  L^{2}\left(  \Omega\right)  \right\vert _{K}:=\left\{  \left.
f\right\vert _{K}:f\in L^{2}\left(  \Omega\right)  \right\}
\]
and use this notation also for subspaces of $L^{2}\left(  \Omega\right)  $. We
introduce index sets for $K\in\mathcal{T}$ and $F\in\mathcal{F}$ by%
\begin{equation}%
\begin{tabular}
[c]{ll}%
$\mathcal{S}_{\ell}\left(  K\right)  :=\left\{  \tau\in\mathcal{S}_{\ell}%
\mid\tau\subset K\right\}  ,$ & $\mathcal{S}\left(  K\right)  :=%
%TCIMACRO{\dbigcup \limits_{\ell=1}^{d}}%
%BeginExpansion
{\displaystyle\bigcup\limits_{\ell=1}^{d}}
%EndExpansion
\mathcal{S}_{\ell}\left(  K\right)  ,$\\
$\mathcal{S}_{\ell}\left(  F\right)  :=\left\{  \tau\in\mathcal{S}_{\ell}%
\mid\tau\subset F\right\}  ,$ & $\mathcal{S}\left(  F\right)  :=%
%TCIMACRO{\dbigcup \limits_{\ell=1}^{d}}%
%BeginExpansion
{\displaystyle\bigcup\limits_{\ell=1}^{d}}
%EndExpansion
\mathcal{S}_{\ell}\left(  F\right)  ,$\\
$\mathcal{I}_{k,\ell}\left(  K\right)  :=\left\{  \left(  \tau,%
%TCIMACRO{\TeXButton{boldalphanew}{\boldsymbol{\alpha}}}%
%BeginExpansion
\boldsymbol{\alpha}%
%EndExpansion
\right)  :\tau\in\mathcal{S}_{\ell}\left(  K\right)  \quad%
%TCIMACRO{\TeXButton{boldalphanew}{\boldsymbol{\alpha}}}%
%BeginExpansion
\boldsymbol{\alpha}%
%EndExpansion
\in\mathbb{N}_{\leq k-\left(  \ell+1\right)  }^{\ell}\right\}  ,$ &
$\mathcal{I}_{k}\left(  K\right)  :=%
%TCIMACRO{\dbigcup \limits_{\ell=1}^{d}}%
%BeginExpansion
{\displaystyle\bigcup\limits_{\ell=1}^{d}}
%EndExpansion
\mathcal{I}_{k,\ell}\left(  K\right)  ,$\\
$\mathcal{I}_{k,\ell}\left(  F\right)  :=\left\{  \left(  \tau,%
%TCIMACRO{\TeXButton{boldalphanew}{\boldsymbol{\alpha}}}%
%BeginExpansion
\boldsymbol{\alpha}%
%EndExpansion
\right)  :\tau\in\mathcal{S}_{\ell}\left(  F\right)  \quad%
%TCIMACRO{\TeXButton{boldalphanew}{\boldsymbol{\alpha}}}%
%BeginExpansion
\boldsymbol{\alpha}%
%EndExpansion
\in\mathbb{N}_{\leq k-\left(  \ell+1\right)  }^{\ell}\right\}  ,$ &
$\mathcal{I}_{k}\left(  F\right)  :=%
%TCIMACRO{\dbigcup \limits_{\ell=1}^{d-1}}%
%BeginExpansion
{\displaystyle\bigcup\limits_{\ell=1}^{d-1}}
%EndExpansion
\mathcal{I}_{k,\ell}\left(  F\right)  ,$\\
\multicolumn{2}{l}{$\mathcal{B}_{k}^{\operatorname*{CR}}\left(  K\right)
:=\left\{  \left.  B_{\tau,%
%TCIMACRO{\TeXButton{boldalphanew}{\boldsymbol{\alpha}}}%
%BeginExpansion
\boldsymbol{\alpha}%
%EndExpansion
}\right\vert _{K}:\left(  \tau,%
%TCIMACRO{\TeXButton{boldalphanew}{\boldsymbol{\alpha}}}%
%BeginExpansion
\boldsymbol{\alpha}%
%EndExpansion
\right)  \in\mathcal{I}_{k}\left(  K\right)  \right\}  \cup\left\{  \left.
B_{k}^{\operatorname*{CR},F}\right\vert _{K}:F\in\mathcal{F}\left(  K\right)
\right\}  .$}%
\end{tabular}
\ \ \label{indexsets}%
\end{equation}

Since $S_{k}\left(  \mathcal{T}\right)  \subset\operatorname*{CR}_{k}\left(
\mathcal{T}\right)  $ it holds%
\begin{equation}
\mathbb{P}_{k}\left(  K\right)  =\left.  S_{k}\left(  \mathcal{T}\right)
\right\vert _{K}=\left.  \operatorname*{CR}\nolimits_{k}\left(  \mathcal{T}%
\right)  \right\vert _{K}=\operatorname*{span}\left\{  \mathcal{B}%
_{k}^{\operatorname*{CR}}\left(  K\right)  \right\}  . \label{pkKincls}%
\end{equation}
Since the number $d+1$ of vertices of a simplex equals the number of facets it
follows from the well known fact (cf. (\ref{defNkhut}), (\ref{indexsets}))%
\[
\left\vert \mathcal{N}_{k}\left(  \hat{K}\right)  \right\vert =\left\vert
%TCIMACRO{\dbigcup \limits_{\ell=0}^{d}}%
%BeginExpansion
{\displaystyle\bigcup\limits_{\ell=0}^{d}}
%EndExpansion
\mathcal{I}_{k,\ell}\left(  K\right)  \right\vert =\dim\mathbb{P}_{k}\left(
K\right)
\]
that $\left\vert \mathcal{B}_{k}^{\operatorname*{CR}}\left(  K\right)
\right\vert =\dim\mathbb{P}_{k}\left(  K\right)  $ and the combination with
Theorem \ref{ThmCRk0drsplit} yields that the functions in $\mathcal{B}%
_{k}^{\operatorname*{CR}}\left(  K\right)  $ are linearly independent. From
elementary linear algebra we conclude that there exists a bidual basis for
$\mathcal{B}_{k}^{\operatorname*{CR}}\left(  K\right)  $, induced by $g_{\tau,%
%TCIMACRO{\TeXButton{boldalphanew}{\boldsymbol{\alpha}}}%
%BeginExpansion
\boldsymbol{\alpha}%
%EndExpansion
}^{K}\in L^{2}\left(  K\right)  $, $\left(  \tau,%
%TCIMACRO{\TeXButton{boldalphanew}{\boldsymbol{\alpha}}}%
%BeginExpansion
\boldsymbol{\alpha}%
%EndExpansion
\right)  \in\mathcal{I}_{k}\left(  K\right)  $ and $g_{F}^{\operatorname*{CR}%
,K}\in L^{2}\left(  K\right)  $, $F\in\mathcal{F}\left(  K\right)  $ so that
the associated functionals
\[
J_{\tau,%
%TCIMACRO{\TeXButton{boldalphanew}{\boldsymbol{\alpha}}}%
%BeginExpansion
\boldsymbol{\alpha}%
%EndExpansion
}^{K}u:=\left(  g_{\tau,%
%TCIMACRO{\TeXButton{boldalphanew}{\boldsymbol{\alpha}}}%
%BeginExpansion
\boldsymbol{\alpha}%
%EndExpansion
},u\right)  _{L^{2}\left(  K\right)  }\quad\text{and\quad}J_{F}%
^{\operatorname*{CR},K}u:=\left(  g_{F}^{\operatorname*{CR},K},u\right)
_{L^{2}\left(  K\right)  }%
\]
satisfy%
\[%
\begin{array}
[c]{ll}%
J_{\tau,%
%TCIMACRO{\TeXButton{boldalphanew}{\boldsymbol{\alpha}}}%
%BeginExpansion
\boldsymbol{\alpha}%
%EndExpansion
}^{K}\left(  \left.  B_{t,%
%TCIMACRO{\TeXButton{boldbeta}{\boldsymbol{\beta}}}%
%BeginExpansion
\boldsymbol{\beta}%
%EndExpansion
}\right\vert _{K}\right)  =\delta_{\tau,t}\delta_{%
%TCIMACRO{\TeXButton{boldalphanew}{\boldsymbol{\alpha}}}%
%BeginExpansion
\boldsymbol{\alpha}%
%EndExpansion
,%
%TCIMACRO{\TeXButton{boldbeta}{\boldsymbol{\beta}}}%
%BeginExpansion
\boldsymbol{\beta}%
%EndExpansion
}, & J_{\tau,%
%TCIMACRO{\TeXButton{boldalphanew}{\boldsymbol{\alpha}}}%
%BeginExpansion
\boldsymbol{\alpha}%
%EndExpansion
}^{K}\left(  \left.  B_{k}^{\operatorname*{CR},F}\right\vert _{K}\right)
=0,\\
J_{F}^{\operatorname*{CR},K}\left(  \left.  B_{t,%
%TCIMACRO{\TeXButton{boldbeta}{\boldsymbol{\beta}}}%
%BeginExpansion
\boldsymbol{\beta}%
%EndExpansion
}\right\vert _{K}\right)  =0, & J_{F}^{\operatorname*{CR},K}\left(  \left.
B_{k}^{\operatorname*{CR},F^{\prime}}\right\vert _{K}\right)  =\delta
_{F,F^{\prime}}.
\end{array}
\]

A second ingredient for the definition of degrees of freedom is the selection
of an assignment function $\operatorname*{mark}\nolimits_{\mathcal{T}%
}:\mathcal{S}\rightarrow\mathcal{T\cup F}_{\partial\Omega}$ which satisfies%
%TCIMACRO{\TeXButton{mark}{\begin{subequations}
%\label{mark}
%\end{subequations}}}%
%BeginExpansion
\begin{subequations}
\label{mark}
\end{subequations}%
%EndExpansion%
\begin{align}
\operatorname*{mark}\nolimits_{\mathcal{T}}\left(  \tau\right)   &
\in\mathcal{T}\quad\text{and }\tau\subset\operatorname*{mark}%
\nolimits_{\mathcal{T}}\left(  \tau\right)  \quad\forall\tau\in\mathcal{S}%
_{\Omega},\tag{%
%TCIMACRO{\TeXButton{mark}{\ref{mark}}}%
%BeginExpansion
\ref{mark}%
%EndExpansion
a}\label{marka}\\
\operatorname*{mark}\nolimits_{\mathcal{T}}\left(  \tau\right)   &
\in\mathcal{F}_{\partial\Omega}\quad\text{and }\tau\subset\operatorname*{mark}%
\nolimits_{\mathcal{T}}\left(  \tau\right)  \quad\forall\tau\in\mathcal{S}%
_{\partial\Omega}. \tag{%
%TCIMACRO{\TeXButton{mark}{\ref{mark}}}%
%BeginExpansion
\ref{mark}%
%EndExpansion
b}\label{markb}%
\end{align}

For the final definition of the degrees of freedom for Crouzeix-Raviart spaces
we need the following lemma. For $K\in\mathcal{T}$ and $F,F^{\prime}%
\in\mathcal{F}\left(  K\right)  $, the restrictions of the basis functions in
$\mathcal{B}_{k}^{\operatorname*{CR}}\left(  K\right)  $ (see (\ref{indexsets}%
)) to a simplex facet are denoted by%
\[
B_{F,\tau,%
%TCIMACRO{\TeXButton{boldalphanew}{\boldsymbol{\alpha}}}%
%BeginExpansion
\boldsymbol{\alpha}%
%EndExpansion
}^{K}:=\left.  \left.  B_{F,%
%TCIMACRO{\TeXButton{boldalphanew}{\boldsymbol{\alpha}}}%
%BeginExpansion
\boldsymbol{\alpha}%
%EndExpansion
}\right\vert _{K}\right\vert _{F}\quad\text{and\quad}B_{F,k}%
^{\operatorname*{CR},F^{\prime}}:=\left.  \left.  B_{k}^{\operatorname*{CR}%
,F^{\prime}}\right\vert _{K}\right\vert _{F}.
\]

\begin{lemma}
\label{LemRestrfacet}Let $k$ be odd. For $K\in\mathcal{T}$ and $F\in
\mathcal{F}\left(  K\right)  $, let $\overset{\bullet}{\mathcal{F}}\left(
K\right)  \subset\mathcal{F}\left(  K\right)  $ be a strict subset, i.e.,
$\left\vert \overset{\bullet}{\mathcal{F}}\left(  K\right)  \right\vert
<\left\vert \mathcal{F}\left(  K\right)  \right\vert ,$ and $F\in
\overset{\bullet}{\mathcal{F}}\left(  K\right)  $. Then the functions in%
\begin{equation}
\overset{\bullet}{\mathcal{B}}_{k}^{\operatorname*{CR}}\left(  K,F\right)
:=\left\{  B_{F,\tau,%
%TCIMACRO{\TeXButton{boldalphanew}{\boldsymbol{\alpha}}}%
%BeginExpansion
\boldsymbol{\alpha}%
%EndExpansion
}^{K}:\left(  \tau,%
%TCIMACRO{\TeXButton{boldalphanew}{\boldsymbol{\alpha}}}%
%BeginExpansion
\boldsymbol{\alpha}%
%EndExpansion
\right)  \in\mathcal{I}_{k}\left(  F\right)  \right\}  \cup\left\{
B_{F,k}^{\operatorname*{CR},F^{\prime}}:F^{\prime}\in\overset{\bullet
}{\mathcal{F}}\left(  K\right)  \right\}  \label{DefF}%
\end{equation}
are linearly independent in $L^{2}\left(  F\right)  $.
\end{lemma}

%

%TCIMACRO{\TeXButton{Proof}{\proof}}%
%BeginExpansion
\proof
%EndExpansion
We first prove that the sum of the spans of the two sets in (\ref{DefF}) is
direct while in the second part of the proof we show the functions in each of
the two sets of (\ref{DefF}) are linearly independent.

Let $u\in\operatorname*{span}\overset{\bullet}{\mathcal{B}}_{k}%
^{\operatorname*{CR}}\left(  K,F\right)  $. Then, there is a splitting
$u=u_{0}+u_{1}$ with $u_{1}\in\operatorname*{span}\left\{  B_{F,k}%
^{\operatorname*{CR},F^{\prime}}:F^{\prime}\in\overset{\bullet}{\mathcal{F}%
}\left(  K\right)  \right\}  $ and%
\[
u_{0}\in\operatorname*{span}\left\{  B_{F,\tau,%
%TCIMACRO{\TeXButton{boldalphanew}{\boldsymbol{\alpha}}}%
%BeginExpansion
\boldsymbol{\alpha}%
%EndExpansion
}^{K}:\left(  \tau,%
%TCIMACRO{\TeXButton{boldalphanew}{\boldsymbol{\alpha}}}%
%BeginExpansion
\boldsymbol{\alpha}%
%EndExpansion
\right)  \in\mathcal{I}_{k}\left(  F\right)  \right\}  .
\]
First, we prove $u=0\implies u_{0}=u_{1}=0$. Since $u_{0}$ vanishes in all
simplex vertices, the condition $u=0$ implies%
\begin{equation}
u_{1}\left(  \mathbf{z}\right)  =\sum_{F^{\prime}\in\overset{\bullet
}{\mathcal{F}}\left(  K\right)  }\gamma_{F^{\prime}}B_{F,k}%
^{\operatorname*{CR},F^{\prime}}\left(  \mathbf{z}\right)  =0\quad
\forall\mathbf{z}\in\mathcal{V}\left(  F\right)  . \label{u1gleich0}%
\end{equation}
In a similar fashion as in (\ref{solQd}) this condition can be formulated as a
linear system for the coefficients $%
%TCIMACRO{\TeXButton{boldgamma}{\mbox{\boldmath$ \gamma$}}}%
%BeginExpansion
\mbox{\boldmath$ \gamma$}%
%EndExpansion
_{K}:=\left(  \gamma_{F^{\prime}}\right)  _{F^{\prime}\in\overset{\bullet
}{\mathcal{F}}\left(  K\right)  }$. We choose a numbering of the facets in
$\overset{\bullet}{\mathcal{F}}\left(  K\right)  $ such that $F_{1}=F$ and
number the vertices in $\mathcal{V}\left(  K\right)  $ such that
$\mathbf{z}_{j}$ is opposite to the facet $F_{j}$. From $F\in\overset{\bullet
}{\mathcal{F}}\left(  K\right)  $ and $B_{F,k}^{\operatorname*{CR},F}=1$ we
conclude that the implication \textquotedblleft condition (\ref{u1gleich0}%
)$\implies%
%TCIMACRO{\TeXButton{boldgamma}{\mbox{\boldmath$ \gamma$}}}%
%BeginExpansion
\mbox{\boldmath$ \gamma$}%
%EndExpansion
_{K}=\mathbf{0}$\textquotedblright\ is equivalent to the condition%
\begin{equation}
\det\mathbf{R}_{d}\left(  d-1\right)  \neq0 \label{detne0}%
\end{equation}
with $\mathbf{R}_{d}$ as in (\ref{defRd}). The relation (\ref{detne0}) is
proved in the appendix (Lemma \ref{LemQR}). and $%
%TCIMACRO{\TeXButton{boldgamma}{\mbox{\boldmath$ \gamma$}}}%
%BeginExpansion
\mbox{\boldmath$ \gamma$}%
%EndExpansion
_{K}=\mathbf{0}$ follows. Hence, $u_{1}=0$ which implies $u_{0}=0$.

From (\ref{detne0}) it follows that the functions $B_{F,k}^{\operatorname*{CR}%
,F^{\prime}}$, $F^{\prime}\in\overset{\bullet}{\mathcal{F}}\left(  K\right)
$, are linearly independent. Finally, we show that the functions $\left\{
B_{F,\tau,%
%TCIMACRO{\TeXButton{boldalphanew}{\boldsymbol{\alpha}}}%
%BeginExpansion
\boldsymbol{\alpha}%
%EndExpansion
}^{K}:\left(  \tau,%
%TCIMACRO{\TeXButton{boldalphanew}{\boldsymbol{\alpha}}}%
%BeginExpansion
\boldsymbol{\alpha}%
%EndExpansion
\right)  \in%
%TCIMACRO{\dbigcup \limits_{\ell=1}^{d-1}}%
%BeginExpansion
{\displaystyle\bigcup\limits_{\ell=1}^{d-1}}
%EndExpansion
\mathcal{I}_{k,\ell}\left(  K\right)  \right\}  $ are linearly independent.
The proof is standard and follows by induction over $\ell=1,2,\ldots$, from

\begin{enumerate}
\item the support properties of $B_{F,\tau,%
%TCIMACRO{\TeXButton{boldalphanew}{\boldsymbol{\alpha}}}%
%BeginExpansion
\boldsymbol{\alpha}%
%EndExpansion
}^{K}$, i.e.: $\left.  B_{F,\tau,%
%TCIMACRO{\TeXButton{boldalphanew}{\boldsymbol{\alpha}}}%
%BeginExpansion
\boldsymbol{\alpha}%
%EndExpansion
}^{K}\right\vert _{\tau^{\prime}}=0$ for all $\tau\in\mathcal{S}_{\ell}\left(
K\right)  $ and $\tau^{\prime}\in\left(
%TCIMACRO{\dbigcup \limits_{r=1}^{\ell}}%
%BeginExpansion
{\displaystyle\bigcup\limits_{r=1}^{\ell}}
%EndExpansion
\mathcal{S}_{r}\left(  K\right)  \right)  \backslash\left\{  \tau\right\}  $ and

\item the linear independence of $\left(  \left.  B_{F,\tau,%
%TCIMACRO{\TeXButton{boldalphanew}{\boldsymbol{\alpha}}}%
%BeginExpansion
\boldsymbol{\alpha}%
%EndExpansion
}^{K}\right\vert _{\tau}\right)  _{%
%TCIMACRO{\TeXButton{boldalphanew}{\boldsymbol{\alpha}}}%
%BeginExpansion
\boldsymbol{\alpha}%
%EndExpansion
\in\mathbb{N}_{\leq k-\ell-1}^{\ell}}$; see the proof of Lemma \ref{LemLinInd}%
.%
%TCIMACRO{\TeXButton{End Proof}{\endproof}}%
%BeginExpansion
\endproof
%EndExpansion

\end{enumerate}

For $F\in\mathcal{F}_{\partial\Omega}$ and adjacent $K\in\mathcal{T}$, let
\begin{equation}
\overset{\bullet}{\mathcal{F}}\left(  K\right)  :=\left\{  F^{\prime}%
\in\mathcal{F}\left(  K\right)  \cap\mathcal{F}_{\partial\Omega}\right\}  .
\label{deffpunkt}%
\end{equation}
Clearly, $F\in\overset{\bullet}{\mathcal{F}}\left(  K\right)  $ and assumption
(\ref{condTnonsingle}) imply $\left\vert \overset{\bullet}{\mathcal{F}}\left(
K\right)  \right\vert <\left\vert \mathcal{F}\left(  K\right)  \right\vert $
so that Lemma \ref{LemRestrfacet} becomes applicable. Elementary linear
algebra tells us that there exists a bidual basis for the functions in
$\overset{\bullet}{\mathcal{B}}_{k}^{\operatorname*{CR}}\left(  K,F\right)  $
(see (\ref{DefF})), induced by $g_{\tau,%
%TCIMACRO{\TeXButton{boldalphanew}{\boldsymbol{\alpha}}}%
%BeginExpansion
\boldsymbol{\alpha}%
%EndExpansion
}^{F}\in L^{2}\left(  F\right)  $, $\left(  \tau,%
%TCIMACRO{\TeXButton{boldalphanew}{\boldsymbol{\alpha}}}%
%BeginExpansion
\boldsymbol{\alpha}%
%EndExpansion
\right)  \in\mathcal{I}_{k}\left(  F\right)  $ and $g_{F^{\prime}%
}^{\operatorname*{CR},F}\in L^{2}\left(  F\right)  $, $F^{\prime}\in
\overset{\bullet}{\mathcal{F}}\left(  K\right)  $ such that the associated
functionals
\begin{equation}
J_{\tau,%
%TCIMACRO{\TeXButton{boldalphanew}{\boldsymbol{\alpha}}}%
%BeginExpansion
\boldsymbol{\alpha}%
%EndExpansion
}^{F}u:=\left(  g_{\tau,%
%TCIMACRO{\TeXButton{boldalphanew}{\boldsymbol{\alpha}}}%
%BeginExpansion
\boldsymbol{\alpha}%
%EndExpansion
}^{F},u\right)  _{L^{2}\left(  F\right)  }\quad\text{and\quad}J_{F^{\prime}%
}^{\operatorname*{CR},F}u:=\left(  g_{F^{\prime}}^{\operatorname*{CR}%
,F},u\right)  _{L^{2}\left(  F\right)  } \label{facetdofs}%
\end{equation}
satisfy%
\begin{equation}%
\begin{array}
[c]{ll}%
J_{\tau,%
%TCIMACRO{\TeXButton{boldalphanew}{\boldsymbol{\alpha}}}%
%BeginExpansion
\boldsymbol{\alpha}%
%EndExpansion
}^{F}\left(  B_{F,t,%
%TCIMACRO{\TeXButton{boldbeta}{\boldsymbol{\beta}}}%
%BeginExpansion
\boldsymbol{\beta}%
%EndExpansion
}^{K}\right)  =\delta_{\tau,t}\delta_{%
%TCIMACRO{\TeXButton{boldalphanew}{\boldsymbol{\alpha}}}%
%BeginExpansion
\boldsymbol{\alpha}%
%EndExpansion
,%
%TCIMACRO{\TeXButton{boldbeta}{\boldsymbol{\beta}}}%
%BeginExpansion
\boldsymbol{\beta}%
%EndExpansion
}, & J_{\tau,%
%TCIMACRO{\TeXButton{boldalphanew}{\boldsymbol{\alpha}}}%
%BeginExpansion
\boldsymbol{\alpha}%
%EndExpansion
}^{F}\left(  B_{F^{\prime\prime},k}^{\operatorname*{CR},K}\right)  =0,\\
J_{F^{\prime}}^{\operatorname*{CR},F}\left(  B_{F,t,%
%TCIMACRO{\TeXButton{boldbeta}{\boldsymbol{\beta}}}%
%BeginExpansion
\boldsymbol{\beta}%
%EndExpansion
}^{K}\right)  =0, & J_{F^{\prime}}^{\operatorname*{CR},F}\left(
B_{F^{\prime\prime},k}^{\operatorname*{CR},K}\right)  =\delta_{F^{\prime
},F^{\prime\prime}}%
\end{array}
\label{bidualpropsfacets}%
\end{equation}
for all $\left(  \tau,%
%TCIMACRO{\TeXButton{boldalphanew}{\boldsymbol{\alpha}}}%
%BeginExpansion
\boldsymbol{\alpha}%
%EndExpansion
\right)  ,\left(  t,%
%TCIMACRO{\TeXButton{boldbeta}{\boldsymbol{\beta}}}%
%BeginExpansion
\boldsymbol{\beta}%
%EndExpansion
\right)  \in\mathcal{I}_{k}\left(  F\right)  $ and $F^{\prime},F^{\prime
\prime}\in\overset{\bullet}{\mathcal{F}}\left(  K\right)  $.

\begin{definition}
Let $k$ be odd and an assignment function $\operatorname*{mark}%
\nolimits_{\mathcal{T}}:\mathcal{S}\rightarrow\mathcal{T}$ be chosen which
satisfies (\ref{mark}). The \emph{degrees of freedom for the Crouzeix-Raviart
spaces} are given

\begin{enumerate}
\item for $\operatorname*{CR}_{k}\left(  \mathcal{T}\right)  $ by%
%TCIMACRO{\TeXButton{doffull}{\begin{subequations}
%\label{doffull}
%\end{subequations}}}%
%BeginExpansion
\begin{subequations}
\label{doffull}
\end{subequations}%
%EndExpansion%
\begin{align}
\mathcal{J}_{k}  &  :=\left\{  J_{\tau,%
%TCIMACRO{\TeXButton{boldalphanew}{\boldsymbol{\alpha}}}%
%BeginExpansion
\boldsymbol{\alpha}%
%EndExpansion
}^{K}\mid\forall\left(  \tau,%
%TCIMACRO{\TeXButton{boldalphanew}{\boldsymbol{\alpha}}}%
%BeginExpansion
\boldsymbol{\alpha}%
%EndExpansion
\right)  \in\mathcal{I}_{k,\Omega}\quad\text{with }\operatorname*{mark}%
\nolimits_{\mathcal{T}}\left(  \tau\right)  =K\in\mathcal{T}\right\} \tag{%
%TCIMACRO{\TeXButton{doffull}{\ref{doffull}}}%
%BeginExpansion
\ref{doffull}%
%EndExpansion
a}\label{doffulla}\\
&  \cup\left\{  J_{\tau,%
%TCIMACRO{\TeXButton{boldalphanew}{\boldsymbol{\alpha}}}%
%BeginExpansion
\boldsymbol{\alpha}%
%EndExpansion
}^{F}:\forall\left(  \tau,%
%TCIMACRO{\TeXButton{boldalphanew}{\boldsymbol{\alpha}}}%
%BeginExpansion
\boldsymbol{\alpha}%
%EndExpansion
\right)  \in\mathcal{I}_{k,\partial\Omega}\quad\text{with }%
\operatorname*{mark}\nolimits_{\mathcal{T}}\left(  \tau\right)  =F\in
\mathcal{F}_{\partial\Omega}\right\} \tag{%
%TCIMACRO{\TeXButton{doffull}{\ref{doffull}}}%
%BeginExpansion
\ref{doffull}%
%EndExpansion
b}\label{doffullb}\\
&  \cup\left\{  J_{F}^{\operatorname*{CR},K}\mid\forall F\in\mathcal{F}%
_{\Omega}\quad\text{with }\operatorname*{mark}\nolimits_{\mathcal{T}}\left(
F\right)  =K\in\mathcal{T}\right\} \tag{%
%TCIMACRO{\TeXButton{doffull}{\ref{doffull}}}%
%BeginExpansion
\ref{doffull}%
%EndExpansion
c}\label{doffullc}\\
&  \cup\left\{  J_{F}^{\operatorname*{CR},F}:\forall F\in\mathcal{F}%
_{\partial\Omega}\right\}  \tag{%
%TCIMACRO{\TeXButton{doffull}{\ref{doffull}}}%
%BeginExpansion
\ref{doffull}%
%EndExpansion
d}\label{doffulld}%
\end{align}

\item and for $\operatorname*{CR}_{k,0}\left(  \mathcal{T}\right)  $ by%
%TCIMACRO{\TeXButton{dofzero}{\begin{subequations}
%\label{dofzero}
%\end{subequations}}}%
%BeginExpansion
\begin{subequations}
\label{dofzero}
\end{subequations}%
%EndExpansion%
\begin{align}
\mathcal{J}_{k,0}  &  :=\left\{  J_{\tau,%
%TCIMACRO{\TeXButton{boldalphanew}{\boldsymbol{\alpha}}}%
%BeginExpansion
\boldsymbol{\alpha}%
%EndExpansion
}^{K}\mid\forall\left(  \tau,%
%TCIMACRO{\TeXButton{boldalphanew}{\boldsymbol{\alpha}}}%
%BeginExpansion
\boldsymbol{\alpha}%
%EndExpansion
\right)  \in\mathcal{I}_{k,\Omega}\quad\text{with }\operatorname*{mark}%
\nolimits_{\mathcal{T}}\left(  \tau\right)  =K\in\mathcal{T}\right\} \tag{%
%TCIMACRO{\TeXButton{dofzero}{\ref{dofzero}}}%
%BeginExpansion
\ref{dofzero}%
%EndExpansion
a}\label{dofzeroa}\\
&  \cup\left\{  J_{F}^{\operatorname*{CR},K}\mid\forall F\in\mathcal{F}%
_{\Omega}\quad\text{with }\operatorname*{mark}\nolimits_{\mathcal{T}}\left(
F\right)  =K\in\mathcal{T}\right\}  . \tag{%
%TCIMACRO{\TeXButton{dofzero}{\ref{dofzero}}}%
%BeginExpansion
\ref{dofzero}%
%EndExpansion
b}\label{dofzerob}%
\end{align}

\end{enumerate}
\end{definition}

\begin{theorem}
Let $k$ be odd and assume that $\mathcal{T}$ contains more than one simplex.
Then, the degrees of freedom (\ref{doffull}) are bidual to the basis in
(\ref{basisfull}) and the set of freedoms (\ref{dofzero}) are bidual to set of
basis function (\ref{basiszero}).
\end{theorem}

%

%TCIMACRO{\TeXButton{Proof}{\proof}}%
%BeginExpansion
\proof
%EndExpansion
We prove the statement only for the space $\operatorname*{CR}_{k}\left(
\mathcal{T}\right)  $ more concretely we show first that for any
$u\in\operatorname*{CR}_{k}\left(  \mathcal{T}\right)  $ the implication%
\begin{equation}
\left(  J\left(  u\right)  =0\quad\forall J\in\mathcal{J}_{k}\right)
\implies\left(  u=0\right)  \label{Ju}%
\end{equation}
holds. Let $J\left(  u\right)  =0$ for all $J\in\mathcal{J}_{k}$. Then,
\[
\left.  u\right\vert _{K}=\sum_{\left(  t,%
%TCIMACRO{\TeXButton{boldbeta}{\boldsymbol{\beta}}}%
%BeginExpansion
\boldsymbol{\beta}%
%EndExpansion
\right)  \in\mathcal{I}_{k}\left(  K\right)  }\gamma_{\left(  t,%
%TCIMACRO{\TeXButton{boldbeta}{\boldsymbol{\beta}}}%
%BeginExpansion
\boldsymbol{\beta}%
%EndExpansion
\right)  }^{K}\left.  B_{t,%
%TCIMACRO{\TeXButton{boldbeta}{\boldsymbol{\beta}}}%
%BeginExpansion
\boldsymbol{\beta}%
%EndExpansion
}\right\vert _{K}+\sum_{F^{\prime}\in\mathcal{F}\left(  K\right)  }%
\delta_{F^{\prime}}\left.  B_{k}^{\operatorname*{CR},F^{\prime}}\right\vert
_{K}.
\]
The biduality property of $J_{\tau,%
%TCIMACRO{\TeXButton{boldalphanew}{\boldsymbol{\alpha}}}%
%BeginExpansion
\boldsymbol{\alpha}%
%EndExpansion
}^{K}$ imply%
\[
\gamma_{\left(  \tau,%
%TCIMACRO{\TeXButton{boldalphanew}{\boldsymbol{\alpha}}}%
%BeginExpansion
\boldsymbol{\alpha}%
%EndExpansion
\right)  }^{K}=0\quad\text{and\quad}\delta_{F}=0
\]
for all $\left(  \tau,%
%TCIMACRO{\TeXButton{boldalphanew}{\boldsymbol{\alpha}}}%
%BeginExpansion
\boldsymbol{\alpha}%
%EndExpansion
\right)  \in\mathcal{I}_{k,\Omega}$ with $K=\operatorname*{mark}%
\nolimits_{\mathcal{T}}\left(  \tau\right)  $ and $F\in\mathcal{F}_{\Omega}$
with $K=\operatorname*{mark}_{\mathcal{T}}\left(  F\right)  $. By applying
this argument for all $K\in\mathcal{T}$ it follows that $u$ has the
representation%
\begin{equation}
u=\sum_{\left(  t,%
%TCIMACRO{\TeXButton{boldbeta}{\boldsymbol{\beta}}}%
%BeginExpansion
\boldsymbol{\beta}%
%EndExpansion
\right)  \in\mathcal{I}_{k,\partial\Omega}}\gamma_{t,%
%TCIMACRO{\TeXButton{boldbeta}{\boldsymbol{\beta}}}%
%BeginExpansion
\boldsymbol{\beta}%
%EndExpansion
}^{F}B_{t,%
%TCIMACRO{\TeXButton{boldbeta}{\boldsymbol{\beta}}}%
%BeginExpansion
\boldsymbol{\beta}%
%EndExpansion
}^{F}+\sum_{F\in\mathcal{F}_{\partial\Omega}}\gamma_{F}B_{F}%
^{\operatorname*{CR},F}. \label{restu}%
\end{equation}
Now, let $F\in\mathcal{F}_{\partial\Omega}$ and $K\in\mathcal{T}$ be the
adjacent simplex. As stated after (\ref{deffpunkt}), $F\in\overset{\bullet
}{\mathcal{F}}\left(  K\right)  $ and $\left\vert \overset{\bullet
}{\mathcal{F}}\left(  K\right)  \right\vert <\left\vert \mathcal{F}\left(
K\right)  \right\vert $ so that Lemma \ref{LemRestrfacet} is applicable and we
can use the bidual properties (\ref{bidualpropsfacets}) and%
\[
\left.  u\right\vert _{F}=\sum_{\left(  t,%
%TCIMACRO{\TeXButton{boldbeta}{\boldsymbol{\beta}}}%
%BeginExpansion
\boldsymbol{\beta}%
%EndExpansion
\right)  \in\mathcal{I}_{k}\left(  F\right)  }\gamma_{t,%
%TCIMACRO{\TeXButton{boldbeta}{\boldsymbol{\beta}}}%
%BeginExpansion
\boldsymbol{\beta}%
%EndExpansion
}^{F}B_{t,%
%TCIMACRO{\TeXButton{boldbeta}{\boldsymbol{\beta}}}%
%BeginExpansion
\boldsymbol{\beta}%
%EndExpansion
}^{F}+\sum_{F^{\prime}\in\overset{\bullet}{\mathcal{F}}\left(  K\right)
\cap\mathcal{F}_{\partial\Omega}}\gamma_{F^{\prime}}B_{F^{\prime}%
}^{\operatorname*{CR},F^{\prime}}%
\]
to see that all the coefficients in (\ref{restu}) are zero.

Since the number of basis functions and the number of degrees of freedom
assigned to the same geometric entity $\tau\in\mathcal{S}$ are the same, it
follows $\left\vert \mathcal{J}_{k}\right\vert =\left\vert \mathcal{B}%
_{k}^{\operatorname*{CR}}\left(  \mathcal{T}\right)  \right\vert $ and
$\mathcal{J}_{k}$ is well-defined as set of degrees of freedom for
$\operatorname*{CR}_{k}\left(  \mathcal{T}\right)  $. The biduality follows by
construction.%
%TCIMACRO{\TeXButton{End Proof}{\endproof}}%
%BeginExpansion
\endproof
%EndExpansion

The dual basis directly gives rise to the definition of a
(quasi-)interpolation operator.

\begin{definition}
Let $k$ be odd and assume that $\mathcal{T}$ contains more than one simplex. A
\emph{quasi-interpolation operator }$I_{\mathcal{T},k}:H^{1}\left(
\Omega\right)  +\operatorname*{CR}_{k}\left(  \mathcal{T}\right)
\rightarrow\operatorname*{CR}_{k}\left(  \mathcal{T}\right)  $ is given by%
\begin{align*}
I_{\mathcal{T},k}u  &  =\sum_{\left(  \tau,%
%TCIMACRO{\TeXButton{boldalphanew}{\boldsymbol{\alpha}}}%
%BeginExpansion
\boldsymbol{\alpha}%
%EndExpansion
\right)  \in\mathcal{I}_{k,\Omega}}\left(  J_{\tau,%
%TCIMACRO{\TeXButton{boldalphanew}{\boldsymbol{\alpha}}}%
%BeginExpansion
\boldsymbol{\alpha}%
%EndExpansion
}^{K}u\right)  B_{\tau,%
%TCIMACRO{\TeXButton{boldalphanew}{\boldsymbol{\alpha}}}%
%BeginExpansion
\boldsymbol{\alpha}%
%EndExpansion
}+\sum_{F\in\mathcal{F}_{\Omega}}\left(  J_{F}^{\operatorname*{CR},K}u\right)
B_{k}^{\operatorname*{CR},F}\\
&  +\sum_{\left(  \tau,%
%TCIMACRO{\TeXButton{boldalphanew}{\boldsymbol{\alpha}}}%
%BeginExpansion
\boldsymbol{\alpha}%
%EndExpansion
\right)  \in\mathcal{I}_{k,\partial\Omega}}\left(  J_{\tau,%
%TCIMACRO{\TeXButton{boldalphanew}{\boldsymbol{\alpha}}}%
%BeginExpansion
\boldsymbol{\alpha}%
%EndExpansion
}^{F}u\right)  B_{\tau,%
%TCIMACRO{\TeXButton{boldalphanew}{\boldsymbol{\alpha}}}%
%BeginExpansion
\boldsymbol{\alpha}%
%EndExpansion
}+\sum_{F\in\mathcal{F}_{\partial\Omega}}\left(  J_{F}^{\operatorname*{CR}%
,F}u\right)  B_{k}^{\operatorname*{CR},F}%
\end{align*}
with the \emph{convention }that $K=\operatorname*{mark}_{\mathcal{T}}\left(
\tau\right)  $ for $\tau\in\mathcal{F}_{\Omega}$ and $F=\operatorname*{mark}%
_{\mathcal{T}}\left(  \tau\right)  $ for $\tau\in\mathcal{F}_{\partial\Omega}$.
\end{definition}

\begin{lemma}
Let $k$ be odd and assume that $\mathcal{T}$ contains more than one simplex.
The quasi-interpolation operator $I_{\mathcal{T},k}$ is a projection and the
restriction to $H_{0}^{1}\left(  \Omega\right)  +\operatorname*{CR}%
\nolimits_{k,0}\left(  \mathcal{T}\right)  $ satisfies
\begin{equation}
I_{\mathcal{T},k}:H_{0}^{1}\left(  \Omega\right)  +\operatorname*{CR}%
\nolimits_{k,0}\left(  \mathcal{T}\right)  \rightarrow\operatorname*{CR}%
\nolimits_{k,0}\left(  \mathcal{T}\right)  . \label{restrprop}%
\end{equation}
For $u\in H_{0}^{1}\left(  \Omega\right)  +\operatorname*{CR}\nolimits_{k,0}%
\left(  \mathcal{T}\right)  $ it holds%
\[
I_{\mathcal{T},k}u=\sum_{\left(  \tau,%
%TCIMACRO{\TeXButton{boldalphanew}{\boldsymbol{\alpha}}}%
%BeginExpansion
\boldsymbol{\alpha}%
%EndExpansion
\right)  \in\mathcal{I}_{k,\Omega}}\left(  J_{\tau,%
%TCIMACRO{\TeXButton{boldalphanew}{\boldsymbol{\alpha}}}%
%BeginExpansion
\boldsymbol{\alpha}%
%EndExpansion
}^{K}u\right)  B_{\tau,%
%TCIMACRO{\TeXButton{boldalphanew}{\boldsymbol{\alpha}}}%
%BeginExpansion
\boldsymbol{\alpha}%
%EndExpansion
}+\sum_{F\in\mathcal{F}_{\Omega}}\left(  J_{F}^{\operatorname*{CR},K}u\right)
B_{k}^{\operatorname*{CR},F}.
\]

\end{lemma}%

%TCIMACRO{\TeXButton{Proof}{\proof}}%
%BeginExpansion
\proof
%EndExpansion
The projection property of $I_{\mathcal{T},k}$ is a direct consequence of the
biduality of the degrees of freedom. The property (\ref{restrprop}) of the
restriction simply follows since all functionals related to geometric entities
on the boundary $\partial\Omega$ are defined via integrals over facets on the
boundary and hence vanish for functions in $H_{0}^{1}\left(  \Omega\right)  $.
The boundary facet functionals $J_{\tau,%
%TCIMACRO{\TeXButton{boldalphanew}{\boldsymbol{\alpha}}}%
%BeginExpansion
\boldsymbol{\alpha}%
%EndExpansion
}^{F}$, $\left(  \tau,%
%TCIMACRO{\TeXButton{boldalphanew}{\boldsymbol{\alpha}}}%
%BeginExpansion
\boldsymbol{\alpha}%
%EndExpansion
\right)  \in\mathcal{I}_{k,\partial\Omega}$, and $J_{F}^{\operatorname*{CR}%
,F}$, $F\in\mathcal{F}_{\partial\Omega}$, vanish for functions $u\in
\operatorname*{CR}_{k,0}\left(  \Omega\right)  $ due to their biduality.%
%TCIMACRO{\TeXButton{End Proof}{\endproof}}%
%BeginExpansion
\endproof
%EndExpansion

The following remark explains why, in general, local degrees of freedom
\textit{do not exist for even polynomial degree} $k$.

\begin{remark}
\label{Remkeven}Let $k$ be even. From Lemma \ref{LemBasis} it follows that a
basis for $\operatorname*{CR}_{k}\left(  \mathcal{T}\right)  $ is given by%
\[
\mathcal{B}_{k}^{\operatorname*{CR}}\left(  \mathcal{T}\right)  =\left\{
B_{k}^{\mathbf{z}}\mid\forall\mathbf{z}\in\mathcal{N}_{k}\left(
\mathcal{T}\right)  \right\}  \cup\left\{  B_{k}^{\operatorname*{CR},K}%
\mid\forall K\in\mathcal{T}^{\prime}\right\}  ,
\]
where $\mathcal{T}^{\prime}\subset\mathcal{T}$ is any submesh with $\left\vert
\mathcal{T}^{\prime}\right\vert =\left\vert \mathcal{T}\right\vert -1$. Let
$\widetilde{\mathcal{T}}\subset\mathcal{T}^{\prime}$ and set $\widetilde
{\Omega}:=%
%TCIMACRO{\dbigcup \limits_{K\in\widetilde{\mathcal{T}}}}%
%BeginExpansion
{\displaystyle\bigcup\limits_{K\in\widetilde{\mathcal{T}}}}
%EndExpansion
K$. Let
\[
\widetilde{\mathcal{B}_{k}^{\operatorname*{CR}}}=\left\{  B_{k}^{\mathbf{z}%
}\mid\forall\mathbf{z}\in\mathcal{N}_{k}\left(  \mathcal{T}\right)  \text{
with }\left.  B_{k}^{\mathbf{z}}\right\vert _{\widetilde{\Omega}}%
\neq0\right\}  \cup\left\{  B_{k}^{\operatorname*{CR},K}\mid\forall
K\in\widetilde{\mathcal{T}}\right\}
\]
be the set of basis functions whose restrictions to $\widetilde{\Omega}$ are
non-trivial. By contradiction, assume that there is a set $\widetilde
{\mathcal{J}_{k}^{\operatorname*{CR}}}$ of basis functions for the dual space
$\left(  \operatorname*{span}\widetilde{\mathcal{B}_{k}^{\operatorname*{CR}}%
}\right)  ^{\prime}$ \emph{with local support property}: $\operatorname*{supp}%
J\subset\widetilde{\Omega}$ for all $J\in\widetilde{\mathcal{J}_{k}%
^{\operatorname*{CR}}}$. Then the corresponding Gram's matrix $\left(
J\left(  B\right)  \right)  _{\substack{J\in\widetilde{\mathcal{J}%
_{k}^{\operatorname*{CR}}}\\B\in\widetilde{\mathcal{B}_{k}^{\operatorname*{CR}%
}}}}$ is regular. However, as in the proof (\textbf{Case 2, Part 1},
\textbf{1st step) }of Lemma \ref{LemBasis} it follows that%
\[
\left.  S_{k}\left(  \mathcal{T}\right)  \right\vert _{\widetilde{\Omega}}%
\cap\operatorname*{span}\left\{  B_{k}^{\operatorname*{CR},K}\mid\forall
K\in\widetilde{\mathcal{T}}\right\}  =\operatorname*{span}\left\{  \left.
\Psi_{k}\right\vert _{\widetilde{\Omega}}\right\}
\]
with $\Psi_{k}$ as in (\ref{Vknccup}). This implies that for the space
$\operatorname*{CR}_{k}\left(  \mathcal{T}\right)  $ and even $k$ there exists
\emph{no} basis of the dual space which has the local support property.%
%TCIMACRO{\TeXButton{End Proof}{\endproof}}%
%BeginExpansion
\endproof
%EndExpansion

\end{remark}

\section{Degrees of freedom and an approximation operator for the case $d=2$
and $k$ is odd\label{SecDofs}}

The degrees of freedom associated to the inner faces for general dimension
$d\geq2$ are defined by integrals over adjacent volume simplices
$K\in\mathcal{T}$. For the important case of two spatial dimensions $d=2$,
these degrees of freedom can be defined alternatively by facet (edge)
integrals with weight functions which are polynomials of maximal degree $k-1$.
In this section we present their construction.

\subsection{Degrees of freedom\label{Dofd2}}

For the following construction we assume that $k\geq1$ is odd and that
$\mathcal{T}$ is a triangulation of a two-dimensional polygonal domain. Recall
the $\Sigma_{1}$ denotes the mesh skeleton (the union of all edges). We will
need the space%
\begin{equation}
\overset{\bullet\bullet}{S_{k,0}}\left(  \mathcal{T}\right)  :=\left\{
u\in\overset{\bullet}{S_{k,0}}\left(  \mathcal{T}\right)  \mid\left.
u\right\vert _{\Sigma_{1}}=0\right\}  \label{DefSkdotdot}%
\end{equation}
which is non-trivial for $k\geq3$ and then consists of polynomial triangle bubbles.

We present a partial basis along degrees of freedom for $\operatorname*{CR}%
\nolimits_{k,0}\left(  \mathcal{T}\right)  $ which consists of conforming and
non-conforming edge functions and the direct sum of their spans with
$\overset{\bullet\bullet}{S_{k,0}}\left(  \mathcal{T}\right)  $ yields
$\operatorname*{CR}\nolimits_{k,0}\left(  \mathcal{T}\right)  $.

This needs some notation. For $d=2$ the set of edges and the set of facets
coincide and we write $\mathcal{E}_{\Omega}$ for $\mathcal{F}_{\Omega}$ (cf.
(\ref{Defboundinfacets})). For $E\in\mathcal{E}_{\Omega}$, we fix a numbering
of the endpoints of $E$: $\mathbf{A}_{1}\left(  E\right)  ,\mathbf{A}%
_{2}\left(  E\right)  $ and define the edge bubble by $b_{E}=\varphi
_{\mathbf{A}_{1}\left(  E\right)  }^{1}\varphi_{\mathbf{A}_{2}\left(
E\right)  }^{1}$ which belongs to $S_{2,0}\left(  \mathcal{T}\right)  $.

The edge basis functions are given by%
\begin{equation}%
\begin{array}
[c]{ll}%
\left.
\begin{tabular}
[c]{l}%
$\forall E\in\mathcal{E}_{\Omega}$\\
$\forall\mu\in\left\{  0,1,\ldots,k-1\right\}  $%
\end{tabular}
\ \ \ \right\}  & B_{\mu}^{E}:=\left\{
\begin{array}
[c]{lc}%
b_{E}P_{\mu}^{\left(  1,1\right)  }\left(  2\varphi_{\mathbf{A}_{2}\left(
E\right)  }^{1}-1\right)  & \mu\leq k-2,\\
B_{k}^{\operatorname*{CR},E} & \mu=k-1.
\end{array}
\right.
\end{array}
\label{edgebasis}%
\end{equation}

\begin{corollary}
\label{CorProp}The basis functions have the following support properties%
\[
\operatorname*{supp}B_{\mu}^{E}\subset\mathcal{T}_{E}\qquad\forall\mu
\in\left\{  0,1,\ldots,k-1\right\}  .
\]

\end{corollary}

Next we define the degrees of freedom for these functions. Since the jumps of
$v\in\operatorname*{CR}_{k,0}\left(  \mathcal{T}\right)  $ satisfy moment
conditions up to an order $k-1$ across each edge, edge functionals of the type
$J_{\nu}^{E}\left(  u\right)  :=\int_{E}g_{\nu}^{E}u$ defined via polynomials
$g_{\nu}^{E}\in\mathbb{P}_{k-1}\left(  E\right)  $ are well-defined for
Crouzeix-Raviart functions.

\begin{definition}
\label{DefEdgeFunc}The edge functionals $J_{\nu}^{E}$ are given for any
$E\in\mathcal{E}_{\Omega}$ and $\nu\in\left\{  1,2,\ldots,k-1\right\}  $ by%
%TCIMACRO{\TeXButton{dofs2d}{\begin{subequations}
%\label{dofs2d}
%\end{subequations}}}%
%BeginExpansion
\begin{subequations}
\label{dofs2d}
\end{subequations}%
%EndExpansion%
\begin{equation}
J_{\nu}^{E}u:=\frac{2}{\left\vert E\right\vert }\int_{E}g_{\nu}^{E}u \tag{%
%TCIMACRO{\TeXButton{dofs2d}{\ref{dofs2d}}}%
%BeginExpansion
\ref{dofs2d}%
%EndExpansion
a}\label{dofs2da}%
\end{equation}
with weight function%
\begin{equation}
g_{\nu}^{E}:=\gamma_{\nu}\left(  P_{\nu}^{\left(  1,1\right)  }\left(
2\varphi_{\mathbf{A}_{2}\left(  E\right)  }-1\right)  -c_{\nu,k}%
P_{k-1}^{\left(  1,1\right)  }\left(  2\varphi_{\mathbf{A}_{2}\left(
E\right)  }-1\right)  \right)  \tag{%
%TCIMACRO{\TeXButton{dofs2d}{\ref{dofs2d}}}%
%BeginExpansion
\ref{dofs2d}%
%EndExpansion
b}\label{dofs2db}%
\end{equation}
and constants\footnote{The constant $\gamma_{\nu}$ takes into account the
relation (cf. \cite[Table 18.3.1]{NIST:DLMF}):%
\[
\int_{-1}^{1}\left(  1-x\right)  \left(  1+x\right)  \left(  P_{\nu}^{\left(
1,1\right)  }\left(  x\right)  \right)  ^{2}dx=\gamma_{\nu}^{-1}.
\]
}%
\begin{equation}
\gamma_{\nu}:=\frac{\left(  2\nu+3\right)  \left(  \nu+2\right)  }{8\left(
\nu+1\right)  }\text{\quad and\ }c_{\nu,k}:=\left\{
\begin{array}
[c]{ll}%
\frac{1-\left(  -1\right)  ^{\nu+1}}{2}\frac{k+1}{\nu+2} & \nu\in\left\{
0,1,\ldots,k-2\right\}  ,\\
\frac{1}{2k+1} & \nu=k-1.
\end{array}
\right.  \tag{%
%TCIMACRO{\TeXButton{dofs2d}{\ref{dofs2d}}}%
%BeginExpansion
\ref{dofs2d}%
%EndExpansion
c}\label{Defcvk}%
\end{equation}

\end{definition}

\begin{lemma}
Let $k\geq1$ be odd. For $\nu,\mu\in\left\{  0,1,\ldots,k-1\right\}  $ and
$E,E^{\prime}\in\mathcal{E}_{\Omega}$ and $u\in\overset{\bullet\bullet
}{S_{k,0}}\left(  \mathcal{T}\right)  $ it holds%
\begin{equation}
J_{\nu}^{E}\left(  B_{\mu}^{E^{\prime}}\right)  =\delta_{\nu,\mu}%
\delta_{E,E^{\prime}},\quad J_{\nu}^{E}\left(  u\right)  =0. \label{funcE1}%
\end{equation}

\end{lemma}%

%TCIMACRO{\TeXButton{Proof}{\proof}}%
%BeginExpansion
\proof
%EndExpansion
First we consider the case $\mu\in\left\{  0,1,\ldots,k-2\right\}  $.

The support properties (Corollary \ref{CorProp}) imply that for $E^{\prime}%
\in\mathcal{E}_{\Omega}\backslash\left\{  E\right\}  $ it holds $\left.
B_{\mu}^{E^{\prime}}\right\vert _{E}=0$ and hence $J_{\nu}^{E}\left(  B_{\mu
}^{E^{\prime}}\right)  =0$ for such $E^{\prime}$. Next, let $E=E^{\prime}$. We
choose an affine bijection $\varphi_{E}:\left[  -1,1\right]  \rightarrow E$
such that $\varphi_{E}\left(  -1\right)  =\mathbf{A}_{1}\left(  E\right)  $.
Then%
\[
\left(  \left.  \left(  2\varphi_{\mathbf{A}_{2}\left(  E\right)  }-1\right)
\right\vert _{E}\circ\varphi_{E}\right)  \left(  x\right)  =x.
\]
For $\mu\in\left\{  0,1,\ldots,k-2\right\}  $ we get%
\[
J_{\nu}^{E}\left(  B_{\mu}^{E}\right)  =\gamma_{\nu}\int_{-1}^{1}\left(
1-x\right)  \left(  1+x\right)  \left(  P_{\nu}^{\left(  1,1\right)  }\left(
x\right)  -c_{\nu}P_{k-1}^{\left(  1,1\right)  }\left(  x\right)  \right)
P_{\mu}^{\left(  1,1\right)  }\left(  x\right)  dx=\delta_{\nu,\mu}\text{.}%
\]

It remains to consider the case $\mu=k-1$. Note that the function $\left.
B_{k-1}^{E}\right\vert _{E^{\prime}}$ is zero for all edges $E^{\prime}%
\subset\mathcal{E}_{\Omega}$ which are not subset of the edge patch
$\omega_{E}$ and hence $J_{\nu}^{E^{\prime}}\left(  B_{k-1}^{E}\right)  =0$
for those edges.

For $E^{\prime}=E$, it holds%
\begin{align*}
J_{\nu}^{E}\left(  B_{k-1}^{E}\right)   &  =\frac{2}{\left\vert E\right\vert
}\int_{E}\gamma_{\nu}\left(  P_{\nu}^{\left(  1,1\right)  }\left(
2\varphi_{\mathbf{A}_{2}\left(  E\right)  }-1\right)  -c_{\nu,k}%
P_{k-1}^{\left(  1,1\right)  }\left(  2\varphi_{\mathbf{A}_{2}\left(
E\right)  }-1\right)  \right)  P_{k}^{\left(  0,0\right)  }\left(  1\right) \\
&  =\int_{-1}^{1}\gamma_{\nu}\left(  P_{\nu}^{\left(  1,1\right)  }\left(
x\right)  -c_{\nu,k}P_{k-1}^{\left(  1,1\right)  }\left(  x\right)  \right)
dx.
\end{align*}
From \cite[18.9.15]{NIST:DLMF} it follows that $\frac{2}{\nu+2}\left(
P_{\nu+1}^{\left(  0,0\right)  }\right)  ^{\prime}=P_{\nu}^{\left(
1,1\right)  }$ so that
\begin{align*}
\int_{-1}^{1}P_{\nu}^{\left(  1,1\right)  }\left(  x\right)  dx  &  =\frac
{2}{\nu+2}\int_{-1}^{1}\left(  P_{\nu+1}^{\left(  0,0\right)  }\right)
^{\prime}=\frac{2\left(  1-\left(  -1\right)  ^{\nu+1}\right)  }{\nu+2},\\
\int_{-1}^{1}P_{k-1}^{\left(  1,1\right)  }\left(  x\right)  dx  &
\overset{k\text{ odd}}{=}\frac{4}{k+1}.
\end{align*}%
%TCIMACRO{\TeXButton{black}{\color{black}}}%
%BeginExpansion
\color{black}%
%EndExpansion
The coefficient $c_{\nu,k}$ in (\ref{Defcvk}) is chosen such that
\[
J_{\nu}^{E}\left(  B_{k-1}^{E}\right)  =\gamma_{\nu}\left(  \frac{2\left(
1-\left(  -1\right)  ^{\nu+1}\right)  }{\nu+2}-c_{\nu,k}\frac{4}{k+1}\right)
=\delta_{k-1,\nu}\qquad\forall\nu\in\left\{  0,1,\ldots,k-1\right\}  .
\]
It remains to consider those edges $E^{\prime}\in\mathcal{E}_{\Omega}$ with
$E^{\prime}\subset\partial\omega_{E}$. Then,%
\[
J_{\nu}^{E^{\prime}}\left(  B_{k-1}^{E}\right)  =\frac{2}{\left\vert
E\right\vert }\int_{E}\left(  P_{\nu}^{\left(  1,1\right)  }\left(
2\varphi_{\mathbf{A}_{2}\left(  E^{\prime}\right)  }-1\right)  -c_{\nu
,k}P_{k-1}^{\left(  1,1\right)  }\left(  2\varphi_{\mathbf{A}_{2}\left(
E^{\prime}\right)  }-1\right)  \right)  P_{k}^{\left(  0,0\right)  }\left(
1-2\lambda_{K,E}\right)  .
\]
Since $\left.  P_{k}^{\left(  0,0\right)  }\left(  1-2\lambda_{K,E}\right)
\right\vert _{E^{\prime}}$ is the (lifted) Legendre polynomial of degree $k$,
the orthogonality properties of Legendre polynomials imply that this integral
is zero.

The second statement in (\ref{funcE1}) simply follows from $\left.
u\right\vert _{E}=0$ for $u\in\overset{\bullet\bullet}{S_{k,0}}\left(
\mathcal{T}\right)  $ and any edge $E\in\mathcal{E}$.%
%TCIMACRO{\TeXButton{End Proof}{\endproof}}%
%BeginExpansion
\endproof
%EndExpansion

The linearly independent functions $B_{\mu}^{E}$ are collected in the set%
\[
\mathcal{B}_{\mathcal{E}}\left(  \mathcal{T}\right)  :=\left\{  B_{\mu}%
^{E}:\forall E\in\mathcal{E}_{\Omega}\quad\forall\mu\in\left\{  0,1,\ldots
,k-1\right\}  \right\}
\]
and span the interelement space%
\[
S_{k,0}^{\mathcal{E}}\left(  \mathcal{T}\right)  :=\operatorname*{span}%
\mathcal{B}_{\mathcal{E}}\left(  \mathcal{T}\right)  .
\]

\begin{corollary}
Let $k\geq1$ be odd. Then%
\[
\operatorname*{CR}\nolimits_{k,0}\left(  \mathcal{T}\right)  =\overset
{\bullet\bullet}{S_{k,0}}\left(  \mathcal{T}\right)  \oplus S_{k,0}%
^{\mathcal{E}}\left(  \mathcal{T}\right)  .
\]

\end{corollary}

\subsection{An approximation operator}

In this section, we will introduce a class of local approximation operators
$I_{\mathcal{T},k}:H_{0}^{1}\left(  \Omega\right)  +\operatorname*{CR}%
_{k,0}\left(  \mathcal{T}\right)  \rightarrow\operatorname*{CR}_{k,0}\left(
\mathcal{T}\right)  $. The construction starts by defining the edge related
interpolation%
\begin{equation}
\mathcal{I}_{\mathcal{T},k}^{\mathcal{E}}u:=\sum_{E\in\mathcal{E}_{\Omega}%
}\sum_{\nu=0}^{k-1}\left(  J_{\nu}^{E}u\right)  B_{\nu}^{E}.
\label{defedgeinterpol}%
\end{equation}
Next, one chooses a bounded linear operator $\overset{\bullet\bullet}{\Pi
}_{\mathcal{T},k}:H_{0}^{1}\left(  \Omega\right)  +\operatorname*{CR}%
_{k,0}\left(  \mathcal{T}\right)  \rightarrow\overset{\bullet\bullet}{S_{k,0}%
}\left(  \mathcal{T}\right)  $ which is local in the sense that for all
$K\in\mathcal{T}$%
\[
\left.  \left(  \overset{\bullet\bullet}{\Pi}_{\mathcal{T},k}u\right)
\right\vert _{K}=\Pi_{k,0}^{K}\left(  \left.  u\right\vert _{K}\right)
\quad\forall u\in H_{0}^{1}\left(  \Omega\right)  +\operatorname*{CR}%
\nolimits_{k,0}\left(  \mathcal{T}\right)  ,
\]
for given local projections $\Pi_{k,0}^{K}:L^{2}\left(  K\right)  \rightarrow
S_{k,0}\left(  K\right)  $. A typical choice of $\Pi_{k,0}^{K}$ is the
$L^{2}\left(  K\right)  $-orthogonal projection onto $S_{k,0}\left(  K\right)
$.

\begin{definition}
A class of \emph{local approximation operators} $I_{\mathcal{T},k}:H_{0}%
^{1}\left(  \Omega\right)  +\operatorname*{CR}_{k,0}\left(  \mathcal{T}%
\right)  \rightarrow\operatorname*{CR}_{k,0}\left(  \mathcal{T}\right)  $ is
given by%
\begin{equation}
I_{\mathcal{T},k}u:=\mathcal{I}_{\mathcal{T},k}^{\mathcal{E}}u+\overset
{\bullet\bullet}{\Pi}_{\mathcal{T},k}\left(  u-\mathcal{I}_{\mathcal{T}%
,k}^{\mathcal{E}}u\right)  . \label{defapproxop}%
\end{equation}

\end{definition}

\begin{remark}
A basis and bidual basis for the Crouzeix-Raviart space $\operatorname*{CR}%
_{k}\left(  \mathcal{T}\right)  $ for $d=2$ and odd $k$ is obtained by
defining basis functions and degrees of freedom for the boundary edges
analogously as in (\ref{edgebasis}), (\ref{DefEdgeFunc}) for inner edges. A
local approximation operator in $\operatorname*{CR}_{k}\left(  \mathcal{T}%
\right)  $ is then obtained by summing in (\ref{defedgeinterpol}) over all
$E\in\mathcal{E}$ and using the resulting edge interpolation operator in
(\ref{defapproxop}).
\end{remark}

\subsection{Non-existence of split facet/simplex degrees of freedom for
$d\geq3$ and $k\neq1$\label{SecNonExist}}

We have already explained in Remark \ref{Remkeven} that there exists no set of
local degrees of freedom for even $k$. In this section we will prove that a
construction as in Section \ref{Dofd2} is possible only for $k=1$ or $d=2$ and
$k$ odd.

First, we consider the case $k=1$.

\begin{lemma}
For $d\geq2$ and $k=1$ the Crouzeix-Raviart basis functions are given by (cf.
Definition \ref{DefNCShapeFct}, Remark \ref{Remk=1}):%
\[
B_{1}^{\operatorname*{CR},F}:=\left\{
\begin{array}
[c]{ll}%
P_{1}^{\left(  0,d-2\right)  }\left(  1-2\lambda_{K,F}\right)  & \text{for
}K\in\mathcal{T}_{F},\\
0 & \text{otherwise,}%
\end{array}
\right.
\]
and the degrees of freedom by%
\[
J_{0}^{F}\left(  u\right)  :=\frac{1}{\left\vert F\right\vert }\int_{F}u.
\]
They form a bidual basis:%
\[
J_{0}^{F}\left(  B_{1}^{\operatorname*{CR},F^{\prime}}\right)  =\delta
_{F,F^{\prime}}\quad\forall F,F^{\prime}\in\mathcal{E}\text{.}%
\]

\end{lemma}%

%TCIMACRO{\TeXButton{Proof}{\proof}}%
%BeginExpansion
\proof
%EndExpansion
It holds%
\[
J_{0}^{F}\left(  B_{1}^{\operatorname*{CR},F}\right)  =\frac{1}{\left\vert
F\right\vert }\int_{F}P_{1}^{\left(  0,d-2\right)  }\left(  1-2\lambda
_{K,F}\right)  =\frac{1}{\left\vert F\right\vert }\int_{F}1=1.
\]
For $F\in\mathcal{E}$, $K\in\mathcal{T}_{F}$, and $F^{\prime}\subset\partial
K\backslash\overset{\circ}{F}$ we get%
\[
J_{0}^{F^{\prime}}\left(  B_{1}^{\operatorname*{CR},F}\right)  =\frac
{1}{\left\vert F^{\prime}\right\vert }\int_{F^{\prime}}P_{1}^{\left(
0,d-2\right)  }\left(  1-2\lambda_{K,F}\right)  \overset
{\text{(\ref{orthofacetpropb})}}{=}0.
\]
If $F^{\prime}$ is outside $\omega_{F}$ the functional $J_{0}^{F^{\prime}%
}\left(  B_{1}^{\operatorname*{CR},F}\right)  $ vanishes since
$\operatorname*{supp}B_{1}^{\operatorname*{CR},F}=\omega_{F}$.%
%TCIMACRO{\TeXButton{End Proof}{\endproof}}%
%BeginExpansion
\endproof
%EndExpansion

Next, we consider the case that $d\geq3$ and $k>1$ is odd. Section
\ref{SecDofs} introduces a construction
%TCIMACRO{\TeXButton{black}{\color{black}}}%
%BeginExpansion
\color{black}%
%EndExpansion
for a bidual basis for $\operatorname{CR}_{k}\left(  \mathcal{T}\right)  $
($\operatorname{CR}_{k,0}\left(  \mathcal{T}\right)  $ resp.), where for all
$\left(  \tau,%
%TCIMACRO{\TeXButton{boldalphanew}{\boldsymbol{\alpha}}}%
%BeginExpansion
\boldsymbol{\alpha}%
%EndExpansion
\right)  \in\mathcal{I}_{k,\ell}$, $\ell\in\left\{  1,2,...,d-1\right\}  $,
the functionals $J_{\tau,%
%TCIMACRO{\TeXButton{boldalphanew}{\boldsymbol{\alpha}}}%
%BeginExpansion
\boldsymbol{\alpha}%
%EndExpansion
}$ are of the form
\[
J_{\tau,%
%TCIMACRO{\TeXButton{boldalphanew}{\boldsymbol{\alpha}}}%
%BeginExpansion
\boldsymbol{\alpha}%
%EndExpansion
}:=\left(  g_{\tau,%
%TCIMACRO{\TeXButton{boldalphanew}{\boldsymbol{\alpha}}}%
%BeginExpansion
\boldsymbol{\alpha}%
%EndExpansion
},\cdot\right)  _{L^{2}\left(  F\right)  }\qquad g_{\tau,%
%TCIMACRO{\TeXButton{boldalphanew}{\boldsymbol{\alpha}}}%
%BeginExpansion
\boldsymbol{\alpha}%
%EndExpansion
}\in\mathbb{P}_{k-1}\left(  F\right)  ,\ \tau\subseteq F.
\]
In this section we prove that this is not possible for $d\geq3$ and $k\geq3$
odd. Let the assignement function $\operatorname{mark}_{\mathcal{T}%
}:\mathcal{S}\rightarrow\mathcal{T}\cup\mathcal{F}$ satisfy
\begin{subequations}
\label{Eq:Properties markT}%
\begin{align}
\operatorname{mark}_{\mathcal{T}}\left(  K\right)   &  =K\qquad\forall
K\in\mathcal{T}\label{SubEq:Properties markT - simplex}\\
\operatorname{mark}_{\mathcal{T}}\left(  \tau\right)   &  \in\mathcal{F}%
\;\text{and}\;\tau\subseteq\operatorname{mark}_{\mathcal{T}}\left(
\tau\right)  \qquad\forall\tau\in\mathcal{S}\setminus\mathcal{T}.
\label{SubEq:Properties markT - faces}%
\end{align}
Further let $J_{\tau,%
%TCIMACRO{\TeXButton{boldalphanew}{\boldsymbol{\alpha}}}%
%BeginExpansion
\boldsymbol{\alpha}%
%EndExpansion
},J_{F}^{\operatorname{CR}}:H^{1}\left(  \mathcal{T}\right)  \rightarrow
\mathbb{R}$ be given by
\end{subequations}
\begin{subequations}
\label{Eq:Entity functionals}%
\begin{align}
J_{\tau,%
%TCIMACRO{\TeXButton{boldalphanew}{\boldsymbol{\alpha}}}%
%BeginExpansion
\boldsymbol{\alpha}%
%EndExpansion
}u  &  :=\left(  g_{\tau,%
%TCIMACRO{\TeXButton{boldalphanew}{\boldsymbol{\alpha}}}%
%BeginExpansion
\boldsymbol{\alpha}%
%EndExpansion
},u\right)  _{L^{2}\left(  \operatorname{mark}_{\mathcal{T}}\left(
\tau\right)  \right)  }\quad g_{\tau,%
%TCIMACRO{\TeXButton{boldalphanew}{\boldsymbol{\alpha}}}%
%BeginExpansion
\boldsymbol{\alpha}%
%EndExpansion
}\in%
\begin{cases}
L^{2}\left(  \operatorname{mark}_{\mathcal{T}}\left(  \tau\right)  \right)  &
\operatorname{mark}_{\mathcal{T}}\left(  \tau\right)  \in\mathcal{T},\\
\mathbb{P}_{k-1}\left(  \operatorname{mark}_{\mathcal{T}}\left(  \tau\right)
\right)  & \operatorname{mark}_{\mathcal{T}}\left(  \tau\right)
\in\mathcal{F},
\end{cases}
\label{SubEq:Entity functionals - conf}\\
J_{F}^{\operatorname{CR}}u  &  :=\left(  g_{F},u\right)  _{L^{2}\left(
F\right)  }\quad g_{F}\in\mathbb{P}_{k-1}\left(  F\right)  ,
\label{SubEq:Entity functionals - nonconf}%
\end{align}
and collect them in the sets
\end{subequations}
\begin{subequations}
\label{Eq:Proposed bidualbasis}%
\begin{align}
\mathcal{J}_{k}^{\operatorname{entity}}\left(  \mathcal{T}\right)   &
:=\left\{  \left.  J_{\tau,%
%TCIMACRO{\TeXButton{boldalphanew}{\boldsymbol{\alpha}}}%
%BeginExpansion
\boldsymbol{\alpha}%
%EndExpansion
}\;\right\vert \;\left(  \tau,%
%TCIMACRO{\TeXButton{boldalphanew}{\boldsymbol{\alpha}}}%
%BeginExpansion
\boldsymbol{\alpha}%
%EndExpansion
\right)  \in\mathcal{I}_{k}\right\}  \cup\left\{  \left.  J_{F}%
^{\operatorname{CR}}\;\right\vert \;F\in\mathcal{F}\right\}
,\label{SubEq:Proposed bidualbasis - full space}\\
\mathcal{J}_{k,0}^{\operatorname{entity}}\left(  \mathcal{T}\right)   &
:=\left\{  \left.  J_{\tau,%
%TCIMACRO{\TeXButton{boldalphanew}{\boldsymbol{\alpha}}}%
%BeginExpansion
\boldsymbol{\alpha}%
%EndExpansion
}\;\right\vert \;\left(  \tau,%
%TCIMACRO{\TeXButton{boldalphanew}{\boldsymbol{\alpha}}}%
%BeginExpansion
\boldsymbol{\alpha}%
%EndExpansion
\right)  \in\mathcal{I}_{k,\Omega}\right\}  \cup\left\{  \left.
J_{F}^{\operatorname{CR}}\;\right\vert \;F\in\mathcal{F}_{\Omega}\right\}  .
\label{SubEq:Proposed bidualbasis - zero bc}%
\end{align}

\end{subequations}
\begin{theorem}
Let $d\geq3$ and $k\geq3$ odd. If the assignment function $\operatorname{mark}%
_{\mathcal{T}}:\mathcal{S}\rightarrow\mathcal{T}\cup\mathcal{F}$ satisfies
(\ref{Eq:Properties markT}), then there exists \emph{no} set of functionals of
the form (\ref{Eq:Entity functionals}) such that

\begin{enumerate}
\item $\mathcal{J}_{k}^{\operatorname{entity}}\left(  \mathcal{T}\right)  $ is
a bidual basis for $\operatorname{CR}_{k}\left(  \mathcal{T}\right)  $.

\item $\mathcal{J}_{k,0}^{\operatorname{entity}}\left(  \mathcal{T}\right)  $
is a bidual basis for $\operatorname{CR}_{k,0}\left(  \mathcal{T}\right)  $ if
there exists $K\in\mathcal{T}$ such that $K\cap\partial\Omega=\emptyset$.
\end{enumerate}
\end{theorem}

\begin{proof}
\textbf{@1.} Let us assume by contradiction that there exists a choice of
functionals of the form \eqref{Eq:Entity functionals}, such that
$\mathcal{J}_{k}^{\operatorname{entity}}\left(  \mathcal{T}\right)  $ is a
bidual basis. First observe that, due to the inclusion
(\ref{Eq:Inclusion Sk in CRk}), the constant function $1\in\operatorname{CR}%
_{k}\left(  \mathcal{T}\right)  $. Hence for any functional $J\in
\mathcal{J}_{k}^{\operatorname{entity}}\left(  \mathcal{T}\right)  $,
$J\left(  1\right)  $ is well defined. Thus we choose any simplex
$K\in\mathcal{T}$ and define the set $\mathcal{I}_{k}\left(  \partial
K\right)  :=\bigcup_{\ell=1}^{d-1}\mathcal{I}_{k,\ell}\left(  K\right)  $ as
well as
\begin{align}
\psi_{\partial K}  &  :=\sum_{\left(  \tau,%
%TCIMACRO{\TeXButton{boldalphanew}{\boldsymbol{\alpha}}}%
%BeginExpansion
\boldsymbol{\alpha}%
%EndExpansion
\right)  \in\mathcal{I}_{k}\left(  \partial K\right)  }J_{\tau,%
%TCIMACRO{\TeXButton{boldalphanew}{\boldsymbol{\alpha}}}%
%BeginExpansion
\boldsymbol{\alpha}%
%EndExpansion
}\left(  1\right)  B_{\tau,%
%TCIMACRO{\TeXButton{boldalphanew}{\boldsymbol{\alpha}}}%
%BeginExpansion
\boldsymbol{\alpha}%
%EndExpansion
}+\sum_{F\in\mathcal{F}\left(  K\right)  }J_{F}^{\operatorname{CR}}\left(
1\right)  B_{k}^{\operatorname{CR},F},\nonumber\\
\phi_{\partial K}  &  :=\psi_{\partial K}-\sum_{F\in\mathcal{F}\left(
K\right)  }B_{k}^{\operatorname{CR},F}. \label{Eq:Def phiPartialK}%
\end{align}
Following the arguments of the proof of Lemma \ref{LemLinInd}, we deduce that
the functions $\left\{  \left.  B_{\tau,%
%TCIMACRO{\TeXButton{boldalphanew}{\boldsymbol{\alpha}}}%
%BeginExpansion
\boldsymbol{\alpha}%
%EndExpansion
}\;\right\vert \;\left(  \tau,%
%TCIMACRO{\TeXButton{boldalphanew}{\boldsymbol{\alpha}}}%
%BeginExpansion
\boldsymbol{\alpha}%
%EndExpansion
\right)  \in\mathcal{I}_{k}\left(  \partial K\right)  \right\}  \cup\left\{
\left.  B_{k}^{\operatorname{CR},F}\;\right\vert \;F\in\mathcal{F}\left(
K\right)  \right\}  $ are linearly independent on $\partial K$. Thus using the
biduality property of $\mathcal{J}_{k}\left(  \mathcal{T}\right)  $ we deduce
that $\left.  \psi_{\partial K}\right\vert _{\partial K}=1$. This and
Definition \ref{DefNCShapeFct} implies that $\left.  \phi_{\partial
K}\right\vert _{\partial K}=\left.  B_{k}^{\operatorname{CR},K}\right\vert
_{\partial K}$ and hence
\[
\phi_{\partial K}\perp\mathbb{P}_{k-1}\left(  F\right)  \qquad\forall
F\in\mathcal{F}\left(  K\right)
\]
by construction. Hence it follows that $J_{\tau,%
%TCIMACRO{\TeXButton{boldalphanew}{\boldsymbol{\alpha}}}%
%BeginExpansion
\boldsymbol{\alpha}%
%EndExpansion
}^{F}\left(  \phi_{\partial K}\right)  =J_{F}^{\operatorname{CR}}\left(
\phi_{\partial K}\right)  =0$ for all $F\in\mathcal{F}\left(  K\right)  $ and
all $\left(  \tau,%
%TCIMACRO{\TeXButton{boldalphanew}{\boldsymbol{\alpha}}}%
%BeginExpansion
\boldsymbol{\alpha}%
%EndExpansion
\right)  \in\mathcal{I}_{k}^{d-1}\left(  K\right)  $ with $\operatorname{mark}%
_{\mathcal{T}}\left(  \tau\right)  =F$. Therefore $\phi_{\partial K}$ can be
written as
\[
\phi_{\partial K}=\sum_{\left(  \tau,%
%TCIMACRO{\TeXButton{boldalphanew}{\boldsymbol{\alpha}}}%
%BeginExpansion
\boldsymbol{\alpha}%
%EndExpansion
\right)  \in\mathcal{R}_{k}\left(  \partial K\right)  }J_{\tau,%
%TCIMACRO{\TeXButton{boldalphanew}{\boldsymbol{\alpha}}}%
%BeginExpansion
\boldsymbol{\alpha}%
%EndExpansion
}\left(  1\right)  B_{\tau,%
%TCIMACRO{\TeXButton{boldalphanew}{\boldsymbol{\alpha}}}%
%BeginExpansion
\boldsymbol{\alpha}%
%EndExpansion
},
\]
where $\mathcal{R}_{k}\left(  \partial K\right)  :=\left\{  \left.  \left(
\tau,%
%TCIMACRO{\TeXButton{boldalphanew}{\boldsymbol{\alpha}}}%
%BeginExpansion
\boldsymbol{\alpha}%
%EndExpansion
\right)  \in\mathcal{I}_{k}\left(  \partial K\right)  \;\right\vert
\;\operatorname{mark}_{\mathcal{T}}\left(  \tau\right)  \in\mathcal{F}%
\setminus\mathcal{F}\left(  K\right)  \right\}  $. As a consequence of Lemma
\ref{LemLinInd} this yields that $\phi_{\partial K}\in\overset{\bullet}{S}%
_{k}\left(  \mathcal{T}\right)  $ and, in turn,%
\begin{equation}
\left.  \phi_{\partial K}\right\vert _{K}\left(  \mathbf{z}\right)
=0\qquad\mathbf{z}\in\mathcal{V}\left(  K\right)  .
\label{Eq:PhiPartialK 0 at vertex}%
\end{equation}
However the combination of \eqref{TechThingsa} and \eqref{Eq:Def phiPartialK}
gives $\left.  \phi_{\partial K}\right\vert _{\partial K}\left(
\mathbf{z}\right)  \neq0$, which contradicts (\ref{Eq:PhiPartialK 0 at vertex}%
). \newline\textbf{@2.} Set $\mathcal{T}_{\Omega}:=\left\{  \left.
K\in\mathcal{T}\;\right\vert \;K\cap\partial\Omega=\emptyset\right\}  $. By
assumption $\left\vert \mathcal{T}_{\Omega}\right\vert \geq1$. By taking some
$K\in\mathcal{T}_{\Omega}$ and repeating the previous construction the second
claim follows.%
%TCIMACRO{\TeXButton{End Proof}{\endproof}}%
%BeginExpansion
\endproof
%EndExpansion

\end{proof}

\begin{remark}
The reason why this proof does not apply for $k=1$ or $d=2$ is that then
\eqref{Eq:Def phiPartialK} and (\ref{Eq:PhiPartialK 0 at vertex}) do not
contradict each other as a consequence of \eqref{TechThingsa}.
\end{remark}

\appendix

\section{Proof of the determinant formula (\ref{defformula})\label{AppDetForm}%
}

In this section, we prove formula (\ref{defformula}).

\begin{lemma}
\label{LemQR}The determinant of the matrix $\mathbf{Q}_{d}\left(  s\right)  $
in (\ref{defQd}) is given by%
\begin{equation}
\det\mathbf{Q}_{d}\left(  s\right)  =\left(  -1\right)  ^{d+1}\left(
1+s\right)  ^{d-1}\left(  d-1-s\right)  \label{defformulalem}%
\end{equation}

\end{lemma}%

%TCIMACRO{\TeXButton{Proof}{\proof}}%
%BeginExpansion
\proof
%EndExpansion
We define the auxiliary matrix%
\begin{equation}
\mathbf{R}_{d}\left(  s\right)  :=\left[
\begin{array}
[c]{cccc}%
1 & 1 & \ldots & 1\\
1 & -s & \ddots & \vdots\\
\vdots & \ddots & \ddots & 1\\
1 & \ldots & 1 & -s
\end{array}
\right]  \label{defRd}%
\end{equation}
and claim%
\begin{equation}
\det\mathbf{R}_{d}\left(  s\right)  =\left(  -1-s\right)  ^{d-1}.
\label{detform2}%
\end{equation}
We prove (\ref{defformulalem}) and (\ref{detform2}) simultaneously by induction.

The Laplace expansion of the determinant of $\mathbf{Q}_{d}$ with respect to
the first row yields%
\[
\det\mathbf{Q}_{d}\left(  s\right)  =-s\det\mathbf{Q}_{d-1}\left(  s\right)
+\sum_{\ell=1}^{d-1}\left(  -1\right)  ^{\ell}\det\mathbf{R}_{d-1,\ell}\left(
s\right)  ,
\]
where the matrix $\mathbf{R}_{d-1,\ell}\left(  s\right)  $ arises by moving
the first row of $\mathbf{R}_{d-1}\left(  s\right)  $ to the $\ell$-th column
and shifting all rows $2,\ldots,\ell-1$ one row up. In other word, an $\left(
\ell-1\right)  $-fold interchange of appropriate rows of $\mathbf{R}%
_{d-1,\ell}\left(  s\right)  $ yields the matrix $\mathbf{R}_{d-1}\left(
s\right)  $. Then, well-known properties of determinants imply%
\[
\det\mathbf{Q}_{d}\left(  s\right)  =-s\det\mathbf{Q}_{d-1}\left(  s\right)
-\left(  d-1\right)  \det\mathbf{R}_{d-1}\left(  s\right)  .
\]
We employ the same reasoning for $\mathbf{R}_{d}\left(  s\right)  $ and obtain%
\[
\det\mathbf{R}_{d}\left(  s\right)  =\det\mathbf{Q}_{d-1}\left(  s\right)
+\sum_{\ell=1}^{d-1}\left(  -1\right)  ^{\ell}\det\mathbf{R}_{d-1,\ell}\left(
s\right)  =\det\mathbf{Q}_{d-1}\left(  s\right)  -\left(  d-1\right)
\det\mathbf{R}_{d-1}\left(  s\right)  .
\]
For the induction start we explicitly compute%
\[%
\begin{array}
[c]{ll}%
\det\mathbf{Q}_{1}\left(  s\right)  =-s, & \det\mathbf{R}_{1}\left(  s\right)
=1
\end{array}
\]
and for the induction step we assume that (\ref{defformulalem}) and
(\ref{detform2}) hold up to $d-1$. We then get%
\begin{align*}
\det\mathbf{Q}_{d}\left(  s\right)   &  =-s\det\mathbf{Q}_{d-1}\left(
s\right)  -\left(  d-1\right)  \det\mathbf{R}_{d-1,\ell}\left(  s\right) \\
&  =-s\left(  \left(  -1\right)  ^{d}(1+s)^{d-2}\left(  d-2-s\right)  \right)
-\left(  d-1\right)  \left(  -1-s\right)  ^{d-2}\\
&  =\left(  -1\right)  ^{d}\left(  1+s\right)  ^{d-2}\left(  s+1\right)
\left(  s-d+1\right)
\end{align*}
and this is (\ref{defformulalem}). In a similar way we get%
\begin{align*}
\det\mathbf{R}_{d}\left(  s\right)   &  =\det\mathbf{Q}_{d-1}\left(  s\right)
-\left(  d-1\right)  \det\mathbf{R}_{d-1}\left(  s\right) \\
&  =\left(  -1\right)  ^{d}(1+s)^{d-2}\left(  d-2-s\right)  -\left(
d-1\right)  \left(  -1-s\right)  ^{d-2}\\
&  =\left(  1+s\right)  ^{d-2}\left(  -1\right)  ^{d}\left(  -s-1\right)
\end{align*}
and this is (\ref{detform2}).%
%TCIMACRO{\TeXButton{End Proof}{\endproof}}%
%BeginExpansion
\endproof
%EndExpansion

\bibliographystyle{abbrv}
\bibliography{nlailu}

\end{document}